\documentclass[11pt]{amsart}
\usepackage{amssymb,amsmath,amsthm,amscd}
\usepackage{latexsym}
\usepackage[mathscr]{eucal}
\usepackage[english]{babel}
\usepackage{graphicx}

\newcommand{\g}{\mbox{${\mathfrak g}$}}

\newcommand{\kf}{\mbox{${\mathfrak k}$}}

\newcommand{\p}{\mbox{${\mathfrak p}$}}

\newcommand{\R}{\mbox{${\mathbb R}$}}
\newcommand{\Z}{\mbox{${\mathbb Z}$}}

\newtheorem{theorem}[subsection]{Theorem}
\newtheorem{lemma}[subsection]{Lemma}
\newtheorem{prop}[subsection]{Proposition}

\newtheorem{corollary}[subsection]{Corollary}

\newtheorem{example}[subsection]{Example}
\newtheorem{remark}[subsection]{Remark}

\newtheorem{notation}[subsection]{Notation}
\newenvironment{rmk}{\begin{remark} \em}{\end{remark}}

\newenvironment{ex}{\begin{example} \em}{\end{example}}
\newenvironment{notn}{\begin{notation} \em}{\end{notation}}

\def\qed{\vbox{\hrule
    \hbox{\vrule\hbox to 6pt{\vbox to 8pt{\vfil}\hfil}\vrule}\hrule}}

\numberwithin{equation}{section}

\begin{document}

\title{Classification of Superpotentials}
\author{A. Dancer}
\address{Jesus College, Oxford University, Oxford, OX1 3DW, United
Kingdom}
\email{dancer@maths.ox.ac.uk}
\author{M. Wang}
\address{Department of Mathematics and Statistics,
McMaster University, Hamilton, Ontario, L8S 4K1, Canada}
\email{wang@mcmaster.ca}
\thanks{The second author was partly supported by NSERC grant
No. OPG0009421.}

\date{revised \today}

\pagestyle{headings}

\begin{abstract}
We extend our previous classification \cite{DW4} of superpotentials of
``scalar curvature type" for the cohomogeneity one Ricci-flat equations.
We now consider the case not covered in \cite{DW4}, i.e., when some
weight vector of the superpotential lies outside (a scaled translate of)
the convex hull of the weight vectors associated with the scalar curvature
function of the principal orbit. In this situation we show that
either the isotropy representation has at most $3$ irreducible summands
or the first order subsystem associated to the superpotential is of the
same form as the Calabi-Yau condition for submersion type metrics
on complex line bundles over a Fano K\"ahler-Einstein product.
\end{abstract}

\maketitle


\setcounter{section}{-1}

\section{\bf Introduction}

In this paper we continue the study we began in \cite{DW4} of
superpotentials for the cohomogeneity one Einstein equations.
These equations are the ODE system obtained
as a reduction of the Einstein equations by requiring that
the Einstein manifold admits an isometric Lie group action whose
principal orbits $G/K$ have codimension one \cite{BB}, \cite{EW}.
As discussed in \cite{DW3}, these equations can be viewed as a
Hamiltonian system with constraint for a suitable Hamiltonian
$\sf H$, in which the potential term depends on the Einstein constant
and the scalar curvature of the principal orbit, and the kinetic
term is essentially the Wheeler-deWitt metric, which is of
Lorentz signature.

For any Hamiltonian system with Hamiltonian $\sf H$ and position
variable $q$, a {\em superpotential} is a {\em globally defined}
function $u$ on configuration space that satisfies the equation
\begin{equation} \label{HJ}
{\sf H}(q, du_q) = 0.
\end{equation}

{}From the classical physics viewpoint, $u$ is a $C^2$ (rather than a
viscosity) solution of a time-independent Hamilton-Jacobi equation.
The literature for implicitly defined first order partial differential
equations then suggests that such solutions are fairly rare. It is
therefore not unreasonable to expect in our case that one can
classify (at least under appropriate conditions) those principal orbits
where the associated cohomogeneity one Einstein equations admit a superpotential.

The existence of such a superpotential $u$ in our case leads
naturally to a subsystem of equations of half the dimension of the full
Einstein system. One way to see this is via generalised first integrals
which are linear in momenta, described in \cite{DW4}. Schematically,
the subsystem may be written as
$\dot{q} =  {\mathcal J} \nabla u$
where $\mathcal J$ is an endomorphism related to the kinetic term of
the Einstein Hamiltonian.

String theorists have exploited the superpotential idea in their
search for explicit metrics of special holonomy (see for example
\cite{CGLP1}, \cite{CGLP2},\cite{CGLP3}, \cite{BGGG} and references
in \cite{DW4}). The point here is that the subsystem defined by the
superpotential often (though not always) represents the condition that
the metric has special holonomy. Also, the subsystem can often be integrated explicitly.

In \cite{DW4}, \S 6,  we obtained classification results for superpotentials
of the cohomogeneity one Ricci-flat equations.  Besides assuming that
$G$ and $K$ are both compact, connected Lie groups such that the isotropy
representation of $G/K$ is multiplicity-free, we also mainly restricted our
attention to superpotentials which are of the same form as the scalar curvature
function of $G/K$, i.e., a finite sum with constant coefficients
of exponential terms. Almost all the known superpotentials are of this kind.

However, the above classification results were further subject to
the technical assumption that the extremal weights for the superpotential
did not lie in the null cone of the Wheeler-de Witt metric.  In \cite{DW4}
we gave some examples of superpotentials which do not satisfy this hypothesis.
These included several new examples which do not seem to be associated
to special holonomy.

In this paper, therefore, we attempt to solve the classification
problem without the non-null assumption on the extremal weights.

As in \cite{DW4}, we use techniques of convex geometry to analyse
the two polytopes naturally associated to the classification problem.
The first is (a rescaled translate of) the convex hull conv$(\mathcal W)$
of the weight vectors appearing in the scalar curvature function of
the principal orbit. The second is the convex hull conv$(\mathcal C)$
of the weight vectors in the superpotential.  In \cite{DW4} we showed
that the non-null assumption forces these polytopes to be equal, so we
could analyse the existence of superpotentials by looking at the
geometry of conv$(\mathcal W)$.

In the current paper, conv$(\mathcal C)$ may be strictly bigger than
conv$(\mathcal W)$ because of the existence of vertices outside
conv$(\mathcal W)$ but lying on the null cone of the Wheeler-de Witt metric.
Our strategy is to consider such a vertex $c$ and project conv$(\mathcal W)$
onto an affine hyperplane separating $c$ from conv$(\mathcal W)$.
We can now analyse the existence of superpotentials in terms of the
projected polytope.

The analysis becomes considerably more complicated because, whereas
in \cite{DW4} we could analyse the situation by looking at the vertices
and edges of conv$(\mathcal W)$, now, because we have projected onto a
subspace of one lower dimension, we have to consider the $2$-dimensional
faces of conv$(\mathcal W)$ also.

We find that in this situation the only polytopes ${\rm conv}(\mathcal W)$
arising from  principal orbits with more than three irreducible summands
in their isotropy representations  are  precisely those coming from
principal orbits which are circle bundles over a (homogeneous) Fano product.
In the latter case, the solutions of the subsystem defined by the superpotential
correspond to Calabi-Yau metrics, as discussed in \cite{DW4}.

After a review of basic material in \S 1, we state the main classification
theorem of the paper in \S2 and give an outline of the strategy
of the proof there.

\section{\bf Review and notation}

In this section we fix notation for the problem and review the set-up
of \cite{DW4}.

Let $G$ be a compact Lie group, $K \subset G$ be a closed subgroup,
and $M$ be a cohomogeneity one $G$-manifold of dimension $n+1$ with
principal orbit type $G/K$, which is assumed to be connected and
almost effective.
A $G$-invariant metric $\overline g$ on $M$ can be written in the form
$\overline{g} = \varepsilon dt^2 + g_t$
where $t$ is a coordinate transverse to the principal orbits,
$\varepsilon = \pm 1,$  and $g_t$ is a $1$-parameter family of
$G$-homogeneous Riemannian metrics on $G/K$. When $\varepsilon = 1$,
the metric ${\overline g}$ is Riemannian, and when $\varepsilon = -1$,
the metric $\overline g$ is spatially homogeneous Lorentzian, i.e.,
the principal orbits are space-like hypersurfaces.

We choose an $Ad(K)$-invariant decomposition $ \g = \kf \oplus \p  $
where $\g$ and $\kf$ are respectively the Lie algebras of $G$
and $K$, and $\p$ is identified with the isotropy representation
of $G/K$. Let
\begin{equation} \label{decomp}
\p =  \p_1 \oplus \cdots \oplus \p_r
\end{equation}
be a decomposition of $\p \approx T_{(K)} (G/K)$ into irreducible
real $K$-representations. We let $d_i$ be the real dimension of
$\p_i$, and $n = \sum_{i=1}^{r} d_i$ be the dimension of $G/K$
(so $\dim M =n+1$). We use $d$ for the vector of dimensions
$(d_1, \cdots, d_r)$. We shall assume that the isotropy
representation of $G/K$ is {\em multiplicity free}, i.e., all the
summands $\p_i$ in (\ref{decomp}) are distinct as $K$-representations.
In particular, if there is a trivial summand it must be 1-dimensional.

We use $q=(q_1, \cdots, q_r)$ to denote exponential coordinates
on the space of $G$-invariant metrics on $G/K$.
The Hamiltonian $\sf H$ for the cohomogeneity one Einstein
equations with principal orbit $G/K$ is now given by:
\[
{\sf H} = {\sf v}^{-1} J + \varepsilon {\sf v} \left( (n-1) \Lambda - {\sf S} \right),
\]
where $\Lambda$ is the Einstein constant, ${\sf v} = \frac{1}{2}e^{d \cdot q}$
is the relative volume and
\begin{equation} \label{KE}
J(p,p)= \frac{1}{n-1} \left( \sum_{i=1}^{r} p_i \right)^2
    -\sum_{i=1}^{r} \frac{p_i^2}{d_i},
\end{equation}
which has signature $(1,r-1)$. The scalar curvature
${\sf S}$ of $G/K$ above can be written as
\[
{\sf S} = \sum_{w \in {\mathcal W}} A_w e^{w \cdot q},
\]
where $A_w$ are {\em nonzero} constants and $\mathcal W$ is a finite
collection of vectors $w \in \Z^r \subset \R^r$. The set $\mathcal W$
depends only on $G/K$ and its elements will be referred to as
{\em weight vectors}. These are of three types

\medskip
(i) type I: one entry of $w$ is $-1$, the others are zero,

(ii) type II: one entry is 1, two are -1, the rest are zero,

(iii) type III: one entry is 1, one is -2, the rest are zero.

\begin{notn} \label{shorthand}
As in \cite{DW4} we use  $(-1^i, -1^j, 1^k)$ to denote the
type II vector $w \in \mathcal W \subset \R^r$ with $-1$ in places $i$ and $j$,
and $1$ in place $k$. Similarly,  $(-2^i, 1^j)$ will
denote the type III vector with $-2$ in place $i$ and $1$ in place $j$, and
$(-1^i)$ the type I vector with $-1$ in place $i$.
\end{notn}

\begin{rmk} \label{facts}
We collect below various useful facts from \cite{DW4} and \cite{WZ1}.
Also, we shall use standard terminology from convex geometry, as given,
e.g., in \cite{Zi}. In particular, a ``face" is not necessarily
$2$-dimensional. However, a vertex and an edge are respectively
zero and one-dimensional. The convex hull of a set $X$ in $\R^r$
will be denoted by conv$(X)$.

(a) For a type I vector $w$, the
coefficient $A_w > 0$ while for type II and type III vectors,
$A_w < 0$.

(b) The type I vector with $-1$ in the $i$th position is absent
from $\mathcal W$ iff
the corresponding summand $\p_i$ is an abelian subalgebra which
satisfies $[\kf, \p_i] = 0$ and $[\p_i, \p_j] \subset \p_j$ for
all $j \neq i$. If the isotropy group $K$  is connected, these
last conditions imply that $\p_i$ is $1$-dimensional, and the
$\p_j, j \neq i,$ are irreducible representations of the (compact)
analytic group whose Lie algebra is $\kf \oplus \p_i$.

(c) If $(1^i, -1^j, -1^k)$ occurs in $\mathcal W$ then its permutations
$(-1^i, 1^j, -1^k)$ and
$(-1^i, -1^j, 1^k)$ do also.

(d) If $\dim \p_i = 1$ then no type III vector with $-2$ in place $i$ is
 present in $\mathcal W$.
If in addition $K$ is connected, then no type II
vector with nonzero entry in place $i$ is present.

(e) If $I$ is a subset of $\{1, \cdots, r \},$ then
each of the equations $\sum_{i \in I} x_i = 1$ and $\sum_{i \in I} x_i =-2$
defines a face (possibly empty) of conv$(\mathcal W)$.
In particular, all type III vectors in $\mathcal W$ are vertices
and $(-1^i, -1^k, 1^j) \in {\mathcal W}$ is a vertex unless
both $(-2^i, 1^j)$ and $(-2^k, 1^j)$ lie in $\mathcal W$.

(f) For $v, w \in \mathcal{W}$ (or indeed for any $v,w$ such that $\sum v_i$
or $\sum w_i=-1$), we have
\begin{equation} \label{Jform}
J(v +d, w+ d) = 1 -\sum_{i=1}^{r} \frac{v_i w_i}{d_i}.
\end{equation}

\end{rmk}

\bigskip
For the remainder of the paper, we shall work in the Ricci-flat
Riemannian case, that is, we take $\varepsilon=1$ and  $\Lambda=0$.
As in \cite{DW4}, any argument that does not use the sign of $A_w$
would be valid in the Lorentzian case. We shall also assume that
conv$({\mathcal W})$ is $r-1$ dimensional. This is certainly the case
if $G$ is semisimple, as $\mathcal W$ spans ${\mathbb R}^r$
(see the proof of Theorem 3.11 in \cite{DW3}).

The superpotential equation (\ref{HJ}) now becomes
\begin{equation} \label{seqn}
       J(\nabla u , \nabla u) = e^{d \cdot q} \ {\sf S},
\end{equation}
where $\nabla$ denotes the Euclidean gradient in $\R^r$.
As in \cite{DW4} we shall look for solutions to Eq.(\ref{seqn}) of the form
\begin{equation} \label{potform}
u = \sum_{\bar{c} \in \mathcal C} F_{\bar{c}} \ e^{\bar{c} \cdot q}
\end{equation}
where $\mathcal C$ is a finite set in $\R^r$, and the $F_{\bar{c}}$
are nonzero constants. Now Eq.(\ref{seqn}) reduces to, for each
$\xi \in \R^r$,
\begin{equation} \label{eqnF}
\sum_{\bar{a}+\bar{c} = \xi} J(\bar{a},\bar{c}) \ F_{\bar{a}} F_{\bar{c}} = \left\{
  \begin{array}{ll}
    A_w  &   \mbox{if} \ \xi = d + w \  \mbox{for some} \ w \in \mathcal{W} \\
     0    &  \mbox{if} \ \xi \notin d + \mathcal W.
   \end{array} \right.
\end{equation}

\medskip
We shall assume henceforth that $r \geq 2$ since the superpotential equation
always has a solution in the $r=1$ case, as was noted in \cite{DW4}, and $J$
is of Lorentz signature only when $r \geq 2$.
The following facts were deduced in \cite{DW4} from Eq.(\ref{eqnF}).

\begin{prop} \label{convinc}
  {\rm conv}$( \frac{1}{2}(d + {\mathcal W})) \subset$ {\rm conv}$(\mathcal{C})$.
\end{prop}

\begin{proof}
\ If $w \in \mathcal{W}$, then Eq.(\ref{eqnF})
implies that $d + w = \bar{a}+\bar{c}$ for some $\bar{a}, \bar{c} \in
\mathcal C$, and hence that $\frac{1}{2} (d + w) =
\frac{1}{2}(\bar{a}+\bar{c}) \in$ conv$(\mathcal C)$.  \ \ \qed
\end{proof}

\begin{prop} \label{E}
If $\bar{a},\bar{c} \in \mathcal C$ and $\bar{a}+\bar{c}$ cannot
be written as the sum of two non-orthogonal
elements of $\mathcal C$ distinct from
$\bar{a},\bar{c}$ then either $J(\bar{a},\bar{c})=0$ or
$\bar{a} + \bar{c} \in d + \mathcal W.$

In particular, if $\bar{c}$ is a vertex of $\mathcal C$, then either
$J(\bar{c},\bar{c})=0,$ or
$2\bar{c} = d + w$ for some $w \in \mathcal W$ and
 $J(\bar{c},\bar{c}) \ F_{\bar{c}}^2 = A_w$.
In the latter case, $J(d +w, d+w)$ has the same sign as $A_w$ so is
$>0$ if $w$ is type I and $<0$ if $w$ is type II or III. \;
\end{prop}

As mentioned in the Introduction, for the classification in \cite{DW4}
we made the assumption that all vertices $\bar{c}$ of $\mathcal C$
are  non-null. Under this assumption, the second assertion of
Prop \ref{E} implies that all vertices of $\mathcal C$ lie in
$\frac{1}{2}(d + \mathcal{W})$. Hence conv$(\mathcal{C})$ is contained in
conv$(\frac{1}{2} ( d + \mathcal{W}))$, and by Prop \ref{convinc}
they are equal.  This meant that in \cite{DW4}, subject to the
non-null assumption, we could study the existence of a superpotential
in terms of the convex geometry of $\mathcal W$.

The aim of the current paper is to drop this assumption.
We still have
\[
{\rm conv}( \frac{1}{2}(d + {\mathcal W})) \subset {\rm conv}(\mathcal{C}),
\]
but can no longer deduce that these sets are equal. The problem is
that a vertex $\bar{c}$ of conv$(\mathcal C)$ may lie outside
{\rm conv}$( \frac{1}{2}(d + {\mathcal W}))$ if it is null.

In fact, it is clear from the above discussion that conv$(\frac{1}{2}(d +
{\mathcal W}))$ is strictly contained in conv$(\mathcal C)$ {\em if and only if}
$\mathcal C$ has a null vertex. For if $\bar{c}$ is a null vertex
of $\mathcal C$ and $2 \bar{c} = d+ w$ for some $w \in \mathcal W$, then
Eq.(\ref{eqnF}) fails for $\xi = d+w$.

\smallskip

We conclude this section by proving an analogue of Proposition 2.5
in \cite{DW4}. The arguments below using Prop \ref{E} are ones which will
recur throughout this paper. Henceforth when we use the term ``orthogonal''
we mean orthogonal with respect to $J$ unless otherwise stated.

\begin{theorem} \label{hplane}
$\mathcal C$ lies in the hyperplane $\{ \bar{x} : \sum \bar{x}_i= \frac{1}{2}(n-1) \}$
$($possibly after subtracting a constant from the superpotential$)$.
\end{theorem}
\begin{proof}
We can assume $0 \notin {\mathcal C}$ by subtracting a constant
from the superpotential. We shall also  use repeatedly below the fact
that as $J$ has signature $(1, r-1)$ there are no null planes,
only null lines.

Denote by $H_{\lambda}$ the hyperplane $\sum \bar{x}_i = \lambda$, so
$\frac{1}{2}(d+ {\mathcal W})$ lies in $H_{\frac{1}{2}(n-1)}$.
Suppose there exist elements of $\mathcal C$ with $\sum \bar{x}_i > \frac{1}{2}(n-1)$.
Let $\lambda_{\rm max}$ denote the greatest value of $\sum \bar{x}_i$ over $\mathcal C$.
If $\tilde{a} \tilde{c}$ is an edge  of
conv$({\mathcal C}) \cap H_{\lambda_{\rm max}}$, then Prop \ref{E} shows
that $\tilde{a}, \tilde{c}$ are null, and that $\tilde{c}$ is orthogonal to the
element of $\mathcal C$ closest to it on the edge. Hence $\tilde{c}$
is orthogonal to the whole edge. Now $J$ is totally null on
Span$\{\tilde{a}, \tilde{c} \}$, so since there are no null planes,
$\tilde{a}, \tilde{c}$ are proportional, which is impossible as they
are both in $H_{\lambda_{\rm max}}$. So ${\mathcal C} \cap H_{\lambda_{\rm max}}$
is a single point $\tilde{c}_{\rm max}$, which is null.

Next we claim that all elements of $\mathcal C$ lying in the
half-space $\sum \bar{x}_i > \frac{1}{2}(n-1)$ must be multiples of
$\tilde{c}_{\rm max}$. If not, let $\lambda_*$ be the greatest
value such that there is an element of $\mathcal C$, not
proportional to $\tilde{c}_{\rm max}$, in $H_{\lambda_*}$.
Let $\tilde{a}$  be a vertex of conv$({\mathcal C}) \cap H_{\lambda_*}$,
not proportional to $\tilde{c}_{\rm max}$. Now, by Prop \ref{E},
$J(\tilde{a}, \tilde{c}_{\rm max})=0$, and so $\tilde{a}$ is not null.
Since $\lambda_* > \frac{1}{2}(n-1)$, we see $\tilde{a} + \tilde{a}$
must be written in another way as a sum of two non-orthogonal elements of
$\mathcal C$. This sum must be of the form
$\mu \tilde{c}_{\rm max} + \tilde{f}$.  But $\tilde{c}_{\rm max}$
is orthogonal to $\tilde{a}$ and to itself, hence to $\tilde{f}$,
a contradiction establishing our claim.

Similarly, all elements of $\mathcal C$ lying in $\sum \bar{x}_i < \frac{1}{2}(n-1)$
are multiples of an element $\tilde{c}_{\rm min}$, should they occur. (Note that
$J$ is negative definite on $H_0$ and we have assumed $0 \notin \mathcal C$
so $\lambda_{\rm min} \neq 0$.)

We denote the sets of elements lying in these open half-spaces by
${\mathcal C}_{+}$ and ${\mathcal C}_{-}$ respectively.
Note that, when non-empty, ${\mathcal C}_{+}$ and ${\mathcal C}_{-}$
are orthogonal to all elements of ${\mathcal C} \cap H_{\frac{1}{2}(n-1)}$.
(For if $\tilde{a} \in {\mathcal C} \cap H_{\frac{1}{2}(n-1)}$ then
$\tilde{a} + \tilde{c}_{\rm max}$ cannot be written in another way as a sum
of two non-orthogonal elements of $\mathcal C$.) In particular,
if $\tilde{c}_{\rm max}$ and $\tilde{c}_{\rm min}$ are orthogonal,
then $\tilde{c}_{\rm max}$ is orthogonal to all of conv$(\mathcal C)$,
which is $r$-dimensional by assumption. So $\tilde{c}_{\rm max}$ is zero,
a contradiction. The same argument implies that
${\mathcal C}_{+}$ and  ${\mathcal C}_{-}$ are both non-empty.

Let $\nu \tilde{c}_{\rm min}$ and $\mu \tilde{c}_{\rm max}$ be
respectively the elements of ${\mathcal C}_{-}$ and ${\mathcal C}_{+}$
closest to $H_{\frac{1}{2}(n-1)}$. Suppose that $\tilde{c}_{\rm max} +
\nu \tilde{c}_{\rm min} = \tilde{c}^{(1)} + \tilde{c}^{(2)}$ with
$\tilde{c}^{(i)} \in \mathcal C$ and $J(\tilde{c}^{(1)},
\tilde{c}^{(2)}) \neq 0$. Non-orthogonality
means the $\tilde{c}^{(i)}$ cannot belong to
the same side of $H_{\frac{1}{2}(n-1)}$ and by the choice of $\nu$,
they cannot belong to opposite sides of $H_{\frac{1}{2}(n-1)}$.  Both
therefore lie in $H_{\frac{1}{2}(n-1)}$. But by the previous
paragraph, $J(\tilde{c}_{\rm max} + \nu \tilde{c}_{\rm min}, \
\tilde{c}^{(1)} + \tilde{c}^{(2)}) = 0.$ This means that
$\tilde{c}_{\rm max} + \nu \tilde{c}_{\rm min}$ is null, which
contradicts $J(\tilde{c}_{\rm max}, \ \tilde{c}_{\rm min}) \neq 0$.
Hence $\tilde{c}_{\rm max} + \nu \tilde{c}_{\rm min}$ lies in $d+
{\mathcal W} \subset H_{n-1}$. Applying the same argument to
$\tilde{c}_{\rm min} + \mu \tilde{c}_{\rm max},$ we find that in fact
$\mu = \nu =1$, i.e., ${\mathcal C}_{+} = \{ \tilde{c}_{\rm max} \}$
and ${\mathcal C}_{-} = \{ \tilde{c}_{\rm min} \}$.

Now $\mathcal C \cap H_{\frac{1}{2}(n-1)}$ (and hence its convex hull)
is contained in the hyperplanes $\tilde{c}_{\rm max}^{\perp}$,
$\tilde{c}_{\rm min}^{\perp}$ in $H_{\frac{1}{2}(n-1)}$. These
hyperplanes are distinct as $\tilde{c}_{\rm max}$ is orthogonal to
itself but not to $\tilde{c}_{\rm min}$. Hence $d + {\mathcal W}
\subset ({\mathcal C} + {\mathcal C}) \cap H_{n-1}$ is contained in
the union of the point $\tilde{c}_{\rm max} + \tilde{c}_{\rm min}$ and
the codimension 2 subspace $\tilde{c}_{\rm max}^{\perp} \cap
\tilde{c}_{\rm min}^{\perp}$ of $H_{n-1}$. So
conv$(d + {\mathcal W})$ is contained in a
codimension 1 subspace of $H_{n-1}$, contradicting our assumption that
$\dim {\rm conv}(d+{\mathcal W}) =r-1$.
\end{proof}

\smallskip

\begin{rmk}\label{barring}
A notational difficulty arises from the fact that, as seen above,
points of $\mathcal C$ are on the same footing as points in
$\frac{1}{2}(d + {\mathcal W})$ rather than points of $\mathcal W$.
Accordingly, we shall use letters $c,u,v,...$ to denote elements
of the hyperplane $\sum u_i =-1$ (such as elements of $\mathcal W$),
and $\bar{c}, \bar{u}, \bar{v},...$ to denote the {\em associated elements}
$\frac{1}{2}(d+c), \frac{1}{2}(d+u), \frac{1}{2}(d+v), \cdots$ of the
hyperplane $\sum \bar{u}_i = \frac{1}{2}(n-1)$ (such as elements of
$\mathcal C$ or of $\frac{1}{2}(d + {\mathcal W}))$.

Note that for any convex or indeed affine sum $\sum \lambda_j \xi^{(j)} $
of vectors $\xi^{(j)}$ in $\R^r$, we have
\[
\overline{ \sum_j \lambda_j \xi^{(j)}} = \sum_j \lambda_j \ \overline{\xi^{(j)}}.
\]
\end{rmk}

Since we now know that the set ${\mathcal C}$, like
$\frac{1}{2}(d + {\mathcal W})$, lies in
$H_{\frac{1}{2}(n-1)}:= \{ \bar{x} : \sum \bar{x}_i = \frac{1}{2}(n-1) \}$,
we will adopt the convention, as in the last paragraph, that
when we refer to hyperplanes such as ${\bar c}^{\perp}$ in the rest
of the paper, we mean ``affine hyperplanes in $H_{\frac{1}{2}(n-1)}$''.

\section{\bf The classification theorem and the strategy of its proof}

We can now state the main theorem of the paper.

\begin{theorem} \label{classthm}
Let $G$ be a compact connected Lie group and $K$ a closed connected
subgroup such that the isotropy representation of $G/K$ is the direct
sum of $r$ pairwise inequivalent $\R$-irreducible summands.
Assume that  $\dim {\rm conv}(\mathcal W) = r-1$, where $\mathcal W$
is the set of weights of the scalar curvature function of $G/K$
$($cf  \S 1$)$. $($This holds, for example, if $G$ is semisimple.$)$

If the cohomogeneity one Ricci-flat equations with $G/K$ as
principal orbit admit a superpotential of form $($\ref{potform}$)$
where $\mathcal C$ contains a $J$-null vertex, then we are in one
of the following situations $($up to permutations of the irreducible
summands$)$:
\begin{enumerate}
\item[$($i$)$] ${\mathcal W} =\{ (-1)^i, (1^1, -2^i): 2 \leq i \leq r \}$,
   $d_1=1$, ${\mathcal C} = \frac{1}{2}(d+\{(-1^1), (1^1, -2^i): 2 \leq i \leq r \})$
   and $r \geq 2$;
\item[$($ii$)$] $r \leq 3$.
\end{enumerate}
\end{theorem}

\begin{rmk}
As mentioned before, the situation where $\mathcal C$ has no null vertex
was analysed in \cite{DW4}. Hence, except for the $r \leq 3 $ case,
Theorem \ref{classthm} completes the classification of superpotentials
of scalar curvature type subject to the above assumptions on $G$ and $K$.
\end{rmk}

\begin{rmk} \label{realization}
The first case of the above theorem is realized by certain circle bundles
over a product of $r-1$ Fano (homogeneous) K\"ahler-Einstein manifolds
(cf. Example 8.1 in \cite{DW4}, and \cite{BB}, \cite{WW}, \cite{CGLP3}), and
the subsystem of the Ricci-flat equations singled out by the superpotential
in these examples corresponds to the Calabi-Yau condition.
For more on the $r=2$ case, see the concluding remarks in \S 10.
\end{rmk}

\begin{rmk} \label{disconnK}
Theorem \ref{classthm} remains true if we replace the connectedness of $G$
and $K$ by the connectedness of $G/K$ and the extra condition on the isotropy
representation given by the second statement in Remark \ref{facts}(d), i.e.,
if $\p_i$ is an irreducible summand of dimension $1$ in the isotropy
representation of $G/K$, then $[\p_i, \p_j] \subset \kf \oplus \p_j$ for
all $j \neq i$.

This weaker property does hold in practice. For example, the exceptional
Aloff-Wallach space $N_{1, 1}$ can be written as
$(SU(3) \times \Gamma)/(U_{1, 1} \cdot \Delta \Gamma)$, where $U_{1,1}$
is the set of diagonal matrices of the form
diag$(\exp(i\theta), \exp(i\theta), \exp(-2i\theta))$ and $\Gamma$ is
the dihedral group with generators
$$ \left( \begin{array}{rrr}
            0  & 1  &  0  \\
            -1 & 0  & 0  \\
            0 & 0 &  1
            \end{array} \right)   \ \ \ \
    \left(\begin{array}{ccc}
           e^{2 \pi i/3} & 0  & 0  \\
           0 & e^{-2 \pi i/3} & 0  \\
           0  &   0   &  1
           \end{array} \right).  $$
\end{rmk}

\medskip

In order to prove Theorem \ref{classthm} we have to analyse the situation
when there is a null vertex $\bar{c} \in \mathcal C$. As discussed in \S 1,
${\rm conv}({\mathcal C})$ now strictly includes
conv$(\frac{1}{2}(d + {\mathcal W}))$ as $\bar{c}$ is not in
${\rm conv}(\frac{1}{2}(d + {\mathcal W}))$. Our strategy is to take an
affine hyperplane $H$ separating $\bar{c}$ from
conv$(\frac{1}{2}(d + {\mathcal W}))$, and consider the projection
$\Delta^{\bar{c}}$ of conv$(\frac{1}{2}(d + {\mathcal W}))$ onto
$H$ from $\bar{c}$.

Roughly speaking, whereas in \cite{DW4} we could analyse the situation
by looking at the vertices and edges of conv$(\frac{1}{2}(d + {\mathcal
  W}))$, now, because we have projected onto a subspace of one lower
dimension, we have to consider the 2-dimensional faces of
conv$(\frac{1}{2}(d + {\mathcal W}))$ also.

This is a natural method of dealing with the situation of a point
outside a convex polytope. It has some relation to the notion of
``lit set'' introduced in a quite different context by
Ginzburg-Guillemin-Karshon \cite{GGK}.

The analysis in the next section will show that the vertices of the
projected polytope $\Delta$ can be divided into three types (Theorem
\ref{Vertex}). We label these types (1A), (1B) and (2). Roughly, these
correspond to vertices orthogonal to $\bar{c}$, vertices $\bar{\xi}$
such that the line through $\bar{c}$ and $\bar{\xi}$ meets
conv$(\frac{1}{2}(d + {\mathcal W}))$ at a vertex, and vertices
$\bar{\xi}$ such that this line meets conv$(\frac{1}{2}(d + {\mathcal
  W}))$ in an edge.

In the remainder of the paper we shall
gradually narrow down the possibilities for each type. In \S 3 we begin
a classification of type (2) vertices. In \S 4 we are able to deduce that
conv$(\frac{1}{2}(d + {\mathcal W}))$ lies in the half space $J(\bar{c}, \cdot)
\geq 0$. We are able to deduce an orthogonality result for
vectors on edges in conv$(\frac{1}{2}(d + {\mathcal W})) \cap \bar{c}^\perp$.
This is analogous to the key result Theorem 3.5 of \cite{DW4} that held
(in the more restrictive situation of that paper) for general
 edges in conv$(\frac{1}{2}(d + {\mathcal W}))$. In \S 5 we exploit this
result and some estimates to classify the possible configurations of
(1A) vertices (i.e. vertices in $\bar{c}^{\perp}$), see Theorem \ref{class-1a}.

In \S 6 we attack the (1B) vertices, exploiting the fact that
adjacent (1B) vertices give rise to a 2-dimensional face of
conv$(\frac{1}{2}(d + {\mathcal W}))$. This is the most laborious part
of the paper, as it involves a case-by-case analysis of such faces.
We show that adjacent (1B) vertices can arise only in a very small number of
situations (Theorem \ref{1B-initial}). In \S 7 we exploit the
listing of $2$-dim faces to show that there is at most one type (2)
vertex, except in two special situations (Theorem \ref{two2}).
In \S 8 and 9, we eliminate more possibilities for adjacent (1B)
and type (2) vertices. We find that if $r \geq 4$ then we are either
in case (i) of the Theorem or there are no type (2) vertices and
no adjacent type (1B) vertices.
Using the results of \S 4, in the latter case we find that
all vertices are (1A) except for a single (1B). Building on the
results of \S 5 for (1A) vertices, we are able to rule out this
situation in \S 10, see Theorem \ref{unique1B} and Corollary \ref{finish}.

\section{\bf Projection onto a hyperplane}

We first present some results about null vectors in $H_{\frac{1}{2}(n-1)}$.
\medskip
\begin{rmk} \label{tangent}
{}From Eq.(\ref{Jform}), the set of null vectors  in the hyperplane
$H_{\frac{1}{2}(n-1)}$ form an ellipsoid $\{ \sum x_i^2 / d_i =1 \}$.
If $\bar{c}$ is null,  then the hyperplane $\bar{c}^{\perp}$
in $H_{\frac{1}{2}(n-1)}$ is the tangent space to this ellipsoid.
So any element $\bar{x} \neq \bar{c}$ of $\bar{c}^{\perp}$
satisfies $J(\bar{x},\bar{x}) < 0$.
\end{rmk}

\begin{lemma} \label{Jpos}
Let $x, y$ satisfy $\sum x_i = \sum y_i =-1$.

Suppose that $J(\bar{x},\bar{x})$ and $J(\bar{y},\bar{y}) \geq 0$.
Then $J(\bar{x},\bar{y}) \geq 0$, with equality iff $\bar{x}$ is null
and $\bar{x}=\bar{y}$. In particular, if $\bar{x}, \bar{y}$
are distinct null vectors then $J(\bar{x},\bar{y}) > 0$.
\end{lemma}
\begin{proof}
This follows from Eq.(\ref{Jform}) and Cauchy-Schwartz.
\end{proof}

\begin{prop} \label{sep-hypl}
Let $H= \{ \bar{x} : h(\bar{x}) = \lambda \}$ be an affine hyperplane,
where $h$ is a linear functional such that
${\rm conv}(\frac{1}{2}(d + {\mathcal W}))$ lies in the open half-space
$\{ \bar{x} : h(\bar{x}) < \lambda \}$. Then there is at most one element
of $\mathcal C$ in the complementary open half-space
$\{ \bar{x} : h(\bar{x}) > \lambda \}$. Such an element is a null
vertex of ${\rm conv}(\mathcal C)$. Hence any element of
$\mathcal C$ outside ${\rm conv}(\frac{1}{2}(d + {\mathcal W}))$ is
a null vertex of $\mathcal C$.
\end{prop}

\begin{proof}
  Suppose the points of $\mathcal C$ with $h(\bar{x}) > \lambda$ are
  $\bar{c}^{(1)}, \cdots, \bar{c}^{(m)}$ with $m>1$.  Our result is stable
  with respect to sufficiently small perturbations of $H$, so we can
  assume that $h(\bar{c}^{(1)}) > h(\bar{c}^{(2)}) \geq h(\bar{c}^{(3)}), \cdots,
  h(\bar{c}^{(m)})$.

  Now $\bar{c}^{(1)} + \bar{c}^{(1)}$ and $\bar{c}^{(1)} + \bar{c}^{(2)}$
  cannot be written in any other way as the sum of two elements of
  $\mathcal C$. Hence, by Prop. \ref{E}, $\bar{c}^{(1)}$ is null and
  $J(\bar{c}^{(1)}, \  \bar{c}^{(2)})=0$.  The only other way
  $\bar{c}^{(2)} + \bar{c}^{(2)}$ can be
  written is as $\bar{c}^{(1)} + \bar{c}$ for some $\bar{c} \in \mathcal
  C$.  But then $\bar{c} = 2 \bar{c}^{(2)} - \bar{c}^{(1)}$, so
  $J(\bar{c}^{(1)}, \bar{c})=0,$ and such sums will not contribute.
  Hence $J(\bar{c}^{(2)}, \bar{c}^{(2)})=0$, contradicting Lemma \ref{Jpos}.
\end{proof}

\begin{corollary} \label{meet}
For distinct elements $\bar{c}, \bar{a}$ of $\mathcal C$, the line
segment $\bar{c} \bar{a}$ meets ${\rm conv}(\frac{1}{2}(d + {\mathcal W})).$
\end{corollary}

This gives us some control over the extent to which conv$(\mathcal C)$
can be bigger than the set ${\rm conv}(\frac{1}{2}(d + {\mathcal W}))$.

\begin{lemma} \label{III44} Let $A \subset H_{\frac{1}{2}(n-1)}$ be
an affine subspace such that
  $ A \cap \ {\rm conv}(\frac{1}{2}(d + {\mathcal W})) $ is a face
  of ${\rm conv}(\frac{1}{2}(d + {\mathcal W}))$. Suppose there
  exists $\bar{c} \in {\mathcal C} \cap A$ with $\bar{c} \notin {\rm
    conv}(\frac{1}{2}(d + {\mathcal W}))$.

  Let $\bar{x} \in A$.
  If $ \bar{x}=\frac{1}{2}(\bar{a} + \bar{a}^{\prime}) $ with $\bar{a},
  \bar{a}^{\prime} \in \mathcal C$, then in fact $\bar{a},
  \bar{a}^{\prime} \in A$.
\end{lemma}

\begin{proof}
If $\bar{a}$ or $\bar{a}^{\prime}$ equals $\bar{c}$ this is clear.

We know by Cor \ref{meet} that if $\bar{a}, \bar{a}^{\prime} \neq \bar{c}$
then the segments $\bar{a} \bar{c}$, $\bar{a}^{\prime} \bar{c}^{\prime}$ meet
${\rm conv}(\frac{1}{2}(d + {\mathcal W}))$. So there exist $0 < s,t \leq 1$
with $t\bar{a} + (1-t)\bar{c}$ and $s \bar{a}^{\prime} + (1-s)\bar{c}$
in ${\rm conv}(\frac{1}{2}(d + {\mathcal W}))$. Hence
\[
\left( \frac{2st}{s+t} \right) \bar{x} + \left(1 - \frac{2st}{s+t} \right)\bar{c} =
\frac{s}{s+t} \left(t\bar{a} + (1-t)\bar{c} \right) + \frac{t}{s+t} \left(s \bar{a}^{\prime} +
(1-s) \bar{c} \right) \in {\rm conv}(\frac{1}{2}(d + {\mathcal W})).
\]
As it is an affine combination of $\bar{x},\bar{c}$ this point also lies in
$A$, so it lies in $A \cap {\rm conv}(\frac{1}{2}(d + {\mathcal W}))$.
Also, it is a convex linear combination of the points $t\bar{a} +
(1-t)\bar{c}$ and $s \bar{a}^{\prime} + (1-s) \bar{c}$ of
${\rm conv}(\frac{1}{2}(d + {\mathcal W}))$.  Hence, by our face assumption,
both these points lie in $A$, so $\bar{a}, \bar{a}^{\prime}$ lie in $A$.
\end{proof}

\begin{rmk} The above lemma will be very useful because
it means that in all our later calculations using Prop \ref{E}
for a face defined by an affine
subspace $A$, we need only consider elements of $\mathcal C$
lying in  $A$.
\end{rmk}

\begin{prop} \label{edgeCC}
Let $vw$ be an edge of conv$({\mathcal W})$ and suppose $\bar{v}, \bar{w}
\in \mathcal C$.

$($i$)$ If there are no points of $\mathcal W$ in the interior of $vw$, then
$J(\bar{v}, \bar{w})=0$.

$($ii$)$ If $u=\frac{1}{2}(v+w)$ is the unique point of $\mathcal W$
in the interior of $vw$, $J(\bar{v}, \bar{w})>0$, and $u$ is type II
or III, then $F_{\bar{v}}, F_{\bar{w}}$ are of opposite signs.
\end{prop}

\begin{proof}
Part(i) is a generalization of Theorem 3.5 in \cite{DW4} and we will
be able to apply the proof of that result after the following argument.
Let the edge $\bar{v} \bar{w}$ of conv$(\frac{1}{2}( d + {\mathcal W}))$
be defined by equations $\langle \bar{x}, u^{(i)} \rangle =
\lambda_i \; : \; i \in {\mathcal I}$ where $\langle \bar{x}, u^{(i)} \rangle \leq \lambda_i$
for $i \in {\mathcal I}$ and $\bar{x} \in$ conv$(\frac{1}{2}( d + {\mathcal W}))$.
(In the above, $\langle \ \ , \ \ \rangle$ is the Euclidean inner product in $\R^r$.)
Note that Span $\{ u^{(i)} : i \in {\mathcal I} \}$ is the
$\langle \ , \ \rangle$-orthogonal complement
of the direction of the edge.

Let $H$ be a hyperplane whose intersection with
conv$(\frac{1}{2}(d+{\mathcal W}))$ is the edge $\bar{v}\bar{w}$.
We can take $H$ to be defined by the equation
$\langle \bar{x}, \sum_{i \in {\mathcal I}} b_i u^{(i)} \rangle = \sum_{i \in {\mathcal I}} b_i \lambda_i$
where $b_i$ are arbitrary positive numbers summing to 1.

If $\bar{a}, \bar{a}^{\prime}$ are elements of $\mathcal C$ whose midpoint
lies in $\bar{v} \bar{w}$, then either they are both in $H$ or one of
them, $\bar{a}$ say, is on the opposite side of $H$ from
conv$(\frac{1}{2}(d+{\mathcal W}))$.  In the latter case $\bar{a}$ is
null and the only element of $\mathcal C$ on this side of $H$ so, by
Prop \ref{E}, is $J$-orthogonal to $\bar{v}, \bar{w}$. Hence, as
$\frac{1}{2}(\bar{a} + \bar{a}^{\prime})$ is an affine combination of
$\bar{v}, \bar{w}$, we see that $J(\bar{a}, \bar{a}^{\prime})=0$,
and so such sums do not contribute in Eq.(\ref{eqnF}). We may
therefore assume that $\bar{a}, \bar{a}^{\prime}$ are in $H$.
But as this is true for all $H$ of the above form, the only sums
that will contribute are those where $\bar{a}, \bar{a}^{\prime}$
are collinear with $\bar{v} \bar{w}$.

Now if $\bar{a}$, say, lies outside the line segment
$\bar{v} \bar{w}$, then it is null and $J$-orthogonal to
$\bar{v}$ or $\bar{w}$, and hence to the whole line.
So the only sums which contribute are those where $\bar{a},\bar{a}^{\prime}$
lie on the line segment $\bar{v} \bar{w}$. Now the proof of Theorem 3.5
in  \cite{DW4} gives (i).

Turning to (ii), note first that the above arguments and
Prop \ref{E} give (ii) immediately if no interior points of
the edge $\bar{v}\bar{w}$ lie in $\mathcal C$. If there are $m$ interior
points in $\mathcal C$, we again proceed as in the proof of Theorem 3.5
in \cite{DW4} and use the notation there. We may assume that Lemma 3.2
(and hence Cor 3.3 and Lemma 3.4) of \cite{DW4}
still holds; for the only issue is the
statement for $\lambda_{m+1}$, but if $c^{(0)} + c^{(\lambda_{m+1})}$ cannot be
written as $c^{(\lambda_j)} + c^{(\lambda_k)}$ ($0< \lambda_j, \lambda_k < m+1$)
then what we want to prove is already true.

Now Lemma 3.4 in \cite{DW4} and our hypothesis $J(\bar{v}, \bar{w}) > 0$
imply that $J_{00} < 0$ and $J_{\lambda_i, \lambda_j} > 0$ except in the
three cases listed there. The proof that the elements of $\mathcal C$ are
equi-distributed in $\bar{v}\bar{w}$ carries over from \cite{DW4} since
the midpoint $\bar{u}$ is not involved in the arguments.

Suppose next that the points in ${\mathcal C} \cap \bar{v}\bar{w}$ are
equi-distributed. In the special case where $m=1$, we have $J(\bar{v}, \bar{u}) =0
=J(\bar{u}, \bar{w})$, which imply $J(\bar{u}, \bar{u}) = 0$. So
the midpoint does not contribute to the equation from $c^{(0)} + c^{(\lambda_{m+1})}$.
If $m > 1$, we write down the equations arising from $c^{(0)} + c^{(\lambda_{m+1})}$
and $c^{(\lambda_{m-1})} + c^{(\lambda_{m+1})}$. The formula for $F_{\lambda_j}$
in \cite{DW4} still holds for $1 \leq j \leq m$, and using this and the second equation
we obtain the analogous formula for $F_{\lambda_{m+1}}$.

Putting all the above information together in the first equation
and using $A_u < 0$, we see that $F_{\lambda_1}^{m+1}/F_0^{m-1}$ is positive
if $m$ is even and negative if $m$ is odd. In either case it follows
immediately that $F_0 F_{\lambda_{m+1}} < 0$, as required.
\end{proof}

\medskip
We shall now set up the basic machinery of the projection of
our convex hull onto an affine hyperplane.

Let $\bar{c}$ be a null vector in $\mathcal C$ and let
$H$ be an affine hyperplane separating $\bar{c}$ from
${\rm conv}(\frac{1}{2}(d + {\mathcal W}))$. Define a map
$P : {\rm conv}(\frac{1}{2}(d + {\mathcal W})) \longrightarrow H$ by letting
$P(\bar{z})$ be the intersection point of the ray $\bar{c} \bar{z}$ with $H$.
We denote by $\Delta$ the image of $P$ in $H$. ($P$ and $\Delta$
of course depend on $\bar{c}$ and the choice of $H$. When considering
projections from several null vertices, we will use the vertices
as  superscripts to distinguish the cases, e.g., $\Delta^{\bar{c}},
\Delta^{\bar{b}}$.)

Let us now consider a vertex $\bar{\xi}$ of $\Delta$.  We know that
$\bar{c}$ and $\bar{\xi}$ are collinear with a subset $P^{-1}(\bar{\xi})$ of
${\rm conv}(\frac{1}{2}(d + {\mathcal W}))$. As $\bar{\xi}$ does not lie in
the interior of a positive-dimensional subset of $\Delta$, we see that no
point of $P^{-1}(\bar{\xi})$ lies in the interior of a subset of
${\rm conv}(\frac{1}{2}(d + {\mathcal W}))$ of dimension $>1$. So
$P^{-1}(\bar{\xi})$ is a vertex or an edge of ${\rm conv}(\frac{1}{2}(d +
{\mathcal W})).$

If $P^{-1}(\bar{\xi})$ is a vertex $\bar{x}$, then $2\bar{x} \in
d + {\mathcal W}$ and in Lemma \ref{III44} we can take the affine subspace
$A$ to be the line through $\bar{c},\bar{\xi},\bar{x}$.
Using this lemma and also Prop \ref{E} and
Cor \ref{meet} we see that either $\bar{x}
\in \mathcal C$ (in which case $ J(\bar{x},\bar{c})=0)$, or $\bar{x}
\notin \mathcal C$ and $\bar{x} = (\bar{a}+ \bar{c})/2$ for some null
element $\bar{a} \in {\mathcal C} \cap A$. We have therefore deduced

\begin{theorem} \label{Vertex}
Let $\bar{\xi}$ be a vertex of $\Delta$.  Then exactly one of
the following must hold:

$($1A$)$ $\bar{\xi}$ $($and hence $P^{-1}(\bar{\xi})$$)$ is orthogonal
to  $\bar{c}$;

$($1B$)$ The line through $\bar{c},\bar{\xi}$ meets ${\rm conv}(\frac{1}{2}(d +
{\mathcal W}))$ in a unique point $\bar{x}$,
 and there exists a null $\bar{a} \in {\mathcal C}$  such that
 $(\bar{a}+\bar{c})/2 = \bar{x}$;

$($2$)$  $\bar{\xi}$ is not orthogonal to $\bar{c}$, and
$\bar{c}$ and $\bar{\xi}$ are collinear with an edge $\bar{v} \bar{w}$ of
 ${\rm conv}(\frac{1}{2}(d + {\mathcal W}))$, $($and hence $c$ and $\xi$
 are collinear with the corresponding edge $vw$ of
 ${\rm conv}(\mathcal W)$$)$. \  $\qed$
\end{theorem}

\begin{rmk} \label{xc}
If $($1B$)$ occurs, then $\bar{a}=2 \bar{x}-\bar{c}$ being null
is equivalent to $J(\bar{x},\bar{x})=J(\bar{x},\bar{c})$, that is,
\begin{equation} \label{xceqn}
\sum_{i=1}^{r} \frac{x_i^2}{d_i} = \sum_{i=1}^{r} \frac{x_i c_i}{d_i}
\end{equation}
\end{rmk}
In particular $x_i$ and $c_i$ are nonzero for some common index $i$.
We will from now on refer to this situation by saying that the
vectors $x$ and $c$ {\it overlap}.

\medskip

We make a preliminary remark about (1A) vertices.

\begin{lemma} \label{IX6} Suppose that $u \in {\mathcal W}$ and
$\bar{u} \in \bar{c}^{\perp}$.

$($a$)$ If $u = (-2^i, 1^j)$  then $c_i \neq 0$.

$($b$)$ Suppose that $K$ is connected. If $u = (-1^i, -1^j, 1^k),$
   then $c_i, c_j, c_k$ are all nonzero.
\end{lemma}
\begin{proof}
After a suitable permutation, we may let $1, \cdots, s$ be the
indices $a$ with $c_a \neq 0$. We need
\[
\sum_{a=1}^{s} \frac{u_a c_a}{d_a} = \sum_{a=1}^{s} \frac{c_a^2}{d_a}= 1.
\]
In case (a) this is impossible if $c_i=0$
(that is, $i \notin \{1, \cdots, s \}$) as then we need
$d_j=1=c_j$ and $c_a =0$ for $a \neq j$, contradicting
$\sum_{k=1}^{r} c_k =-1$.

Next, Cauchy-Schwartz on $(\frac{u_a}{\sqrt{d_a}})_{a=1}^{s},
(\frac{c_a}{\sqrt{d_a}})_{a=1}^{s}$
shows $\sum_{a=1}^{s} \frac{u_a^2}{d_a} \geq 1$.
 In  case (b), if, say, $c_k =0,$ then since
 $\frac{1}{d_i} + \frac{1}{d_j} \geq 1$ and $d_i, d_j \geq 2$
(see Remark \ref{facts}(d)) we must have $d_i = d_j =2$.
The equations then imply $c_i=c_j =-1$ and $c_a =0$
for $a \neq i,j$, also giving a contradiction. Similar arguments
rule out $c_i = 0$ or $c_j = 0$.
\end{proof}

\smallskip
In the next two sections we shall get stronger results on (1A) vertices.
Let us now consider type (2) vertices.

\begin{theorem} \label{cedge}
Consider a type $($2$)$ vertex $\bar{\xi}$ of $\Delta$. So $c$ and
$\xi$ are collinear with an edge $vw$ of ${\rm conv}(\mathcal W)$.
Suppose there are no points of $\mathcal W$ in the interior of
$vw$. Then we have

$($i$)$  $c= 2v-w$  or

$($ii$)$ $c = (4v-w)/3.$

In $($i$)$ the points of $\mathcal C$ on the line through $\bar{c}, \bar{\xi}$
are $\bar{c}$ and $\bar{w}$. In $($ii$)$ they are
$\bar{c},  \bar{w}$ and $\bar{c}^{(1)} = (2\bar{v}+\bar{w})/3 =
(\bar{c}+\bar{w})/2$. We need $J(\bar{c}^{(1)},\bar{w})=0$.
\end{theorem}

\begin{proof}
This is very similar to the arguments of \S 3 in \cite{DW4}.
We apply Lemma \ref{III44} to the line through $\bar{v}, \bar{w}$.

(A) We write the elements of $\mathcal C$ on the line as $\bar{c}=
\bar{c}^{(0)}, \bar{c}^{(1)}, \cdots, \bar{c}^{(m+1)}$ with $m \geq 0$.
So $\bar{c}^{(m+1)}$ is either null or is $\bar{w}$. No other
$\bar{c}^{(j)}$ can lie beyond $\bar{w}$, by Cor \ref{meet}.

By assumption  $\bar{c}=\bar{c}^{(0)}$ is not orthogonal to the whole line.
 As $\bar{c}$ is null, this means $\bar{c}$ is not orthogonal to
 {\em any} other point on the line. So $\bar{c}^{(0)} + \bar{c}^{(j)}$ is either
$2 \bar{v}$, $2 \bar{w}$ or else is  a sum of two other $\bar{c}^{(i)}$.
In particular, $\bar{c}^{(0)} + \bar{c}^{(1)} = 2\bar{v}$.
In fact $\bar{c}^{(0)} + \bar{c}^{(j)}$ is never $2\bar{w}$;
for the only possibility is for $\bar{c}^{(0)} + \bar{c}^{(m+1)} =2\bar{w},$
in which case $\bar{c}^{(m+1)}$ is null, and so
 $\bar{c}^{(m)} + \bar{c}^{(m+1)}= 2 \bar{w}$, contradicting $\bar{v} \neq \bar{w}$.

We deduce that for $j > 1$, we have $\bar{c}^{(0)} + \bar{c}^{(j)}
=\bar{c}^{(k)} + \bar{c}^{(p)}$ for some $1 \leq k,p \leq j-1$.

(B) Let $\bar{c}^{(m+1)}$ be null. Since the segment
$\bar{c}^{(0)} \bar{c}^{(m+1)}$ lies in the interior of the null ellipsoid,
Lemma \ref{Jpos} implies that $J(\bar{c}^{(i)},\bar{c}^{(j)}) > 0$
unless $i=j=0$ or $m+1$.  Arguments very similar to those in \S 3
of \cite{DW4} enable us  to determine the signs of the $F_{\bar{c}^{(j)}}$ in
(\ref{potform}) and show that the contributions from the pairs summing
to $\bar{c}^{(1)} + \bar{c}^{(m+1)}$ cannot cancel. So we have a
contradiction unless $\bar{c}^{(1)} + \bar{c}^{(m+1)}= \bar{w}$, which can only
happen if $m=1$, i.e., $\bar{c}^{(0)} + \bar{c}^{(1)} =2 \bar{v}, \;
 \bar{c}^{(1)} + \bar{c}^{(2)} =2 \bar{w}$ and $\bar{c}^{(0)} + \bar{c}^{(2)} =
2 \bar{c}^{(1)}$ (otherwise $\bar{c}^{(0)} + \bar{c}^{(2)}$ cannot cancel).
Hence we have
\[
c = (3v-w)/2 \;\; ; \;\;   c^{(1)} = (v+w)/2   \;\; ; \;\
c^{(2)} =(3w-v)/2.
\]
Writing $F_j$ for $F_{\bar{c}^{(j)}}$, we need $2 F_0 F_2
J(\bar{c},\bar{c}^{(2)}) + F_1^2 J(\bar{c}^{(1)}, \bar{c}^{(1)})=0$ so that the
contributions from $\bar{c}^{(0)} + \bar{c}^{(2)}$ and $\bar{c}^{(1)} + \bar{c}^{(1)}$
cancel.  As $J(\bar{c}, \bar{c}^{(2)})$ and $J(\bar{c}^{(1)}, \bar{c}^{(1)}) > 0$,
we need $F_0$ and $F_2$ to have opposite signs.  Now, as
$J(\bar{c},\bar{c}^{(1)}), J(\bar{c}^{(1)}, \bar{c}^{(2)})>0$, we see that $A_v$
and $A_w$ have opposite signs.  So we may let $w$ be type I and
$v$ be type II or III, as long as the asymmetry between $\bar{c}^{(0)}$
and $\bar{c}^{(2)}$ is removed.  Note that $v,w$ cannot overlap
if $v$ is type II, as then Remark \ref{facts}(c) means $w$ is not a vertex.
The possibilities are (up to permutation)
\[
\begin{array}{|l|c|c|l|l|}
    \hline
    &   v       &    w    &    c^{(0)} = \frac{1}{2}(3v-w)  &   c^{(2)} = \frac{1}{2}(3w-v)  \\ \hline
(1) &   (-2,1,0,\cdots) & (-1,0,\cdots)      &  (-\frac{5}{2},\frac{3}{2},0,\cdots) &
(-\frac{1}{2},-\frac{1}{2},0,\cdots) \\
(2) &     (-2,1,0, \cdots) & (0,-1,0,\cdots)   & (-3,2,0,\cdots) &
(1,-2,0,\cdots) \\
(3) &     (-2,1,0, \cdots) & (0,0,-1,0,\cdots) &  (-3, \frac{3}{2},\frac{1}{2},\cdots)
& (1, -\frac{1}{2}, -\frac{3}{2} , 0, \cdots) \\
(4) &    (1,-1,-1,0,\cdots) & (0,0,0,-1,\cdots) & (\frac{3}{2}, -\frac{3}{2},-\frac{3}{2}, \frac{1}{2}, \cdots)
 & (-\frac{1}{2}, \frac{1}{2}, \frac{1}{2}, -\frac{3}{2}, \cdots) \\ \hline
\end{array}
\]
Now, it is clear in (1) and (2) that $\bar{c}^{(0)}$ and $\bar{c}^{(2)}$
can't both give null vectors.  For (3) and (4), we find that the nullity
equations for $\bar{c}^{(0)}$ and $\bar{c}^{(2)}$ have no integral solutions
in $d_i$ (in fact $d_3$ (resp. $d_4$) must be $5/2$).

Therefore in fact $\bar{c}^{(m+1)}$ cannot be null.

(C) Now suppose that $\bar{c}^{(m+1)}=w$ and $m >0$.
Since $J(\bar{c}, \bar{v}) \neq 0$, $\bar{v}$ must lie between $\bar{c}^{(0)}$
and $\bar{c}^{(1)}$. So $J(\bar{c}^{(0)}, \cdot)$ and $J(\ \cdot \ , \bar{c}^{(m+1)})$
are affine functions on the line, vanishing at $\bar{c}^{(0)}$ and $\bar{c}^{(m)}$
respectively. Hence $J(\bar{c}^{(0)}, \bar{c}^{(i)}) \; (i \geq 1)$
and $J(\bar{c}^{(i)}, \bar{c}^{(m+1)}) \; (0 \leq i \leq m-1)$
are the same sign as $J(\bar{c}^{(0)}, \bar{c}^{(1)})$.
It follows that  $J(\bar{c}^{(i)}, \cdot)$ is an affine function on
the line, taking the same sign as $J(\bar{c}^{(0)},\bar{c}^{(1)})$
at $\bar{c}^{(0)}, \bar{c}^{(m+1)}$ (for $1 \leq i \leq m-1$).
Thus $J(\bar{c}^{(i)}, \bar{c}^{(j)})$ is the
same sign as $J(\bar{c}^{(0)}, \bar{c}^{(1)})$ except for the cases
\[
J(\bar{c}^{(0)}, \bar{c}^{(0)})=0=J(\bar{c}^{(m)}, \bar{c}^{(m+1)})
\;\;\; : \;\;\; {\rm sign} \ J(\bar{c}^{(m+1)}, \bar{c}^{(m+1)}) =
-{\rm sign} \ J(\bar{c}^{(0)}, \bar{c}^{(1)}).
\]
It then follows that the sign and non-cancellation arguments
of (B) (taken from \S 3 of \cite{DW4}) still hold, except in the case
$m=1$.

These give the two cases of the Theorem.  If $m=0$, we have
$c^{(1)} = w$ and $c^{(0)} =2v-w$ as $c^{(0)} + c^{(1)} =2v$.
If $m=1$, then $c^{(2)} =w, \;c^{(0)} + c^{(1)} =2v$ and
$c^{(0)} + c^{(2)} =2c_1$ (for cancellation). Hence
$c^{(0)} = (4v-w)/3,$  $c^{(1)} = (2v+w)/3$, as well as
$J(\bar{c}^{(1)}, \bar{c}^{(2)})=0$.
\end{proof}

\begin{rmk} \label{2int} If there are points of $\mathcal W$ in the
interior of $vw$, we can still conclude that $c^{(0)} + c^{(1)} =2v$.  Hence
$c = \lambda v + (1-\lambda)w$ for $1 < \lambda \leq 2$, since if
$\lambda > 2$ then $\bar{c}^{(1)}$ is beyond $\bar{w}$. It must then be null,
and $m=0$, so there is no way of getting $2\bar{w}$ as a sum of two
elements in $\mathcal C$.
\end{rmk}

\begin{lemma} \label{w}
For case $($i$)$ in Theorem \ref{cedge} $($i.e., $c=2v-w$$)$,
either $w$ is type I, or $w$ is type III and $v_i=-1, w_i=-2$
for some index $i$.
\end{lemma}

\begin{proof}
It follows from above  that $J(\bar{w}, \bar{w}) F_{\bar{w}}^2 =
A_w$ so $J(\bar{w}, \bar{w})$ is positive if $w$ is type I and
negative if $w$ is type II or III. In the latter case,
 $\sum_i \frac{w_i^2}{d_i} > 1,$
but by nullity, $c=2v-w$ satisfies $\sum \frac{c_i^2}{d_i} =1$. Hence
for some $i$ we have $|w_i| > |c_i| = |2v_i - w_i|$. As $v_i, w_i \in
\{-2,-1,0,1\}$, it follows that $v_i=-1, w_i=-2$.
\end{proof}

\smallskip

We are now able to characterise the case where $c$ is a type I vector.

\begin{theorem} \label{Fano}
If $c$ is a type I vector, say $(-1, 0, \cdots)$ for definiteness,
then ${\mathcal W}$ is given by $ \{(-1)^i, (1^1, -2^i) : i =2, \cdots, r \}$.
\end{theorem}

\begin{rmk}
Equivalently, $\mathcal W$ is as in Ex 8.1 of \cite{DW4}, where the
hypersurface in the Ricci-flat manifold is a circle bundle over a
product of K\"ahler-Einstein Fano manifolds. A superpotential was found for
this example in \cite{CGLP3}.
\end{rmk}

\begin{proof}
Nullity of $\bar{c}$ implies $d_1=1$, so $(-2^1, 1^i) \notin \mathcal W$.
Also $(-1^1, -1^j, 1^k) \notin \mathcal W$, as then $c$ would be in conv$({\mathcal W})$.
Let us consider the vertices $\bar{\xi}$ in $\Delta$. $\bar{\xi}$ cannot be of
type (1A); otherwise $\xi_1 = -1$, which implies the existence of a type II
vector in $\mathcal W$ with a nonzero first component, contradicting the above.
There can also be no $\bar{\xi}$ of type (1B)  since by Remark \ref{xc} the vector
$\bar{x}$ satisfies $0 < -x_1$, which we ruled out above.

Hence all vertices of $\Delta$ are of type (2), i.e., correspond to edges $vw$
of conv$(\mathcal W)$ such that $c = \lambda v + (1- \lambda)w$ and
  $\lambda > 1$. From this equation it follows that $v,w$ are of the form
\[
v = (  -1^i ), \ w = (1^1, -2^i)
\]
for some $i>1$. As $\Delta$ (being a $(r-2)$-dimensional
polytope in an ($r-2$)-dimensional affine space)
has at least $r-1$ vertices, such vectors occur for all $i \neq 1$.

Now no type II vector can be in $\mathcal W$, otherwise $v$ would not
be a vertex. Also $(1^i, -2^j)$ with $i,j \neq 1$ cannot be in
$\mathcal W$, as then $(-1^j)$ would not be a vertex. We have already seen
 $(-2^1, 1^i)$ is not in $\mathcal W$.  So $\mathcal W$ is as
claimed.
\end{proof}

\medskip

We shall henceforth exclude this case, i.e. case (i) of Theorem
\ref{classthm}, from our discussion.

We conclude this section by giving a preliminary listing of the
possibilities for $c$ when we have a type (2) vertex. These are given by
 cases (i) and (ii) of Theorem \ref{cedge}, as well as the possible cases
when there is a point of $\mathcal W$ in the interior of $vw$.

For Theorem \ref{cedge}(i) the possible $v,w,c$ are:
\[
\begin{array}{|c|c|c|c|} \hline
    &        v          &           w         &     c=2v-w        \\ \hline
(1) & (-1,1,-1,\cdots)  &    (-2,1, \cdots)  &     (0,1,-2,\cdots) \\
(2) & (-1,-1,1, \cdots) &   (-2,1,\cdots)   &      (0,-3,2,\cdots) \\
(3) & (-1,0,-1,1,\cdots) &  (-2,1,\cdots)   &      (0,-1,-2,2,\cdots) \\
(4) & (-2,1,\cdots)     &   (-1,0,\cdots)   &      (-3,2,\cdots) \\
(5) & (-2,1,\cdots)     &    (0,0,-1,\cdots) &     (-4,2,1,\cdots)\\
(6) & (-1,0,\cdots)     &    (0,-1,\cdots)   &     (-2,1,\cdots) \\
(7) & (1,-1,-1,\cdots)  &    (0,0,0,-1,\cdots) &   (2,-2,-2,1, \cdots) \\ \hline
\end{array}
\]
\centerline{Table 1: $c = 2v-w$ cases}

\smallskip

\noindent{where $\cdots$  denotes zeros as usual.} To arrive at this list,
recall from Lemma \ref{w} that $w$ is either type I or type III with
$v_i=-1, w_i=-2$ for some $i$. Note also that if $w$ is type I and $v$
is type II then $v$ cannot overlap with $w$ as $w$ cannot then be a vertex.
Furthermore, the other possibility with $w$ type I and $v$ type III is
excluded as we are assuming in Theorem \ref{cedge} that there are no points of
$\mathcal W$ in the interior of $vw$.  Finally, the case $w=(-2,1,\cdots),
v=(-1,0,\cdots)$ can be excluded as this just gives the example in
Theorem \ref{Fano}.

In order to list the possibilities under Theorem \ref{cedge}(ii), recall that we
need $J(\bar{c}^{(1)}, \bar{w})=0$ where $c^{(1)} = (2v+w)/3$. Equivalently,
we need
\begin{equation} \label{4/3cond}
2 J(\bar{v}, \bar{w}) + J(\bar{w}, \bar{w})=0.
\end{equation}
This puts constraints on the possibilities for $v, w$. For instance,
$w$ cannot be type I, as for such vectors $J(\bar{v}, \bar{w}) \geq 0$ and
 $J(\bar{w},  \bar{w})>0$. Also,  if $w$ is type II or III, then from
 the superpotential equation we need $J(\bar{w},  \bar{w})<0,$ so
 $J(\bar{v}, \bar{w})>0$. If $w$ is type III, say $(-2, 1, 0, \cdots)$,
 then since $d_1 \geq 2,$ we have $\frac{4}{d_1} + \frac{1}{d_2} \leq 3$,
and the above equation gives $J(\bar{v}, \bar{w}) \leq \frac{1}{4}$ with equality iff
$d_1=2, d_2 =1$.

By the above remarks and the nullity of $\bar{c}$, after a moderate amount
of routine computations, we arrive at the following possibilities, up to
permutation of entries. In the table we have listed only the minimum
number of components for each vector and all unlisted components are zero.
Note that the entries (12)-(16) can occur only if $K$ is not connected (cf.
Remark \ref{Kconn}).

\[
\begin{array}{|c|c|c|c|c|} \hline
     &    v     &      w     &      c^{(1)} = (2v+w)/3      & c = (4v-w)/3 \\ \hline
(1)  &   (0,0,-2,1)  & (-2,1,0,0) & ( -\frac{2}{3}, \frac{1}{3}, -\frac{4}{3}, \frac{2}{3}) & (\frac{2}{3}, -\frac{1}{3}, -\frac{8}{3} , \frac{4}{3}) \\
(2)  &   (-2,0,1) & (-2,1,0) & (-2, \frac{1}{3}, \frac{2}{3}) & (-2, -\frac{1}{3}, \frac{4}{3}) \\
(3) & (-1,0,0,) & (-2,1,0) & (-\frac{4}{3}, \frac{1}{3},0) & (-\frac{2}{3}, -\frac{1}{3},0) \\
(4) & (0,0,-1) & (-2,1,0) & (-\frac{2}{3}, \frac{1}{3}, -\frac{2}{3}) & (\frac{2}{3}, -\frac{1}{3}, -\frac{4}{3}) \\
(5) & (-1,0,1,-1) & (-2,1,0,0) &  (-\frac{4}{3}, \frac{1}{3}, \frac{2}{3}, -\frac{2}{3}) & ( -\frac{2}{3},-\frac{1}{3},
\frac{4}{3}, -\frac{4}{3}) \\
(6) &   (-1,1,-1) & (-2,1,0)   &  (-\frac{4}{3}, 1, -\frac{2}{3}) & (-\frac{2}{3}, 1,-\frac{4}{3}) \\
(7) & (0,0,1,-1,-1) & (-2,1,0,0,0) & (-\frac{2}{3}, \frac{1}{3}, \frac{2}{3}
-\frac{2}{3},-\frac{2}{3}) & (\frac{2}{3},-\frac{1}{3},\frac{4}{3},-\frac{4}{3},-\frac{4}{3}) \\
(8) & (0,1,-1,-1) & ( -2,1,0,0) & (-\frac{2}{3}, 1, -\frac{2}{3}, -\frac{2}{3}) &
(\frac{2}{3}, 1, -\frac{4}{3},-\frac{4}{3}) \\
(9) & (0,-1,-1,1) & (1,-1,-1,0) & (\frac{1}{3}, -1, -1, \frac{2}{3}) & (-\frac{1}{3}, -1,-1,
\frac{4}{3}) \\
(10) & (1,-1,0,-1) & (1,-1,-1,0) & (1,-1,-\frac{1}{3}, -\frac{2}{3}) & (1,-1,\frac{1}{3},
-\frac{4}{3}) \\
(11) & (1,-2,0) & ( 1,-1,-1)  & (1, -\frac{5}{3}, -\frac{1}{3}) & (1, -\frac{7}{3}, \frac{1}{3})\\
 \hline
(12) & (1, 0, 0, -2) & (1, -1, -1, 0) & (1, -\frac{1}{3}, -\frac{1}{3}, -\frac{4}{3})
           & (1, \frac{1}{3}, \frac{1}{3}, -\frac{8}{3}) \\
(13) & (0, 0, 0, -2, 1) & (1, -1, -1, 0, 0) & (\frac{1}{3}, -\frac{1}{3}, -\frac{1}{3}, -\frac{4}{3},
         \frac{2}{3})  & (-\frac{1}{3}, \frac{1}{3}, \frac{1}{3}, -\frac{8}{3},
         \frac{4}{3})  \\
(14) & (0,-1,0,1,-1)  &  (1, -1, -1, 0,0) & (\frac{1}{3}, -1, -\frac{1}{3}, \frac{2}{3}, -\frac{2}{3})
          & (-\frac{1}{3}, -1, \frac{1}{3}, \frac{4}{3}, -\frac{4}{3})\\
(15) & (1, 0,0,-1,-1) & (1, -1,-1,0,0) & (1, -\frac{1}{3}, -\frac{1}{3}, -\frac{2}{3}, -\frac{2}{3})
          & (1, \frac{1}{3}, \frac{1}{3}, -\frac{4}{3}, -\frac{4}{3}) \\
(16) & (0,0,0,1,-1,-1) & (1,-1,-1,0,0,0) & (\frac{1}{3}, -\frac{1}{3}, -\frac{1}{3},
        \frac{2}{3}, -\frac{2}{3}, -\frac{2}{3}) & (-\frac{1}{3}, \frac{1}{3}, \frac{1}{3}, \frac{4}{3},
        -\frac{4}{3}, -\frac{4}{3}) \\ \hline
\end{array}
\]

\centerline{Table 2: $c= \frac{1}{3}(4v-w)$ cases}

\smallskip

We will also need a listing of those  cases
for which $vw$ has interior points lying in conv$(\mathcal W)$.

\[
\begin{array}{|c|c|c|c|} \hline
    &           v       &         w             &            c  \\   \hline
(1) &   (1,-2,\cdots)   &   (-2,1, \cdots)      &      (3 \lambda -2, \ 1-3 \lambda, \cdots)          \\
(2) &   (1,-2,\cdots)   &   (-1,0,\cdots)       &     (2 \lambda-1, \ -2 \lambda, \cdots)             \\
(3) &   (-1,0, \cdots) &   (1,-2,\cdots)        &     (1-2 \lambda,\  2 \lambda -2,  \cdots)          \\
(4) &   (-2,1,0, \cdots) &  (0,1,-2,\cdots)     &     (-2 \lambda, 1, \ 2 \lambda -2,\cdots)         \\
(5) &   (1,-1,-1,\cdots) & (-1,1,-1,\cdots)     &     (2 \lambda -1, 1 - 2 \lambda, -1, \cdots) \\ \hline
\end{array}
\]

\centerline{Table 3: Cases with interior points}

\smallskip

{\noindent Recall from} Remark \ref{2int} that $1 < \lambda \leq 2$ and $\cdots$ denote zeros.
Note that except in (4) all interior points which may lie in $\mathcal W$
actually do.

\section{\bf The sign of $J(\bar{c}, \bar{w})$}

\begin{theorem} \label{Jsign}
 ${\rm conv}(\frac{1}{2}(d+ {\mathcal W}))$  lies  in the closed half-space
$J(\bar{c}, \cdot) \geq 0$, i.e., the same closed half-space in which the
null ellipsoid lies.
\end{theorem}

\begin{proof}
We know that if $\bar{\xi}$ is a vertex of $\Delta^{\bar{c}}$ then
there are three possibilities, given by (1A), (1B) and (2) of Theorem
\ref{Vertex}.  If (1A) occurs, then by definition $J(\bar{c},
\bar{\xi}) = 0$.  If (1B) occurs, let $\bar{a}$ be the null vector
in Theorem \ref{Vertex}. Then by Lemma \ref{Jpos}, $J(\bar{c}, \bar{a})>0,$
which in turn implies that $J(\bar{c}, \bar{\xi}) > 0$.  It is now enough
to show that $J(\bar{c}, \bar{\xi}) \geq 0$ if $\bar{\xi}$ is a type
(2) vertex of $\Delta^{\bar{c}}$, since it then follows that
$\Delta^{\bar{c}}$, and hence ${\rm conv}(\frac{1}{2}(d+ {\mathcal W}))$,
lies in the half-space  $J(\bar{c}, \cdot) \geq 0$.

\smallskip

Suppose then that $\bar{\xi}$ is a type (2) vertex with $J(\bar{c},
\bar{\xi}) < 0$. By Remark \ref{2int}, $c = \lambda v + (1- \lambda)w$
for some $v, w \in  \mathcal W$ with $1< \lambda \leq 2$,  and both
$J(\bar{c},\bar{v}), J(\bar{c}, \bar{w}) < 0$.
In particular, from Remark \ref{tangent} and Lemma \ref{Jpos},
  $J(\bar{v}, \bar{v}), J(\bar{w}, \bar{w}) < 0$ since
$\bar{v}, \bar{w}$ lie on the side of $\bar{c}^{\perp}$
opposite to the null ellipsoid.

 But
\begin{eqnarray*}
0 = 4J(\bar{c}, \bar{c})&=& J(d + \lambda v + (1- \lambda)w, \
d + \lambda v + (1- \lambda)w) \\
&=&  J(\lambda(d+v) + (1-\lambda)(d+w), \ \lambda(d+v) + (1-\lambda)(d+w))\\
&=& \lambda^2 J(d+v,d+v) + 2 \lambda(1-\lambda)J(d+v, d+w) + (1-\lambda)^2
                           J(d+w, d+w).
\end{eqnarray*}
It follows from the above remarks that $J(d+v, d+w) <0$, that is
\[
\sum_i^{} \frac{v_i w_i}{d_i} > 1.
\]
One then checks that this condition is only satisfied in the
following cases (up to permutation of indices and interchange of
$v$ and $w$):

(a)   $v = (-2,1,0, \cdots),  \ w = (-2,0,1,0, \cdots)$  with $ 1< d_1  < 4$;

(b)   $v = (-2,1,0,\cdots   ), \ w = (-1,1,-1,0, \cdots )$ with $d_1=2$, or
$(d_1,d_2)=(3,2)$, or $d_2=1$;

(c)   $v = (1,-1,-1, 0, \cdots ),  w = (1,-1,0, -1, 0, \cdots)$ with $d_1=1$ or
     $d_2=1$;

(d) $v=(1, -1, -1, 0, \cdots), w=(0, -1, -1, 1, 0, \cdots)$ with $d_2 =1$ or
     $d_3=1$.

In case (a), $c = (-2,\lambda, 1-\lambda,0, \cdots)$. The condition
$d_1 < 4$ is incompatible with the nullity of $\bar{c}$. Interchanging
$v$ and $w$ reverses only the role of $\lambda$ and $1-\lambda$.

A similar argument rules out case (b) with $v,w$ as shown,
 as here $c = (-\lambda -1, 1,
\lambda-1)$. If we interchange $v$ and $w$, then  $c= (\lambda-2, 1, -\lambda, \cdots)$.
Theorem \ref{cedge} tells us $\lambda = 4/3$ or $2$, so
$ c = (-2/3, 1, -4/3, \cdots) \;\; {\rm or} \;\; (0,1,-2,\cdots). $

In the former case $c^{(1)} := (2v+w)/3 = (-4/3, 1, -2/3, \cdots),$ so the
condition $J(\bar{w}, \bar{c}^{(1)})=0$ gives $8/3d_1 + 1/d_2 = 1$. Thus $(d_1, d_2) =
(3,9)$ or $(4,3)$ but in neither case is $\bar{c}$ null.  In the latter case
nullity means $1/d_2 + 4/d_3 = 1$, so $J(\bar{c},\bar{v}) =
\frac{1}{4}(1 - 1/d_2 - 2/d_3) >0,$
a contradiction.

In case (c), $c=(1, -1, -\lambda, \lambda-1)$ and if $v$ and $w$ are interchanged,
the last two components of $c$ are interchanged. But $\bar{c}$ cannot be null if
$d_1=1$ or $d_2 =1$. A similar argument works for case (d).
\end{proof}

\begin{corollary} \label{face}
${\rm conv}(\frac{1}{2}(d + {\mathcal W})) \cap \bar{c}^\perp$ is a
$($possibly empty$)$ face of ${\rm conv}(\frac{1}{2}((d + {\mathcal W})).$  \ $\qed$
\end{corollary}

This enables us to adapt Theorem 3.5 of \cite{DW4} to the
elements of $\bar{c}^{\perp}$.

\begin{corollary} \label{orthogonal}
Let $vw$ be an edge of conv$(\mathcal W)$ and suppose
 $\bar{v}$ and $\bar{w}$ are in $\bar{c}^\perp$.
  Suppose further that there are no elements of $\mathcal W$ in the
  interior of $vw$. Then $J(\bar{v}, \bar{w})=0$.
\end{corollary}

\begin{proof}
  This is essentially the same as the proof of Theorem 3.5 of
 \cite{DW1}.  As ${\rm conv}(\frac{1}{2}(d + {\mathcal W})) \cap
  \bar{c}^\perp$ is a face of ${\rm conv}(\frac{1}{2}(d + {\mathcal W}))$,
  Lemma \ref{III44} shows that for calculations in $\bar{c}^{\perp}$
  we need only consider elements of $\mathcal C$ in  this hyperplane.
Note that by Cor \ref{meet}, no elements of $\mathcal C$ lie on the
opposite side of $\bar{c}^{\perp}$ to
${\rm conv}(\frac{1}{2}(d + {\mathcal W})$.

Any vertex of conv$({\mathcal C}) \cap \bar{c}^{\perp}$ outside
${\rm conv}(\frac{1}{2}(d + {\mathcal W})) \cap \bar{c}^{\perp}$ is,
by Prop \ref{E}, null, so must be $\bar{c}$ by Lemma \ref{Jpos}.
Now Cor \ref{meet} shows that $\bar{c}$ is the only element of
conv$({\mathcal C}) \cap \bar{c}^{\perp}$ outside
${\rm conv}(\frac{1}{2}(d + {\mathcal W})) \cap \bar{c}^{\perp}$.
But any sum $\bar{c} + \bar{a}$ with $\bar{a} \in \bar{c}^{\perp}$ does not
contribute, so in fact we are in the situation of Theorem 3.5 of \cite{DW4}.
\end{proof}

\medskip
We introduce the following sets:
\begin{eqnarray*}
\hat{S}_1 &=& \{ i \in \{1, \cdots, r\}: \exists \ {\rm unique \;} w \in {\mathcal W}
\; {\rm with \;} \bar{w} \in \bar{c}^\perp \; {\rm and \;} w_i=-2 \}  \\
\hat{S}_{\geq 2} &=& \{i \in \{1, \cdots, r\} : \exists \ {\rm  more \; than \; one \;}
 w \in {\mathcal W} \; {\rm with \;} \bar{w} \in \bar{c}^\perp \; {\rm and \;}
 w_i=-2 \}
\end{eqnarray*}
These are similar to the sets $S_1$, $S_{\geq 2}$ of \cite{DW4}, but now we
require that the vectors $w$ to lie in $\bar{c}^\perp$.

It is immediate from Cor \ref{orthogonal} that $d_i=4$ if $i \in \hat{S}_{\geq 2}$,
(cf Prop 4.2 in \cite{DW4}).

We next prove a useful result about which elements of $\frac{1}{2}(d +
{\mathcal W})$ can be orthogonal to $\bar{c}$. This will give us
information about when (1A) vertices can occur.

\begin{lemma} \label{Zperp}
Assume that we are not in the situation of Theorem \ref{Fano}
 $($i.e., $c$ is not of
type I $)$. Let $u \in {\mathcal W}$ be such that $\bar{u} \in \bar{c}^{\perp}$.
 Then:

$($a$)$ there exists $i$ with $c_i \neq 0$ and $-2 < c_i < 1$;

$($b$)$ if $c \in \Z^r$ then there is at most one such $u$, and hence at
most one $($1A$)$ vertex $($wrt $c$$)$.
\end{lemma}

\begin{proof}
(a) The condition $J(\bar{u}, \bar{c})=0$ means
$\sum_{i} \frac{u_i c_i}{d_i} = 1$,
and nullity of $\bar{c}$ means
$\sum_{i} \frac{c_i^2}{d_i} = 1$.
As $u_i \in \{-2,-1,0,1 \}$, if the condition in (a) does not hold,
then $u_i c_i \leq c_i^2   \;\;\; {\rm for \; all} \; i$
so we must have equality for all $i$. Now $c_i = u_i$ for all $i$ with
$c_i$ nonzero. As $\sum c_i =-1$ and $c \neq u$ (since $c \notin \mathcal W$
by definition), this means $c$ is a type I vector
and we are in the situation of Theorem \ref{Fano}.

(b) We see from the previous paragraph that we need $u_i c_i > c_i^2$
for some $i$.  If $c \in \Z^r$ this means $c_i=-1$ and $u_i=-2$.
The orthogonality condition is now $\frac{2}{d_i} + \frac{c_j}{d_j} = 1$
where $u_j=1$. As $d_i \neq 1$ we see $c_j \geq 0$.

If $c_j=0$ then $d_i=2$. If $c_j >0$ then $d_i \geq 3$ so
$\frac{1}{3} \leq \frac{c_j}{d_j} < \frac{1}{c_j}$, where the second inequality
is due to the nullity requirement $\frac{1}{d_i} + \frac{{c_j}^{2}}{d_j}
 \leq 1$.  So $c_j =1$ or $2$. Moreover, the latter implies $(d_i,d_j)=(3,6)$
and $c = (-1^i, 2^j)$, which contradicts $\sum c_i=-1$.

We see that
either $c_j=1$ and $(d_i, d_j)=(4,2)$ or $(3,3)$, or $c_j=0$ and $d_i=2$.
Cor \ref{orthogonal} implies that if there is more than one such $u$
(say $(-2^i, 1^j)$ and $(-2^i, 1^k)$) for a given $i$, then $d_i=4$, so
$(d_i,d_j,d_k)=(4,2,2)$, and $(c_i,c_j, c_k)=(-1,1,1)$,
contradicting the nullity of $c$.

It now readily follows that the nullity condition prevents there being
more than one $u \in \mathcal W$ with $\bar{u} \in \bar{c}^{\perp}$
except when $c = (-1,-1,1,0,\cdots)$ with $d=(4,4,2,\cdots)$ or
$(3,3,3,\cdots)$ and $u= (-2,0,1,0,\cdots), (0,-2,1,0,\ldots)$. But in this
case if both $u$ occur then $c \in {\rm conv} ({\mathcal W})$, a contradiction.
\end{proof}

\medskip
We shall study (1A) vertices for non-integral $c$ in the next section.
The following results will be useful.

\begin{prop} \label{Lmod}
Let $v=(-2^i, 1^j)$ and $w=(-2^k, 1^l)$ be elements of $\mathcal W$ such that
$\bar{v}, \bar{w} \in \bar{c}^\perp$. Suppose that $i \in \hat{S}_1$ and
$\{i,j \} \cap \{k,l \} = \emptyset$. Then $k \in \hat{S}_{\geq 2}$
and $(d_i, d_k, d_l)=(2,4,2)$.
\end{prop}

\begin{proof}
By Remark \ref{facts}(e) the affine subspace
 $\{ \bar{x} :x_i + x_k = -2, x_j + x_l = 1 \} \cap \bar{c}^{\perp}$
meets ${\rm conv}(\frac{1}{2}(d +  {\mathcal W}))$ in a face, whose possible
elements are $v,w, u = (-2^k, 1^j), y= (-1^i, 1^j, -1^k)$ and $z = (-1^i, -1^k, 1^l)$
  (since $i \in \hat{S}_1$).

As $J(\bar{v}, \bar{w}) = \frac{1}{4}$, we see from Thm \ref{orthogonal} that
 $vw$ is not an edge so $z$ is present in the face.
 Now Cor \ref{orthogonal} on $vz$
implies $d_i=2$.  Also, $u$ must be present, otherwise $y$ is present and
Cor \ref{orthogonal} on $zw$ and $yw$ gives a contradiction. So
$k \in \hat{S}_{\geq 2}$, and Cor \ref{orthogonal}
on $uw$ implies $d_k=4$. Now considering $zw$ implies $d_l=2$.
\end{proof}

\begin{rmk}
This is similar to the proof of Prop 4.6 in \cite{DW4}.
But we cannot now deduce that $d_j=1$ as the proof of this in \cite{DW4}
relied on the existence of $t= (-1^i,-1^j,1^k)$, and although we
know this is in $\mathcal W$ we do not know if $\bar{t}$
 lies in $\bar{c}^\perp$.
\end{rmk}

\begin{prop} \label{AC}
If $i \in \hat{S}_{1}$ and $v=(-2^i, 1^j)$ gives an element of
$\bar{c}^\perp$ then $w=(-1^i, -1^j, 1^k)$ cannot give an element of
$\bar{c}^\perp$.
\end{prop}

\begin{proof}
This is similar to Prop. 4.3 in \cite{DW4}. Since $i \in \hat{S}_1$,
the vectors $\bar{v},\bar{w}$ lie on an edge in the face
 $\{ \bar{x} : 2x_i + x_j = -3 \} \cap \bar{c}^{\perp}$ of
${\rm conv}(\frac{1}{2}(d+\mathcal W))$, and $J(\bar{v}, \bar{w}) =
\frac{1}{4}( 1 -\frac{2}{d_i}+
\frac{1}{d_j}) \neq 0$ since $d_i\neq 1$.
\end{proof}

\begin{corollary} \label{AD}
With $v$ as in Prop \ref{AC}, there are no elements $w = (-2^j, 1^k)$
with $\bar{w}$ in $\bar{c}^\perp$.
\end{corollary}

\begin{proof}
This is similar to Prop 4.4 in \cite{DW4}. If $k=i$, then the type I
vector $u:=(-1^i)= \frac{1}{3}(2v+w)$ lies in $\mathcal W$ and
$\bar{u} \in \bar{c}^{\perp}$. By Lemma \ref{Iperp} below,  $u=c$,
contradicting $c \notin {\mathcal W}$.

We can therefore take $k \neq i$. Now $\bar{v},\bar{w}$
 lie on an edge in the face
$\{ \bar{x} : 3x_i + 2x_j =-4 \} \cap \bar{c}^{\perp}$
 (this is a face by Prop \ref{AC} and the assumption
$i \in \hat{S}_1$). But $J(\bar{v},\bar{w}) =
\frac{1}{4}( 1 + \frac{2}{d_j}) \neq 0$.
\end{proof}

\section{\bf Vectors orthogonal to a null vertex}

In this section we analyse the possibilities for
 $\frac{1}{2}(d + {\mathcal W}) \cap \bar{c}^\perp$.
This will give us an understanding of the vertices of type (1A).

We first dispose of the case of type I vectors.

\begin{lemma} \label{Iperp}
If $u$ is a type $I$ vector and $\bar{u} \in \bar{c}^\perp$ then $c=u$, so
we are in the situation of Theorem \ref{Fano}.
\end{lemma}

\begin{proof}
Up to a permutation we may let $u = (-1,0, \cdots)$. The
orthogonality condition implies $c_1 = -d_1$. But
nullity implies $\sum c_i^2/ d_i = 1$, so $d_1 =1$ and $c_i=0$ for $i > 1$.
(Note that in particular $u \notin {\mathcal W}$.)
\end{proof}

\medskip
We shall therefore assume from now on there are no type I vectors giving
points of $\bar{c}^{\perp}$.

\begin{lemma} \label{coplane}

$($i$)$ Two type II vectors whose nonzero entries lie in the same set of three
 indices cannot both give elements of $\bar{c}^\perp$.

$($ii$)$ Two type III vectors $(-2^i, 1^j)$ and $(1^i, -2^j)$ cannot
both give elements in $\bar{c}^\perp$.

$($iii$)$ Three type III vectors  whose nonzero entries all lie in the same
set of three indices cannot all give rise to elements in $\bar{c}^\perp$.
\end{lemma}

\begin{proof}
These all follow from Lemma \ref{Iperp} by exhibiting an affine combination
of the given vectors which is of type I.
\end{proof}

\medskip
Let $u, v \in \mathcal W$ be such that $\bar{u}$ and $\bar{v} \in \bar{c}^\perp$.
It follows that  $\lambda \bar{u} + (1- \lambda) \bar{v} \in \bar{c}^\perp$ for all
$\lambda$. Hence Remark \ref{tangent} shows that {\em for all} $\lambda$
\begin{eqnarray*}
 0 &\geq& J(d + \lambda u + (1-\lambda)v, d + \lambda u + (1-\lambda v)) \\
  &=&  J(\lambda (d+u) + (1-\lambda)(d+v), \lambda (d+u) + (1- \lambda)(d+v)) \\
  &=&  \lambda^2 (J(d+u,d+u) + J(d+v, d+v)  -2 J(d+u, d+v) ) + \\
  & &  2 \lambda(  J(d+u, d+v) - J(d+v, d+v))  + J(d+v, d+v).
\end{eqnarray*}
Equality occurs if and only if $\lambda u + (1- \lambda)v = c$, as
$\bar{c}$ is the only null vector in $\bar{c}^\perp$.

Multiplying by $-1$, using Eq.(\ref{Jform}), and recalling that the minimum
value of a quadratic $\alpha \lambda^2 + \beta \lambda + \gamma$ with
$\alpha > 0$ is $\gamma - (\beta^2/4 \alpha)$, we deduce
the following result.

\begin{lemma}  \label{estimate}
If $u, v \in \mathcal W$ and $\bar{u}, \bar{v} \in c^\perp$ then
\begin{equation} \label{est}
\sum_{i=1}^{r} \frac{u_i^2}{d_i} +
\sum_{i=1}^{r} \frac{v_i^2}{d_i}
-\left( \sum_{i=1}^{r} \frac{u_i^2}{d_i} \right) \left( \sum_{i=1}^{r} \frac{v_i^2}{d_i} \right)
\leq
 \sum_{i=1}^{r} \frac{u_i v_i}{d_i} \left(2- \sum_{i=1}^{r} \frac{u_i v_i}{d_i} \right).
\end{equation}
Moreover, equality occurs if and only if $c = \lambda u + (1- \lambda)v$
for some $\lambda$. \; $\qed$
\end{lemma}

\begin{rmk} \label{equal}
 By definition, $c$ does not lie in ${\rm conv}(\mathcal W)$. So in the case of
 equality in Eq.(\ref{est}) we cannot have $0 \leq \lambda \leq 1$. This observation
 will in many cases show that equality cannot occur.
\end{rmk}

\begin{rmk}
The right-hand side of Eq.(\ref{est}) is maximised when
$\sum \frac{u_i v_i}{d_i} = 1$ (i.e., $J(\bar{u}, \bar{v})=0$).
In this case Eq.(\ref{est}) just follows from
$\sum \frac{u_i^2}{d_i} , \sum \frac{v_i^2}{d_i} \geq 1$, which
is true for any two vectors in $\bar{c}^\perp$. If $J(\bar{u}, \bar{v})
 \neq 0$, we get sharper information.
\end{rmk}

\begin{corollary} \label{II/II} Suppose that $K$ is connected.
If $u, v$ are type II vectors in $\mathcal W$ with $\bar{u}, \bar{v} \in \bar{c}^\perp$ then
\begin{equation*} \label{estII}
\frac{1}{2} \leq \sum_{i=1}^{r} \frac{u_i v_i}{d_i} \leq \frac{3}{2} \ ,
\end{equation*}
with equality if and only if $c = \lambda u + (1- \lambda)v$ for
some $\lambda$, in which case all the $d_i=2$ whenever $i$ is an index
such that $u_i$ or $v_i$ is nonzero.
\end{corollary}

\begin{proof}
Writing $X = \sum \frac{u_i^2}{d_i}$ and $Y= \sum \frac{v_i^2}{d_i}$ we see that
$1 \leq X,Y \leq \frac{3}{2}$. The lower bound arises from
$\bar{u},\bar{v}$ being in $\bar{c}^\perp$, while the upper bound follows from
Remark \ref{facts}(d) and the assumption that $u,v$ are type II vectors.

Now $X+Y -XY=1-(1-X)(1-Y)$ is minimised for $X,Y$ in this range if
$X=Y = \frac{3}{2}$, when it takes the value $\frac{3}{4}$.  The
inequality Eq.(\ref{est}) now gives the result.
\end{proof}

\medskip
When $K$ is connected, it follows that any two such type II vectors must
overlap. Moreover, if they have only one common index then we are in the
case of equality in Cor \ref{II/II}. The nullity of $\bar{c}$ implies that
$\lambda = \frac{1}{2}$ in this case, contradicting Remark \ref{equal}.

Combining this remark with Cor \ref{estII} and Lemma \ref{coplane} (i)
, we deduce the following result.

\begin{corollary} \label{IIposs}
Assume $K$ is connected. If  $u, v$ are type II vectors in $\mathcal W$
with $\bar{u}, \bar{v} \in \bar{c}^\perp$, then either
$u = (-1^a,-1^b, 1^i), v = (-1^a, -1^b, 1^j)$ or
$u= (1^a, -1^b, -1^i), v= (1^a, -1^b, -1^j)$.

Hence the collection of all such type II vectors is of the form, for some
fixed $a,b$:

$($i$)$  $(-1^a,-1^b,1^i) :  i \in I$ for some set $I$; or

$($ii$)$ $(1^a, -1^b, -1^i) : i \in I$ for some set $I$; or

$($iii$)$ $(1,-1,-1,0, \cdots),\ (1,-1,0,-1,\cdots), \ (1,0,-1,-1,\cdots).$
\;  $\qed$
\end{corollary}

We now investigate type III vectors.

\begin{lemma} \label{II/III} Suppose $K$ is connected.
If $u$ is a type II vector and $v$ a type III vector in $\mathcal W$ with
$\bar{u}, \bar{v} \in \bar{c}^\perp$, then $\sum_{i=1}^{r} \frac{u_i v_i}{d_i} > 0.$
\end{lemma}

\begin{proof}
With the notation of Cor \ref{II/II} we have $1 \leq X \leq \frac{3}{2}$
and $1 \leq Y \leq 3$. So $X+Y-XY = 1- (1-X)(1-Y) \geq 0$, and Eq.(\ref{est})
gives the desired inequality. Also, the case of equality (i.e., $X= \frac{3}{2},
Y=3$) leads to $\lambda = \frac{2}{3}$, again contradicting Remark \ref{equal}.
\end{proof}

\begin{rmk} \label{Kconn}
While Cor \ref{II/II} - Lemma \ref{II/III} are stated under the
assumption that $K$ is connected, the actual property we used is that
in Remark \ref{disconnK}. By contrast, the next two results do not
require this property.
\end{rmk}

\begin{lemma} \label{III/III}
Any two type III vectors $u,v$ giving elements of
$\bar{c}^\perp$ must overlap.
\end{lemma}

\begin{proof}
Write $u=(-2^i, 1^j)$ and $v=(-2^k, 1^l)$. By Cor \ref{orthogonal},
if $i, k \in \hat{S}_{\geq 2}$ then $d_i = d_k =4$. Since
$J(\bar{c}, \bar{u})= 0$ we have (by Cauchy-Schwartz)
$$ 1= \left(-\frac{2c_i}{d_i} +\frac{c_j}{d_j}\right)^2 \leq \left(\frac{4}{d_i}+\frac{1}{d_j}\right)
\left(\frac{c_i^2}{d_i} + \frac{c_j^2}{d_j}\right).$$
Hence
$$ \frac{c_i^2}{d_i} + \frac{c_j^2}{d_j} \geq \frac{d_j}{1+d_j} \geq \frac{1}{2}.$$
If $u$ and $v$ do not overlap, then the above and the analogous result
from considering $J(\bar{c}, \bar{v})=0$, together with the nullity of $\bar{c}$,
imply that $d_j=1=d_l$ and the only nonzero components of $c$ are
$c_i=c_k= -1, c_j=c_l = \frac{1}{2}$. But then $c$ is the midpoint of
$uv$,  contradicting $c \notin {\rm conv}(\mathcal W)$.

So if $u$ and $v$ do not overlap, we can take
 $i \in \hat{S}_1$.  Proposition \ref{Lmod} shows that $k \in \hat{S}_{\geq 2}$
and $(d_i, d_k, d_l) =(2,4,2)$. Hence $2 < X \leq 3$ and $ Y = \frac{3}{2}$,
so $X+Y-XY \geq 0$ and $\sum \frac{u_i v_i}{d_i} \geq 0$. Non-overlap means
that equality holds. But then $\lambda = 1/3$, contradicting Remark \ref{equal}.
\end{proof}

\medskip
Lemma \ref{III/III}, together with Lemmas \ref{coplane} and \ref{AD},
implies the following corollary.

\begin{corollary} \label{IIIposs}
The type III vectors associated to elements of $\frac{1}{2}(d + {\mathcal W})
 \cap \ \bar{c}^\perp$
are, up to permutation of indices, either of the form

$(a) \; (-2^1, 1^i), \ \; i \in I$,  $($with $d_1=4$ if $|I| \geq 2$$)$, or

$(b) \; (1^1, -2^i), \; i \in I$,

\noindent{for some subset} $I \subset \{2, \cdots, r\}$. \; $\qed$
\end{corollary}

Having found the possible configurations for type III vectors in
$\bar{c}^{\perp}$, we start to analyse the type II vectors for each
such configuration. For the rest of this section we will assume that
$K$ is connected (cf Remark \ref{Kconn}).

\begin{rmk} \label{IIrmk} Lemma \ref{II/III} now shows that in case (a)
of Cor \ref{IIIposs}, if $|I| \geq 2$, then every type II vector associated
to an element of $\bar{c}^\perp$ must have ``-1'' in place 1.
Similarly, in case (b), if $|I| \geq 3$, then every such type II vector
has ``1'' in place 1.  (So if a type II is present then $d_1 \neq 1$).
If $|I|=2$, the only possible type II vectors with ``0'' in place 1 are
$(0^1,-1^2,-1^3, 1^i)$ where $i \geq 4$, and all type II vectors whose
first entry is nonzero actually must have first entry equal to $1$.
\end{rmk}

\begin{lemma} \label{37}
 In case $($a$)$ of Cor. \ref{IIIposs} with $|I| \geq 2$ there are no
 type II vectors associated to elements of $\bar{c}^\perp$.
\end{lemma}

\begin{proof}
  Let $v = (-2^1, 1^k)$ and $w= (-1^1, 1^i, -1^j)$
give elements of $\bar{c}^{\perp}$ with $k \neq i,j$.
  Consider the face $\{\bar{x} : x_i + x_k=1, \; x_1 + x_j =-2\}
\cap \bar{c}^{\perp}$.
 Other than $v,w$
  the possible elements in this face come from $u=(-1^1, 1^k, -1^j)$ and $s=
  (-2^1, 1^i)$. As $d_1=4$, $J(\bar{v}, \bar{w}) \neq 0$, so $vw$ is not an
  edge and $u$ must be present. But $J(\bar{u}, \bar{w}) =
  \frac{1}{4}(1 -\frac{1}{d_1} -
  \frac{1}{d_j}) \neq 0$ since $d_1=4$, giving a contradiction. So
  $k=i$ or $j$ for every such $v,w$.

Hence if such a $w$ exists there are at most two type III
vectors. Now if $|I|=2$ and the type IIIs are $(-2,1,0,\cdots), (-2,0,1,\cdots),$
we cannot have $w = (-1,1,-1, \cdots)$ or $(-1,-1,1, \cdots)$ as then a
suitable affine combination of the above vectors give a type I vector.
(cf Lemmas \ref{Iperp}, \ref{coplane}).
So in fact no type II vectors give rise to elements of $\bar{c}^\perp$.
\end{proof}

\begin{lemma} \label{50}
The vectors $v=(-2,1,0, \cdots)$ and $w=(0,1,-1,-1,0, \cdots)$ are not
both associated to elements of $\bar{c}^\perp$, unless $(0,1,-2,0,\cdots)$
or $(0,1,0,-2,0,\cdots)$ is also.
\end{lemma}

\begin{proof}
Suppose $(0,1,-2,0,\cdots), (0,1,0,-2,0,\cdots)$ are absent.
Consider the face $\{\bar{x} : x_2 = 1, \; x_1 + x_3 + x_4 =-2\}
 \cap \bar{c}^{\perp}$.
 The other possible
elements of this face come from
 $t= (-1,1,-1,0,\cdots)$ and $y=(-1,1,0,-1,0,\cdots)$.
Both these must be present, as $J(\bar{v}, \bar{w}) \neq 0$. Applying
Cor \ref{orthogonal} to $wt,vt$ and $wy$ we obtain $(d_1,d_2,d_3,d_4)=(4,2,2,2)$.

Now we have equality (for $y,t$) in Eq.(\ref{est}), as both sides equal $15/16$.
We find that $\lambda = 1/2$, giving a contradiction again to Remark \ref{equal}.
\end{proof}

\medskip
Combining this with Lemma \ref{II/III} (and using Lemma \ref{AC}) yields:
\begin{corollary}  \label{uniqueIII}
If there is a unique type III vector $u=(-2^1,1^2)$
with $\bar{u}$ in $\bar{c}^\perp$,
then the type II vectors associated to elements of $\bar{c}^{\perp}$
 all have ``-1'' in place 1. Moreover $(-1^1,-1^2,1^i)$
cannot be present. Also, if $(-1^1,1^2, -1^i)$ is present for some $i \geq 3$
then $(d_1, d_2)=(4,2)$ or $(3,3)$ and the index $i$ is unique.
\end{corollary}

For the last assertion, observe that $(-1^1,1^2,-1^i)$ and the type III
vector are joined by an edge, so Cor \ref{orthogonal} shows the dimensions
are as stated. If we have two such type II for $i_0$ and $i_1$  then
Eq.(\ref{est}) implies $d_{i_0} + d_{i_1} \leq 4$. Hence since $K$ is connected,
$d_{i_0}= d_{i_1}=2$ and we have equality in Eq.(\ref{est}) with
$\lambda= \frac{1}{2}$, giving a contradiction.

\begin{lemma} \label{34}
Let the type III vectors be as in  Cor \ref{IIIposs}$($b$)$, i.e., they are
$(1^1, -2^a),  a \in I$. Assume that  $|I| \geq 2$.
If we have a type II vector $w=(1^1, -1^i, -1^j)$
with $\bar{w}$ in $\bar{c}^{\perp}$ then $i,j \in I$
\end{lemma}

\begin{proof}
  Suppose for a contradiction that $w = (1^1, -1^i, -1^j)$ is present
  (so $d_1 \neq 1$) and $(1^1, -2^j)$ absent (i.e. $j \notin I)$.
  Since $|I|\geq 2,$ we can consider $v = (1^1,-2^k)$ where $k \in I$
  (so $k \neq j$) and $k \neq i$.  Consider the face
 $\{\bar{x} : x_1=1, \; x_i
  + x_j + x_k = -2\} \cap \bar{c}^{\perp}$.
 As well as $v,w$ the possible elements of
  $\mathcal W$ in the face giving elements of $\bar{c}^\perp$ are
  $y=(1^1, -1^i, -1^k)$, $t= (1^1, -1^j, -1^k)$ and $u= (1^1, -2^i)$.
  As $d_1 \neq 1$, $vw$ is not an edge so $t$ is present. Now Cor
  \ref{orthogonal} applied to $vt$ and $tw$ gives $d_1 = d_j = 2$ and
  $d_k =4$.

Moreover, if $i \in I$ then $u$ is present, so the edge $wu$ gives $d_i=4$.
Thus we have shown that $d_a = 4$ for all $a \in I$.

Now considering $(1^1, -2^a)$ and $(1^1, -2^b)$ with $a,b \in I$, we see that
we have equality in Eq.(\ref{est}) (both sides equal $\frac{3}{4}$). In fact
$c$ is the average of these two vectors (i.e., $\lambda=\frac{1}{2}$), so as
in Remark \ref{equal} we have a contradiction.
\end{proof}

\begin{lemma} \label{55}
Let the type III vectors be as in Cor \ref{IIIposs}$($b$)$, i.e., they are
$(1^1, -2^a),  a \in I$. Assume that  $|I| \geq 3$. Then $d_1 = 1$.
\end{lemma}

\begin{proof}
Each pair $v,w$ of type III vectors gives an edge, and if $d_1 \neq 1$,
then we have $J(\bar{v}, \bar{w})>0$. By Theorem \ref{orthogonal}
all the midpoint vectors $(1^1, -1^a, -1^b)$ are present for $a, b \in I$.
Now Prop \ref{edgeCC} shows that $F_{\bar{v}}$ and $F_{\bar{w}}$ have
opposite signs, so we have a contradiction if $|I| \geq 3$.
\end{proof}

\medskip
Putting together our results so far, we obtain a description of the possibilities
for $ \bar{c}^\perp \cap \frac{1}{2}(d + {\mathcal W})$.

\begin{theorem} \label{class-1a}
Assume that $r \geq 3$ and $K$ is connected, and
that we are not in the situation of Thm \ref{Fano}. Up to permutation of the irreducible
summands, the following are the possible configurations of vectors in $\mathcal W$
associated to elements of $\frac{1}{2}(d + {\mathcal W}) \cap \bar{c}^\perp$.

\noindent{$($1$)$}  $\{(-2^1, 1^i), \ 2\leq i \leq m\}$ for fixed $m \geq 2$.
    There are no type II vectors, and $d_1=4$ if $m \geq 3$.

\noindent{$($2$)$}  $\{(1^1, -2^i), \  2\leq i \leq m\}$ for fixed
     $m \geq 3$ and $d_1=1.$ There are no type II vectors.

\noindent{$($3$)$}$($i$)$  $\{(1^1,-2^2), (1^1,-2^3), (-1^2,-1^3,1^i), \ 4 \leq i \leq m\}$
   with $d_1=1, \ d_2 = d_3=2$.

$($ii$)$ $\{(1^1, -1^2, -1^3), (1^1, -2^2), (1^1, -2^3), (-1^2, -1^3, 1^i), \ 4 \leq i \leq m\}$,
  \ $d_1 \neq 1, d_2=d_3=2$.

\noindent{$($4$)$} $\{(1,-2,0,0, \cdots), (1,0,-2, 0,\cdots), (1,-1,-1, 0,\cdots)\}$ with $d_1 \neq 1$.

\noindent{$($5$)$} A unique type III $(-2,1,0, \cdots)$. Possible type II vectors are

$($i$)$ $(-1,1,-1, 0, \cdots)$ with either $(d_1,d_2) = (4,2)$ or $(3,3)$; or

$($ii$)$ $\{(-1^1,1^3,-1^i), \ 4\leq i \leq m \}$ for fixed $m\leq r$ and  with $d_1=2$; or

$($iii$)$ $\{(-1^1,-1^3,1^i), \ 4 \leq i \leq m \}$ for fixed $m \leq r$ and with $d_1=2$.

\noindent{$($6$)$}  No type III vectors.   Possible type II vectors are

$($i$)$  $\{(-1^1,-1^2,1^i), \ 3\leq i \leq m\}$ for fixed $m\leq r$, with $d_1=d_2=2$ if $m \geq 4$; or

$($ii$)$ $\{(1^1, -1^2, -1^i), \ 3\leq i \leq m \}$ for fixed $m \leq r$, with $d_1=d_2=2$ if $m \geq 4$; or

$($iii$)$ $\{(1^1,-1^2,-1^3), (1^1,-1^2,-1^4), (1^1,-1^3,-1^4)\}$ with $d_1 = d_2 = d_3= d_4=2$.
\end{theorem}

\begin{proof}
Cor \ref{IIIposs} gives the possibilities for the type III vectors in
$\bar{c}^\perp$. If there are none then Cor \ref{IIposs} gives the
possibilities in (6). If there is a unique type III vector, then
Cor \ref{uniqueIII} and Cor \ref{IIposs} give us the cases listed in (5)
(or (1) with $m=2$ if there are no type II).
If we have two or more type III vectors with $-2$ in the same place then
Lemma \ref{37} shows we are in case (1).

If we have more than two type III vectors with  $1$ in the same place $a$,
then $d_a =1$ by Lemma \ref{55}. Remark \ref{IIrmk} then implies there
are no type II vectors and we are in case (2).

If we have exactly two type III vectors with $1$ in the same place, e.g.,
$(1, -2, 0, \cdots)$ and $(1, 0, -2, 0, \cdots)$, then the proof of
Lemma \ref{55} shows that  if the type II vector $(1, -1, -1, 0, \cdots)$
is absent we must have $d_1=1$. On the other hand, if $d_1=1$ we are,
by Remark \ref{IIrmk} and the connectedness of $K$, in case (2) or (3)(i).
If $d_1 \neq 1$, then by the above, Remark \ref{IIrmk}, and
Cor \ref{IIposs}, we are in case (3)(ii) or (4).

The statements about values of the $d_i$ follow from straightforward
applications of Cor \ref{orthogonal} to the obvious edges of
 ${\rm conv}(\frac{1}{2}(d +{\mathcal W})) \cap \bar{c}^\perp$.
\end{proof}

\begin{rmk} \label{suppl-1a}
The possibilities in Theorem \ref{class-1a} can be somewhat sharpened.
In cases (1), (2), and (3), $m$ cannot be $r$; in other words the maximum
number of vectors is not allowed. This follows easily from looking at
the system of equations expressing the nullity of $\bar{c}$, the orthogonality
of the vectors to $\bar{c}$ and the fact that the entries of $c$ sum up to $-1$.
Similarly, $r \neq 3$ in (5)(i) and $r \neq 4$ in (6)(iii).

When $m \geq 5$ in (5)(ii) or (5)(iii), the segment joining two type II vectors
is an edge, so Cor \ref{orthogonal} gives $d_3 =2$.
\end{rmk}

\section{\bf Adjacent (1B) vertices}

We now turn to (1B) vertices.
Let $\bar{\xi}, \bar{\xi}^{\prime}$ be adjacent (1B) vertices of $\Delta$.
Then there exist vertices $\bar{x}, \bar{x}^{\prime}$ of
${\rm conv}(\frac{1}{2}(d + {\mathcal W}))$ such that
$\bar{c},\bar{\xi},\bar{x}$ are collinear and $\bar{c},
\bar{\xi}^{\prime}, \bar{x}^{\prime}$ are collinear.  Moreover, there
exist null vectors $\bar{a}, \bar{a}^{\prime}$ such that
$\bar{x} = (\bar{a} + \bar{c})/2$ and $\bar{x}^{\prime}= (\bar{a}^{\prime} + \bar{c})/2$.
By Cor \ref{meet}, there must be an element $\bar{y}$ of
${\rm conv}(\frac{1}{2}(d + {\mathcal W}))$ on $\bar{a} \bar{a}^{\prime}$,
so $P^{-1}(\bar{\xi} \bar{\xi^{\prime}})$ contains the convex hull of
 $\bar{x}, \bar{x}^{\prime}, \bar{y}$ and hence is 2-dimensional.
As $\bar{\xi} \bar{\xi}^{\prime}$ is by assumption an edge of $\Delta$,
$P^{-1}(\bar{\xi} \bar{\xi}^{\prime})$ is a 2-dimensional face of
${\rm conv}(\frac{1}{2}(d + {\mathcal W}))$.

So we need to analyse the 2-dimensional faces of  conv$(\mathcal W)$
containing vertices $x, x^{\prime}$ such that
\begin{equation} \label{1Bcond}
x = (a + c)/2, \;\;\;  x^{\prime} = (a^{\prime} + c)/2, \;\;\;\
\bar{a}, \bar{a}^{\prime} \; {\rm null},
\end{equation}
and such that $c$ lies in the 2-dimensional plane defining this face.
The lines through $x,c$ (resp. $x^{\prime},c$) only meet ${\rm conv}({\mathcal W})$
at $x$ (resp. $x^{\prime}$).

Most 2-faces of conv$(\mathcal W)$ are triangular. We list below (up
to permutation of components) all the possible non-triangular faces.
For further details regarding how this listing is arrived at, see \cite{DW5}.
We emphasize that only the {\em full} faces are being listed, i.e.,
configurations formed by all the possible elements of $\mathcal W$ in
a given $2$-dimensional plane.  As the set of weight vectors for a
given principal orbit may be a subset of the full set of possible
weight vectors, these full faces may degenerate to subfaces or even
lower-dimensional faces (see Remark \ref{subshape}).

\medskip

\noindent{\it Listing convention:} In the interest of economy and clarity,
we make the convention that when we list vectors in $\mathcal W$ belonging
to a $2$-face we will use the freedom of permuting the summands to
place nonzero components of the vectors first and we will only put down the
minimum number of components necessary to specify the vectors.

\smallskip

\noindent{\bf Hexagons:}  There are $3$ possibilities.
\smallskip

\noindent{(H1)} This is the face in the plane $\{x_1 + x_2 + x_3 = -1; \ x_a=0, \
{\rm for} \ \ a >3 \}$.  Points of $\mathcal W$ are
$(-2^i, 1^j), (-1^i, 1^j, -1^k), (-1^i)$ where  $i,j,k \in \{1,2,3\}$.  The  type  III
vectors form the vertices of the hexagon.

\noindent{(H2)} The plane here is $\{x_1 + x_2 = -1, \;
x_3 + x_4=0, \; x_i=0 \;(i > 4)\}.$  Points of $\mathcal W$ are vertices
\[
u=(-2,1,0,0), \ v=(1,-2,0,0), \ y=(-1,0,1,-1), \ y^{\prime}=(0,-1,1,-1),
\]
\[
z=(-1,0,-1,1), \ z^{\prime}=(0,-1,-1,1),
\]
and the interior points
\[
 \alpha=(-1,0,0,0), \ \beta=(0,-1,0,0).
\]

\noindent{(H3)} The plane is $\{x_2=-1, \ x_1+x_3+x_4 =0, \ x_i=0 \;(i > 4)\}$. Points
of $\mathcal W$ are the vertices
\[
u=(-1,-1,1,0), \ v=(0,-1,1,-1), \ w=(1,-1,0,-1),
\]
\[
 x=(1,-1,-1,0), \ y=(0,-1,-1,1), \ z=(-1,-1,0,1)
\]
and the centre
\[
t=(0,-1,0,0).
\]

\noindent{\bf Square:} (S) with midpoint $t=(0, -1, 0, 0, 0)$  and vertices
$$ v=(-1,-1,1,0,0), \ u=(0,-1,0,1,-1),
$$
$$ s=(0,-1,0,-1,1), \ w=(1,-1,-1,0,0).$$


\noindent{\bf Trapezia}: We have vertices $v,u,s,w,t$ with $2v-s = 2u-w$
and $t = \frac{1}{2}(s+w),$ i.e., these are symmetric trapezia.  Below we
list the possible $v, u, s, w$.
\[
\begin{array}{|c|c|c|c|c|} \hline
     &       v          &        u         &         s        &  w    \\ \hline
(T1) &(-2,1,0,0) &(-2,0,1,0) &(0,0,-2,1) & (0,-2,0,1)  \\
(T2) &(-2,0,1,0) &(-2,1,0,0) &(0,-1,1,-1) & (0,1,-1,-1) \\
(T3) &(-1,-1,0,1) &(0,-1,-1,1) &(-2,1,0,0) & (0,1,-2,0)  \\
(T4) & (0,0,1,-1,-1) &(1,0,0,-1,-1) &(-2,1,0,0,0) &(0,1,-2,0,0) \\
(T5) &(-1,0,0,1,-1) &(0,0,-1,1,-1) & (-2,1,0,0,0) &(0,1,-2,0,0)  \\
(T6) &(1,-1,-1,0,0) &(1,-1,0,-1,0) & (0,0,-1,1,-1) &(0,0,1,-1,-1)  \\ \hline
\end{array}
\]
\centerline{Table 4: Possible trapezoidal faces}

Note that the configuration with vertices $(-1,-1,1,0,0),(-1,-1,0,1,0),
(0,0,1,-1,-1),$ and  $(0,0,-1,1,-1)$ is equivalent to (T6) under the composition of
a permutation and a $J$-isometric involution.


\smallskip

\noindent{\bf Parallelograms}: We have vertices $v, u, s, w$ with $v-u=s-w$.
\[
\begin{array}{|c|c|c|c|c|} \hline
     &  v          &       u       &        s   &     w   \\ \hline
(P1) &  (-2,1,0,0)  &  (-1,0,-1,1) &  (-2,0,1,0) & (-1,-1,0,1) \\
(P2) & (-2,1,0,0,0) & (-2,0,1,0,0) &( 0,1,0,-1,-1) & (0,0,1,-1,-1) \\
(P3) & (-2,1,0,0,0) & (-2,0,1,0,0) & (0,0,-1,-1,1) & (0,-1,0,-1,1) \\
(P4) & (-2,1,0,0) &  (-1,0,1,-1) & (-1,-1,0,1) & (0,-2,1,0)  \\
(P5) & (-2,1,0,0,0) &   (-1,0,0,1,-1) &  (-1,0,1,-1,0) &  (0,-1,1,0,-1) \\
(P6) & (-2,1,0,0,0) &   (-1,0,0,-1,1) &  (-1,0,-1,1,0) &  (0,-1,-1,0,1) \\
(P7) & (-2,1,0,0,0) &  (0,1,-1,-1,0) & (-1,0,0,1,-1) &(1,0,-1,0,-1) \\
(P8) & (1,-1,-1,0,0,0) & (0,0,0,1,-1,-1) & (1,-1,0,-1,0,0) & (0,0,1,0,-1,-1) \\
(P9) & (1,-1,-1,0,0,0) & (1,-1,0,0,-1,0) & (0,0,-1,1,0,-1) & (0,0,0,1,-1,-1)\\
(P10) & (1,-1,-1,0,0,0) & (1,0,0,0,-1,-1) &(0,-1,-1,1,0,0) &(0,0,0,1,-1,-1)\\
(P11) & (0,0,1,-1,-1,0) & (0,-1,0,-1,0,1) &( 1,0,0,0,-1,-1) & (1,-1,-1,0,0,0)\\
(P12) & (1,0,-1,0,-1) & (1,-1,-1,0,0) & ( 0,0,-1,1,-1) & (0,-1,-1,1,0)\\
(P13) & (-1,0,-1,0,1) & ( 0,0,-1,-1,1) & (0,-1,-1,1,0) & (1,-1,-1,0,0)\\
(P14) &(-1,0,1,0,-1) & (-1,-1,1,0,0) &(0,0,-1,1,-1) & (0,-1,-1,1,0) \\
(P15) & (-1,0,1,0,-1) & ( -1,-1,1,0,0) & (0,0,1,-1,-1) & (0,-1,1,-1,0)\\
(P16) & (-2,1,0,0) &  (0,1,-2,0) & (-1,0,1,-1) & (1,0,-1,-1) \\
(P17)  & (-2,1,0,0)  &(0,1,0,-2) &(-2,0,1,0) &(0,0,1,-2) \\ \hline
\end{array}
\]

\centerline{Table 5: Possible parallelogram faces}

\begin{rmk}
(P1), (P2), (P3), and (P17) are actually rectangles.
(P16) also includes the midpoints $y= (u+v)/2 = (-1,1,-1,0)$ and
$z= (s+w)/2 = (0,0,0,-1)$. The rectangle (P17) also includes the midpoints
$y=(u+v)/2 = (-1,1,0,-1)$ and  $z=(s+w)/2=(-1,0,1,-1).$
\end{rmk}

\begin{rmk} \label{subshape}
We must also consider subshapes of the above. Each symmetric
trapezium contains two parallelograms. The two rectangles
with midpoints (P17), (P16) will contain asymmetric trapezia.
(P17) also contains parallelograms and squares. (For (P16), note
that $s$ is present iff $w$ is.) Furthermore, there are numerous
subshapes of the hexagons. The regular hexagon (H3) contains
rectangles with midpoint (by omitting opposite pairs of vertices). Besides triangles,
the hexagon (H2) contains pentagons, rectangles and squares
(with midpoints), and  kite-shaped quadrilaterals (e.g. $y^{\prime} u z^{\prime} v)$.
For (H1) see the discussion before Theorem \ref{triangle}.
Finally, the triangle with midpoints of all sides (where the
vertices are the three type III vectors with $1$ in the same place)
contains a trapezium (by omitting one vertex) and hence parallelograms.
\end{rmk}

\begin{rmk} \label{nonface} We also note for future reference that
there are examples where we can have four or more coplanar elements
of $\mathcal W$ but the plane cannot be a face. These examples are not
of course relevant to the case of adjacent (1B) vertices, but some
will be relevant when we consider multiple vertices of type (2).
The examples which we will need in that context are the following
three trapezia

\[
\begin{array}{|c|c|c|c|c|} \hline
         &     v     &      u        &    s      &     w \\ \hline
(T^{*}1) & (0,1,-1,-1) & (1,0,-1,-1) & (-2,1,0,0) & (1,-2,0,0) \\
(T^{*}2) &(0,-1,1,-1) & (1,-1,0,-1) & (-2,1,0,0) & (0,1,-2,0)  \\
(T^{*}3) & (1,-1,-1,0) & (1,-1,0,-1) &(-1,0,-1,1) & (-1,0,1,-1) \\ \hline
\end{array}
\]

\centerline{Table 6: Further trapezia}

\smallskip

In (T*2),(T*3), as in (T1)-(T7), we have $2v-s = 2u-w$.
In these examples $t = \frac{1}{2}(s+w)$ may also be present.
In (T*1) we have $s-w = 3(v-u)$, and the vectors
$t=(2s+w)/3=(-1,0,0,0)$ and $r=(s+2w)/3=(0,-1,0,0)$ will also
be present.

As an example, we explain why the trapezium (T*2) can never be a face.
As $u$ is present in $\mathcal W$, so are $u^{\prime}=(-1,1,0,-1)$ and
$u^{\prime \prime}= (-1,-1,0,1)$. Now $(2 u^{\prime} + u^{\prime \prime})/3 = (2s +u)/3=
(-1, \frac{1}{3},0,-\frac{1}{3})$ is in the plane,
but $u^{\prime}$ is not, so this plane cannot give a face.
Similar arguments involving  $(1,0,-1,-1), (-1,0,0,0)$
(resp. $(-1,0,-1,1),(-1,0,1,-1)$) show (T*1) (resp. (T*3)) cannot
be faces.

These arguments also show several parallelograms cannot be faces,
 but these will not be relevant for our purposes.
\end{rmk}

We now begin to classify the possible $2$-faces
which arise from adjacent (1B) vertices. We shall repeatedly
use Prop \ref{E}, Cor \ref{meet}, and Lemma \ref{III44}. Let $E$ denote
the affine $2$-plane determined by the $2$-face being studied.

\begin{theorem} \label{parallel}
Suppose we have adjacent $($1B$)$ vertices corresponding to a parallelogram
face $vusw$ of ${\rm conv}(\mathcal W)$. So we have
$\bar{u}= (\bar{a}+\bar{c})/2$ and
$\bar{w}= (\bar{a}^{\prime} +\bar{c})/2$ for null $\bar{a}, \bar{a}^{\prime}$.
Suppose the vertices $v,u,s,w$ are the only elements of $\mathcal W$ in
the face. Then  $u,w$ are adjacent vertices of the parallelogram, and either

$($i$)$
 ${\mathcal C}\cap E = \{\bar{c},\bar{a},\bar{a}^{\prime}, \bar{e} \}$
 where $\bar{e}$ is null with
$v=(a+e)/2$ and $s = (a^{\prime} + e)/2$; or

$($ii$)$ $\bar{v},\bar{s} \in \mathcal C$ and $J(\bar{a},\bar{v}) =
 J(\bar{a}^{\prime}, \bar{s})= J(\bar{s},\bar{v})=0$.

Moreover, if none of $v,u,s,w$ is type I, then $($i$)$ cannot occur.
\end{theorem}

\begin{proof}
We may introduce coordinates in the 2-plane $E$ using the sides $sv$ and $sw$ to define
the coordinate axes. In this way we can speak of ``left" or ``right",
``up" or ``down". If we extend the sides of the parallelogram to
infinite lines, these lines divide the part of the plane outside the
parallelogram into 8 regions, and $\bar{c}$ must be in the interior
of one such region.

We first observe that if $\bar{c}$ is in one of the four regions
which only meet the parallelogram at a vertex, then $\bar{a} \bar{a}^{\prime}$
does not meet the parallelogram, contradicting Lemma \ref{meet}.

\noindent{(A)} Let $\bar{c}$ then lie in a region which meets the parallelogram in an edge.
Without loss of generality we may assume the edge is $uw$.
By Cor \ref{meet}, all elements of $\mathcal C \cap E$
lie on or between the rays from $\bar{c}$ through
$\bar{a}, \bar{a}^{\prime}$. Hence, by Lemma \ref{Jpos}, $J(\bar{b},\bar{c}) >0$
for all $\bar{b} \in {\mathcal C} \setminus \{ \bar{c} \}$. If $\bar{b}$ is a rightmost
element of $ ({\mathcal C}\cap E) \setminus \{ \bar{c} \}$, then as $\bar{b} +\bar{c}$
cannot be written in another way as a sum of two elements of $\mathcal C$,
we deduce from Prop \ref{E} that $\bar{b} + \bar{c} \in d + \mathcal W$.
So $\bar{b}$ is either  $\bar{a}$ or $\bar{a}^{\prime}$. All other elements
of ${\mathcal C} \cap E$ lie to the left of $\bar{a} \bar{a}^{\prime}$.
Note also that a rightmost element of
$({\mathcal C}\cap E) \setminus \{ \bar{c}, \bar{a}, \bar{a}^{\prime} \}$ satisfies
$b+c = a+a^{\prime}, 2v$ or $2s$.

\noindent{(B)} Next let $\bar{e} = 2\bar{v}- \bar{a}$. Observe that as well as
$\bar{v}=(\bar{a}+\bar{e})/2,$ we have $\bar{s}= (\bar{a}^{\prime} + \bar{e})/2$,
since $2\bar{v}-\bar{a} =
2(\bar{v}-\bar{u}) +\bar{c} = 2(\bar{s}-\bar{w}) +\bar{c} =2\bar{s} - \bar{a}^{\prime}$.

If $\bar{e} \in \mathcal C$, then it must be null, and the same argument as
above shows that no elements of $(\mathcal C \cap E) \setminus \{ \bar{e} \}$
lie to the left of $\bar{a},\bar{a}^{\prime}$, so we are in  case (i).
Now, Lemma \ref{Jpos} shows $J(\bar{h},\bar{k}) >0$ for all
$\bar{h} \neq \bar{k} \in \mathcal C \cap E$. If $v,u,s,w$ are all type
II/III, we see that $F_{\bar{c}},F_{\bar{e}}$ are of one sign
 and $F_{\bar{a}},F_{\bar{a}^{\prime}}$ the other sign.
But now the contributions from $\bar{a} + \bar{a}^{\prime}$ and
$\bar{c}+\bar{e}$ in the superpotential equation cannot cancel.

If $\bar{e} \notin \mathcal C$ then, as in the argument before Theorem
\ref{Vertex}, $\bar{s},\bar{v} \in \mathcal C$ and we are in case
(ii). Prop \ref{edgeCC} shows $\bar{v}, \bar{s}$ are orthogonal.
Moreover, note that the remark at the end of (A) shows that $v + c$ or $s+c$
is left of $a + a^{\prime}$.
\end{proof}

\begin{lemma} \label{parallel2}
In case $($ii$)$ of Theorem \ref{parallel}, we have
$J(\bar{v},\bar{v})=J(\bar{s},\bar{s})$.
\end{lemma}

\begin{proof}
As $\bar{c}$ and $\bar{a}=2\bar{u}-\bar{c}$ are both null, and
similarly $\bar{c}$ and $\bar{a}^{\prime}=2\bar{w}-\bar{c}$ are both
null, we deduce (cf Remark \ref{xc})
\begin{equation} \label{rel1}
J(\bar{u},\bar{u}) = J(\bar{u},\bar{c}) \;\;\; : \;\;\;
J(\bar{w},\bar{w}) = J(\bar{w},\bar{c}).
\end{equation}
We also have
\begin{equation} \label{rel2}
2J(\bar{u},\bar{v}) = J(\bar{c},\bar{v}) \;\;\; : \;\;\;
2J(\bar{w},\bar{s}) = J(\bar{c},\bar{s})
\end{equation}
from the orthogonality conditions on $\bar{a},\bar{v}$ and
$\bar{a}^{\prime},\bar{s}$.

Now
$J(\bar{s},\bar{s}) - J(\bar{v},\bar{v}) = J(\bar{s},\bar{s}) -
J(\bar{w}-\bar{u}-\bar{s}, \bar{w}-\bar{u}-\bar{s})$,
 which, on expanding out and
using the second relations of Eqs.(\ref{rel1}),(\ref{rel2}), becomes
$J(2\bar{u}-\bar{c},\bar{w}-\bar{s}) -J(\bar{u},\bar{u})$. Now
\[
J(2\bar{u}-\bar{c},\bar{w}-\bar{s}) -J(\bar{u},\bar{u}) =
J(2\bar{u}-\bar{c}, \bar{u}-\bar{v}) - J(\bar{u},\bar{u}) =
J(2\bar{u}-\bar{c},\bar{u}) - J(\bar{u},\bar{u}) =
J(\bar{u}-\bar{c},\bar{u})=0.
\]
We have used the first relations of Eqs.(\ref{rel2}), (\ref{rel1}) in the
second and fourth equalities.
\end{proof}

\begin{rmk}
We must also consider the case when the midpoint of one side or a
pair of opposite sides of the parallelogram face is in $\mathcal W$.
This can happen for (P16) and (P17). Note that $v,u,s,w$ are
type II/III in these cases.

In fact, the argument of Theorem \ref{parallel} is still valid if
one or both of the midpoints of $vu, sw$ is in $\mathcal W$ and
$c$ lies in the region to the right of $uw$ (or the left of $vs$).

Keeping $c$ in the region to the right of $uw$, we now need to
consider the case where one or both of the midpoints of $vs, uw$ is
in $\mathcal W$.  The conclusions (in \ref{parallel}(ii)) still hold
except that we no longer have $J(\bar{v}, \bar{s})=0$.

However, we  have to make slight modifications to the arguments
as $\frac{1}{2}(\bar{a} + \bar{a}^{\prime})$ may be in $\mathcal C \cap E$.
If $\bar{e} \in \mathcal C$, then, as $\bar{a} + \bar{a}^{\prime}$
is not in $d + \mathcal W$, the usual sign argument shows that the
terms in the superpotential equation summing up to
$\bar{a}+\bar{a}^{\prime}$ do not cancel, which is a contradiction.
So $\bar{e} \notin \mathcal C$ and our previous arguments hold
except for the use of Prop \ref{edgeCC}.

Note that we also have to consider the possibility that $a, a^{\prime},$
and $e$ lie on the line through $vs$. But now the midpoint of $uw$ must
be present and ${\mathcal C} \cap E = \{ \bar{c}, \bar{a}, \bar{a}^{\prime},
\frac{1}{2}(\bar{a}+\bar{a}^{\prime}) \}$, with $v+s=a + a^{\prime}.$
The usual sign argument then forces the midpoints of $uw$ and $vs$
to be present and of type I. Hence this special configuration cannot
occur in (P16) or (P17).

Lastly, since the proof of Lemma \ref{parallel2} makes no mention
of midpoints, it remains valid if midpoints are present.
\end{rmk}

\smallskip

The conditions of Theorem \ref{parallel} and Lemma \ref{parallel2},
together with the nullity of $\bar{a}, \bar{a}^{\prime}, \bar{c},$ put
very strong constraints on $vusw$ and the dimensions. In fact, one can
check that these constraints cannot be satisfied for any of our parallelograms
(including those of Remark \ref{subshape}) with one exception.
This is the rectangle $yy^{\prime}z^{\prime}z$ in (H2) with
$c=(-2, 1, 0, \cdots)$ and $\frac{1}{d_3}+ \frac{1}{d_4}=\frac{1}{d_1},$
which will be dealt with in Lemma \ref{subH2}.
We now give an example of how to apply the above conditions in a specific
case.

\begin{ex} \label{example-parallel}
Consider parallelogram (P8). The equation of the $2$-plane $E$ containing
the parallelogram is
\begin{equation} \label{c}
x_2 = -x_1, \; x_5=x_6,  \; x_2 + x_5=-1, \; x_1 + \cdots + x_6 =-1
 \end{equation}
and $x_i=0$ for $i > 6$.  As all vertices are type II/III, we must
be in case (ii) of Theorem \ref{parallel}.

(A) Take $c$ to face the side $uw$. Note that $vs$ and $uw$ have equation
$x_1 = 1, x_1=0$ respectively, so $c_1 < 0$. Also, the remarks at the end of
parts (A) and (B) in the proof of Theorem \ref{parallel} shows that
$c_1 > -\frac{1}{3}$, as $v+c$ or $s+c$ is left of $a + a^{\prime}$
so $1 + c_1  > -2c_1$

The condition $J(\bar{v},\bar{s})=0$ implies $d_1=d_2 =2$ and Lemma
\ref{parallel2} implies $d_3=d_4$. Eqs.(\ref{rel1}) and (\ref{rel2})
give four linear equations in $c_i$. Now $d_3=d_4$ and Eq.(\ref{rel2})
show $c_3 =c_4$, so the equations for the plane give
$ c = (\frac{1}{2}-c_4, \ -\frac{1}{2} + c_4, \ c_4, \ c_4, -\frac{1}{2} -c_4,
\ -\frac{1}{2} - c_4 ).$
Next $d_1 =d_2 =2$ and Eq.(\ref{rel2}) show
 $c_4 = 3d_4/(2d_4 +2)$ and $c_1 = (1-2d_4)/(2d_4 + 2)$.

But the condition $-\frac{1}{3} < c_1 < 0$ now implies $d_3 =d_4=1$,
and it follows that $c$ cannot be null.

(B) The argument if $c$ faces $vs$ is very similar. We have
$d_3=d_4$ and $d_5 =d_6 =2$, and the orthogonality equations imply
$c_3 = c_4$. So $c$ has the same form  as
in the second paragraph of (A) above. We find
$c_4 = -3d_4/(2d_4 + 2))$ and $c_1 = (1+4d_4)/(2d_4 + 2)$.
But we now have the inequality $1 < c_1 < \frac{4}{3},$ so again
$d_3=d_4=1$, violating nullity.

(C)   If $c$ faces $vu$ or $sw$ then we need $J(\bar{s}, \bar{w})=0$
(resp. $J(\bar{v}, \bar{u})=0$), which is impossible.
\end{ex}

\begin{ex} \label{example-square}
The example of the square (S) with midpoint can be treated in essentially
the same way as the parallelograms. By symmetry, we may assume that $c$
lies in the region that intersects $uw$. However, because
$\frac{1}{2}(a + a^{\prime})$ may now be the midpoint and hence in
$\mathcal W$, the configuration of Theorem \ref{parallel}(i) can occur,
even though all vertices are type II.  We have $\mathcal C \cap E = \{\bar{c},\bar{a},
\bar{a}^{\prime}, \bar{e} \}$ with $a= (-1,-1,1,1,-1,\cdots),
 \ a^{\prime}=(1,-1,-1,-1,1, \ldots)$, \ $c= (1,-1,-1,1,-1,\cdots),$
and $e= (-1,-1,1,-1,1,\cdots)$ with nullity condition
$\sum_{i=1}^{5} \frac{1}{d_i} =1$. We will be able to rule this case out in \S 7.
On the other hand, the configuration of Theorem \ref{parallel}(ii) cannot
occur, as one easily checks.
\end{ex}

\smallskip
Next assume that adjacent (1B) vertices in $\Delta^{\bar{c}}$ determine
a trapezium $vusw$ as shown in the diagram below:

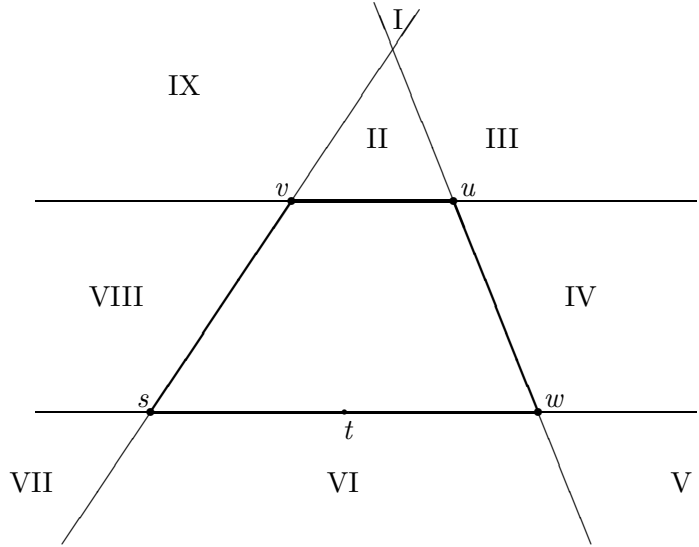
\begin{picture}(200, 210)(0, 30)
\put(90, 80){\line(1,0){250}}
\put(90, 160){\line(1,0){250}}
\put(100, 30){\line(2,3){135}}
\put(300, 30){\line(-2,5){82}}
\put(133.33,80){\circle*{3}}
\put(186.66,160){\circle*{3}}
\put(280,80){\circle*{3}}
\put(248,160){\circle*{3}}
\put(206.66,80){\circle*{2}}
\put(180, 162){\em v}
\put(250, 162){\em u}
\put(128, 82){\em s}
\put(282, 82){\em w}
\put(206, 70){\em t}
\put(225, 225){I}
\put(215, 180){II}
\put(260, 180){III}
\put(290,120){IV}
\put(330, 50){V}
\put(200,50){VI}
\put(80, 50){VII}
\put(110, 120){VIII}
\put(140, 200){IX}
\thicklines
\put(186.66,160){\line(1,0){61.34}}
\put(133.33,80){\line(1,0){146.67}}
\put(133.33,80){\line(2,3){53.33}}
\put(280,80){\line(-2,5){32}}
\end{picture}

{\noindent where $t$ is} the midpoint of $sw$ and $vu$ is parallel to $sw$. We assume
that $v, u, s, w \in {\mathcal W}$ but our conclusions hold whether or not $t$
lies in $\mathcal W$. We will now derive constraints on the $2$-face and
$E \cap {\mathcal C}$ resulting from having $c$ lie in one of the
regions shown above. For theoretical considerations, we need only treat
the cases where $c$ lies in regions I to VI. In practice, for an asymmetric
trapezium, we must consider $c$ lying in the remaining regions as well.
In the following we will adopt the convention that $\bar{a}, \bar{a}^{\prime}$
always denote null vectors in $\mathcal C$.

{\noindent (I)} $c$ {\it in region I}: \; This is impossible because then $\bar{s}=
\frac{1}{2}(\bar{c}+\bar{a})$ and $\bar{w}=\frac{1}{2}(\bar{c}+\bar{a}^{\prime})$
for some $\bar{a}, \bar{a}^{\prime},$ and so  $\bar{a}\bar{a}^{\prime}$
would not intersect ${\rm conv}(\frac{1}{2}(d+{\mathcal W})$,
a contradiction to Cor \ref{meet}.

{\noindent (II)} $c$ {\it in region II}: \; Then $\bar{v}=\frac{1}{2}(\bar{c}+\bar{a}), \bar{u}=
\frac{1}{2}(\bar{c}+\bar{a}^{\prime})$ for some $\bar{a}, \bar{a}^{\prime}$.
We get a contradiction to Cor \ref{meet} if $\bar{a}, \bar{a}^{\prime}$ lie
below the line $sw$. They also cannot lie on the line $sw$ since the argument
in (A) in the proof of Theorem \ref{parallel} and Cor \ref{meet} imply that
${\mathcal C} \cap E =\{\bar{c}, \bar{a}, \bar{a}^{\prime} \}$, and the
terms corresponding to $\bar{s}, \bar{w}$  in the
superpotential equation would be unaccounted for.

Let $e=2s-a, \; e^{\prime} = 2w-a^{\prime}$. These points lie in region VI, and
since we have a trapezium, $e \neq e^{\prime}$. We may now apply Theorem \ref{Vertex}
to $\bar{a}$ and $\bar{a}^{\prime}$ to obtain the possibilities:

(i) $\bar{s}, \; \bar{w} \in {\mathcal C}; \; J(\bar{a},{s})=0=J(\bar{a}^{\prime}, \bar{w})$,

(ii) $\bar{s} \in {\mathcal C}, \; J(\bar{a}, \bar{s}) = 0; \; w \notin {\mathcal C},$
     \; $\bar{e}^{\prime} \in {\mathcal C}$ is null, $J(\bar{e}^{\prime}, \bar{s}) = 0$,

(iii) $\bar{w} \in {\mathcal C}, \; J(\bar{w}, \bar{a}^{\prime}) = 0; \; \bar{s} \notin {\mathcal C},$
     $\bar{e} \in {\mathcal C}$ is null, $J(\bar{e}, \bar{w})=0$.

{\noindent Note that} the last condition in (ii) (resp. (iii)) results from
applying Theorem \ref{Vertex} to $\bar{e}^{\prime}$ (resp. $\bar{e}$).

{\noindent (III)} $c$ {\it in region III}: \; We have $\bar{v}=\frac{1}{2}(\bar{c}+\bar{a}),
  \bar{w}=\frac{1}{2}(\bar{c}+\bar{a}^{\prime})$ for some $\bar{a}, \bar{a}^{\prime}$
  lying respectively in regions VIII and VI (in view of Cor \ref{meet}).
  Applying Theorem \ref{Vertex} we obtain the possibilities:

(i) $\bar{s} \in {\mathcal C}, \; J(\bar{a}, \bar{s})=0=J(\bar{a}^{\prime}, \bar{s})$,

(ii) $\bar{s} \notin {\mathcal C}, \; 2\bar{s} = \bar{a}+\bar{a}^{\prime}$
     (which implies $c+s = v+w$).

{\noindent (IV)} $c$ {\it in region IV}: \; We have $\bar{u}=\frac{1}{2}(\bar{c}+\bar{a}),
\; \bar{w}=\frac{1}{2}(\bar{c}+\bar{a}^{\prime})$ for some $\bar{a}, \bar{a}^{\prime}
\in {\mathcal C} \cap E$.

 If $a$ lies in region IX, then Cor \ref{meet} implies that $\bar{a}^{\prime}$
lies in region VI. Applying Theorem \ref{Vertex} to $\bar{a}$ and $\bar{a}^{\prime}$
we obtain the possibilities:

(i) $\bar{s} \in {\mathcal C}, \; J(\bar{a},\bar{s})=0=J(\bar{a}^{\prime}, \bar{s})$,

(ii) $2s=a + a^{\prime}$, i.e., $c+s=u+w$.

If $a$ lies on the line $sv$, then we may apply Theorem \ref{Vertex} to
$\bar{a}^{\prime}$. We cannot have $2\bar{s}=\bar{a}^{\prime}+\bar{e}^{\prime}$
with $\bar{e}^{\prime} \in {\mathcal C}$ and null, otherwise $\bar{a}\bar{e}^{\prime}$
would not intersect ${\rm conv}(\frac{1}{2}(d+{\mathcal W}))$. So we have

(iii) $\bar{s} \in {\mathcal C}$ and $J(\bar{s}, \bar{a}^{\prime})=0$.

If $a$ lies in region II, then $\bar{a}^{\prime}$ lies in region VI.
Let $\bar{e}=2\bar{v}-\bar{a}$ and $\bar{e}^{\prime} = 2\bar{s}-\bar{a}^{\prime}$.
As we have a trapezium, $\bar{e} \neq \bar{e}^{\prime}$. Now $e$ lies
in region VII or VIII while $e^{\prime}$ lies in region VIII or IX,
so by Cor \ref{meet} $\bar{e}$ and ${\bar e}^{\prime}$ cannot both lie
in $\mathcal C$ and hence be null. Theorem \ref{Vertex} now gives
the possibilities:

(iv) $\bar{v}, \; \bar{s} \in {\mathcal C}, \; J(\bar{a}, \bar{v})=0=J(\bar{a}^{\prime}, \bar{s})$,
   (and by Prop \ref{edgeCC} $J(\bar{v}, \bar{s})=0$),

(v) $\bar{v} \in {\mathcal C}, \; J(\bar{a}, \bar{v})=0,$ \; $\bar{e}^{\prime} \in {\mathcal C}$
 is null, and $J(\bar{e}^{\prime}, \bar{v})=0$,

(vi) $\bar{s} \in {\mathcal C}, \; J(\bar{a}^{\prime}, \bar{s})=0$, \; $\bar{e} \in {\mathcal C}$
is null, and $J(\bar{e}, \bar{s})=0$.

{\noindent (V)} $c$ {\it in region V}: \; We have $\bar{u}=\frac{1}{2}(\bar{c}+\bar{a}^{\prime}),
   \bar{s}=\frac{1}{2}(\bar{c}+\bar{a})$ for some $\bar{a}, \bar{a}^{\prime}$
   lying respectively in regions VIII and II (by Cor \ref{meet}). Theorem
   \ref{Vertex} now gives the possibilities:

(i) $\bar{v} \in {\mathcal C}, \; J(\bar{a}, \bar{v})=0=J(\bar{a}^{\prime}, \bar{v})$,

(ii) $\bar{v} \notin {\mathcal C}, \; 2\bar{v}=\bar{a} + \bar{a}^{\prime}$
      (which implies $c+v = u+s$).

{\noindent (VI)} $c$ {\it in region VI}: \; We have $\bar{s}=\frac{1}{2}(\bar{c}+\bar{a}),
 \; \bar{w}=\frac{1}{2}(\bar{c}+\bar{a}^{\prime})$ for some $\bar{a}, \bar{a}^{\prime}$
  lying respectively in regions VIII and IV (by Cor \ref{meet}).
  (To rule out $\bar{a}, \bar{a}^{\prime}$ lying in the line $vu$, we
  proceed as in case (II), except that when $t \in {\mathcal W}$, we conclude
  instead that ${\mathcal C} \cap E =\{\bar{c}, \bar{a}, \bar{a}^{\prime},
  \frac{1}{2}(\bar{a}+\bar{a}^{\prime}) \}$. One can still check that $\bar{v},
  \bar{u}$ cannot be both accounted for.) Now let $\bar{e}=2\bar{v}-\bar{a}$ and
  $\bar{e}^{\prime} = 2\bar{u}-\bar{a}^{\prime}$. Again, having a trapezium means
  $\bar{e} \neq \bar{e}^{\prime}$ and Theorem \ref{Vertex} now gives the possibilities:

(i) $\bar{u}, \bar{v} \in {\mathcal C}, \; J(\bar{a}, \bar{v})=0=J(\bar{a}^{\prime}, \bar{u})$,
    \; (and $J(\bar{u},\bar{v})=0$ by Prop \ref{edgeCC}),

(ii) $\bar{v} \in {\mathcal C}, \; J(\bar{a}, \bar{v})=0, \; \bar{u} \notin {\mathcal C}$,
   $\bar{e}^{\prime} \in {\mathcal C} $ is null,

(iii) $\bar{u} \in {\mathcal C}, \; J(\bar{a}^{\prime}, \bar{u})=0, \; \bar{v} \notin {\mathcal C}$,
   $\bar{e} \in {\mathcal C} $ is null.

\begin{rmk} \label{spineq}
We mention a useful inequality which holds in (II) and (VI) above,
as well as in parallelogram faces with the same configuration
(cf Example \ref{example-parallel}(A)).

Let us consider (II), where we choose in $E$ coordinates such that
the first coordinate axis is parallel to $\bar{s}\bar{w}$
(assumed to be horizontal) and the second coordinate axis is arbitrary,
with the second coordinate increasing as we go up. As in (A) in the
proof of Theorem \ref{parallel}, all points in
$({\mathcal C} \cap E)\setminus \{\bar{c}, \bar{a}, \bar{a}^{\prime}\}$
must lie below the line $\bar{a}\bar{a}^{\prime}$. Let $\bar{b}$ be a
point among these with largest second coordinate. Since we have seen
above that either $\bar{s}$ or $\bar{w}$ lies in ${\mathcal C} \cap E$,
we have $s_2 \leq b_2$.  Furthermore,  as $\bar{b}+\bar{c}$ cannot
lie in $d+ {\mathcal W}$ it must be balanced by sums of elements
in ${\mathcal C} \cap E$, with the limiting configuration given by
$\bar{a}+\bar{a}^{\prime}$. So we have
$\frac{1}{2}(b_2 +c_2) \leq a_2 =a_2^{\prime} = 2v_2 -c_2$.
Combining the two inequalities we get $ 3c_2 \leq 4v_2 -s_2.$

Equality in the above holds iff $\bar{b}$ lies in $\bar{s}\bar{w}$
and $\bar{b}+\bar{c}=\bar{a}+\bar{a}^{\prime}$. In particular,
$\bar{b}$ is unique, so in II(i), the inequality above is strict.

Note that we only need $\bar{v}\bar{u}$ and $\bar{s}\bar{w}$ to
be parallel and the presence or absence of $t$ in $\mathcal W$
is immaterial. Hence in Theorem \ref{parallel}(ii) we also have
an analogous strict inequality, which we have already used, e.g.,
in (B) of Example \ref{example-parallel}. (For a parallelogram,
there may be midpoints on the pair of non-horizontal sides lying
in $\frac{1}{2}(d + {\mathcal W})$, but $\frac{1}{2}(\bar{b}+\bar{c})$
can never equal these midpoints, so we still get the inequality
we want.)

For the configuration in (VI), we still have an analogous
inequality, but since $\frac{1}{2}(\bar{a}+\bar{a}^{\prime}) \in {\mathcal C}$,
we lose uniqueness of $\bar{b}$ and hence the strict inequality.

We will also have occasion to apply the above analysis to appropriate trapezoidal
regions in  hexagon (H3).
\end{rmk}

\smallskip

The method described above together with Remark \ref{spineq} can now
be used to rule out the trapezia (T1)-(T6) as well as those mentioned in
Remark \ref{subshape}.

\begin{ex} \label{example-trapezium}
For the trapezium (T3), the vectors
$v, u, s, w$ are given in Table 4, and lie in the $2$-plane
$\{x_1+x_2+x_3+x_4=-1, x_2+2x_4 =1\}$. $vu$ is given by $x_4 =1$
while $sw$ is given by $x_4 =0$. $sv$ is given by $x_3 =0$
and $wu$ is given by $x_1 =0$. The vector $c$ that we are looking for
has the form $(-c_3 +c_4 -2, 1-2c_4, c_3, c_4)$. Since the trapezium
is symmetric, an explicit symmetry being induced by interchanging $x_1$
and $x_3$, we need only consider $c$ lying in regions II-VI.

(A) If $c$ lies in region III, then $c_1>0, \; c_4 >1$. Since $a=2v-c$,
we obtain $a=(c_3 -c_4, -3+2c_4, -c_3, 2-c_4)$. Similarly,
 $a^{\prime}=(c_3 -c_4+2, 1+2c_4, -4-c_3, -c_4)$.
If we are in case (ii), then $c=v+w-s=(1, -1, -2, 1)$, which violates
$c_4 > 1$. So we must be in case (i).

It follows from $J(\bar{a}, \bar{s})=0=J(\bar{a}^{\prime}, \bar{s})$ that
$d_1+3=-2c_3 + 4c_4$ and $J(\bar{w}, \bar{s})=J(\bar{v},\bar{s})$.
The second equality implies that $d_1=d_2$. Using this together with
the first equality and the null condition for $\bar{a}^{\prime}$
(in the form $J(\bar{w}, \bar{w})=J(\bar{w}, \bar{c})$, see Remark \ref{xc})
 we get $c_4 =d_1(d_1 -1)/(4d_1+2d_3)$. Since $c_4 >1$, we have $d_1(d_1 -5)>2d_3$,
so $d_1 >5$. But by Remark \ref{tangent}, $J(\bar{s}, \bar{s})< 0$,
which gives $d_1 < 5$ (since $d_1=d_2$), a contradiction.

(B) Let $c$ lie in region IV, so that $c_1>0, \; 0<c_4 <1$.
We obtain $a=(2+c_3 -c_4, 2c_4 -3,-2-c_3, 2-c_4)$ and
$a^{\prime}=(2+c_3 -c_4, 1+2c_4, -4-c_3, -c_4)$.  We claim that
$a_3 > 0$, so that $a$ lies in region IX. To see this, we solve for
$c_3, c_4$ using the null conditions $J(\bar{u}, \bar{u})=J(\bar{c}, \bar{u})$
and $J(\bar{w}, \bar{w})=J(\bar{c}, \bar{w})$ for $\bar{a}, \bar{a}^{\prime}$
respectively. We obtain $c_4 = (d_2 d_3 +2d_3 d_4 -d_2 d_4)/(d_3 (d_2+3d_4))$
and $a_3 = -2-c_3 =(d_2 d_3+2d_3 d_4 -d_2 d_4)/(d_2(d_2  +3d_4))$. Since
$c_4 > 0$ we obtain our claim.

Since $a$ lies in region IX, we first check if (ii) holds.
In this case, $c=(2, -1, -3, 1)$ which contradicts $c_4 <1$.
The equations in (i) together imply the contradiction $0=-4/d_2$.

(C) Suppose $c$ lies in region V, so that $c_1 > 0, \; c_4 <0$.
We obtain $a=(c_3 -c_4 -2, 1+2c_4, -c_3, -c_4)$ and
$a^{\prime}=(c_3 -c_4 +2, 2c_4 -3, -2-c_3, 2-c_4)$. If (ii) holds
then $c=(-1, 1, -1, 0)$ and this contradicts $c_1 > 0$.
Hence (i) must hold.

By Remark \ref{xc},  the null condition for $\bar{a}$ is
$J(\bar{s}, \bar{s})=J(\bar{s}, \bar{c})$, which is
$\frac{c_3}{d_1}-(\frac{1}{d_1}+\frac{1}{d_2}) c_4 =0$.
The two equations in (i) imply $J(\bar{u}, \bar{v})=J(\bar{s}, \bar{v})$
and $J(\bar{v}, \bar{v})< 0$, which in turn give $d_1 =2$. Using this,
the null condition for $\bar{a}$, and $J(\bar{a}^{\prime}, \bar{v})=0$
we obtain $c_4 = \frac{-1}{1+d_2}$ and $c_3= -\frac{d_2 +2}{d_2(d_2 +1)}$.
But $c_1=c_4-c_3 -2 > 0$, which simplifies to $1>d_2(d_2 +1)$, a contradiction.

(D) Let $c$ lie now in region II. Then $c_1 < 0, \; c_3 < 0, \; 1<c_4 \leq \frac{4}{3}$
where the last upper bound comes from the inequality in Remark \ref{spineq}.
We obtain $a=(c_3 -c_4, 2c_4 -3, -c_3, 2-c_4)$
and $a^{\prime}=(2+c_3 -c_4, 2c_4 -3, -2-c_3, 2-c_4)$. The null conditions
for $\bar{a}, \bar{a}^{\prime}$ then give
$$ c_3 = -\frac{2d_1 d_4 + d_1 d_2 -d_2 d_4}{(d_1 +d_3)(d_2 +2d_4)-d_2 d_4},
\ \ \ \ \ c_4=\frac{(d_1 +d_3)(d_2 + 2d_4)}{(d_1+d_3)(d_2+2d_4)-d_2 d_4}. $$

Suppose we are in case (i). The two equations and the above values of
$c_3, c_4$ combine to give $(d_1 -d_3)((d_1 +d_3)(d_2 +2d_4)-d_2 d_4) =0$.
However, the upper bound $c_4 \leq \frac{4}{3}$ translates into
$(d_1+d_3)(d_2 +2d_4) \geq 4d_2d_4$. So the second factor is positive
and we have $d_1=d_3$. Putting this information into the equation
$J(\bar{a}, \bar{s})=0$, we get
$$d_1 d_2 d_4(d_2 + 15)=2d_1^2 d_2(d_2 +1)-6d_1d_2^2 +2d_4(2d_1^2 d_2 +2d_1^2 +d_2^2).$$
By Remark \ref{tangent} we also have $J(\bar{s}, \bar{s}) < 0$, i.e.,
$1<\frac{4}{d_1} + \frac{1}{d_2}$, so either $d_2 =1$ or $d_1 < 8$.
Substituting these values into the equation above and using $c_4 \leq \frac{4}{3}$
we obtain in each instance a contradiction.

If we are in case (ii), then by adding the equations $J(\bar{a},\bar{s})=0$
and $2J(\bar{w}, \bar{s})-J(\bar{a}^{\prime}, \bar{s})=0$ (equivalent to
$J(\bar{e}^{\prime}, \bar{s})=0$), we obtain $1=\frac{2}{d_1}+\frac{1}{d_2}$.
Hence $(d_1, d_2)=(4,2)$ or $(3,3)$. One then checks that these values are
incompatible with the null condition for $\bar{e}^{\prime}$, $J(\bar{a}, \bar{s})=0$,
 and the bound $c_3 < 0$.

An analogous argument works to eliminate case (iii), where we now need
the bound $c_4 \leq \frac{4}{3}$ instead.

(E) Lastly suppose $c$ lies in region VI, so $c_1, c_3 < 0$ and $-\frac{1}{3}
\leq c_4 < 0$, where the lower bound for $c_4$ results from Remark \ref{spineq}.
We have $a=(c_3 -c_4 -2, 1+2c_4, -c_3, -c_4)$ and $a^{\prime}=(2+c_3-c_4, 1+2c_4,
-4-c_3, -c_4)$. Using the null conditions for $\bar{a}, \bar{a}^{\prime}$,
we obtain
$$c_1=-\frac{2(d_2+d_3)}{d_1+d_2+d_3}, \ \ c_2=  \frac{d_1+5d_2 +d_3}{d_1+d_2+d_3}, \ \
c_3 = \frac{-2(d_1 + d_2)}{d_1 + d_2 +d_3}, \ \  c_4 = \frac{-2d_2}{d_1 + d_2 +d_3}.$$

If we are in case (i), $J(\bar{u}, \bar{v})=0$ gives $d_2=d_4=2$. The
other two equations and the above values of $c_3, c_4$ then give
$3(d_1 + d_3 +2)(d_1 + d_3 -4) = 4(3d_1 +3d_3 -2)$. The lower bound
$-\frac{1}{3} \leq c_4$ becomes $d_1 +d_2 +d_3 \geq 6d_2$.
Using this inequality in the above Diophantine relation leads to
a contradiction. (Alternatively, observe the relation is a quadratic in
$d_1 + d_3$ with no rational roots).

For case (ii), using the two equations and the above values for
$c_3, c_4$, we arrive at the relation
$$(d_1+d_2+d_3)((d_1 -5)d_2 d_4+ d_1 d_4 + d_2 d_3 +2d_3 d_4)=
   2d_2(d_1 d_2 +2d_1 d_4 -d_2 d_4 +d_2 d_3 +2d_3 d_4).$$
Using the lower bound $-\frac{1}{3} \leq c_4$ in the above relation
we see that  $d_1 \leq 3$. By direct substitution, we further
obtain $d_1 \neq 3$. Finally, if $d_1 =2$, the null condition
for $\bar{c}$ gives $1 > \frac{1}{2}c_1^2$ and so  $d_2 + d_3 \leq 4$.
The lower bound on $c_4$ now implies $d_2 =1$. Since $c_2 > 1$,
the null condition for $\bar{c}$ is violated.

Case (iii) reduces to case (ii) upon interchanging the first and
third summands. Therefore, the trapezium (T3) has been eliminated.
\end{ex}

\smallskip

We discuss next the  hexagons (H1)-(H3). As the three cases
are similar, we will focus on (H3) and refer to the following
(schematic) diagram:

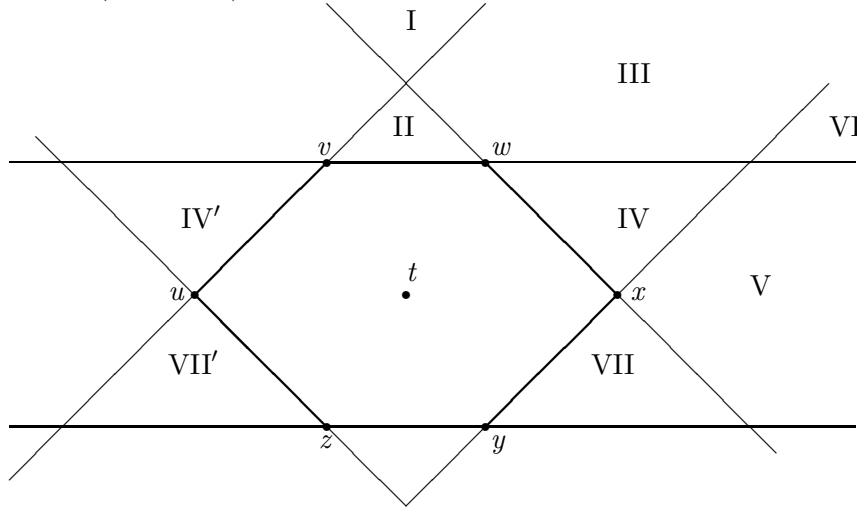
\begin{picture}(300, 200)(0, 10)
\put(80, 50){\line(1, 0){320}}
\put(80, 150){\line(1, 0){320}}
\put(230, 20){\line(1,1){160}}
\put(230, 20){\line(-1,1){140}}
\put(80, 30){\line(1,1){180}}
\put(370, 40){\line(-1, 1){170}}
\put(196, 153){\it v}
\put(262, 153){\it w}
\put(315, 98){\it x}
\put(140, 98){\it u}
\put(197, 42){\it z}
\put(262, 42){\it y}
\put(230, 105){\it t}
\put(230, 100){\circle*{3}}
\put(150, 100){\circle*{3}}
\put(200, 150){\circle*{3}}
\put(260, 150){\circle*{3}}
\put(310, 100){\circle*{3}}
\put(260, 50){\circle*{3}}
\put(200, 50){\circle*{3}}
\put(230, 200){I}
\put(225, 160){II}
\put(310, 180){III}
\put(310, 125){IV}
\put(145, 125){${\rm IV}^{\prime}$}
\put(390, 160){VI}
\put(360, 100){V}
\put(140, 70){${\rm VII}^{\prime}$}
\put(300, 70){VII}
\thicklines
\put(200, 50){\line(1,0){60}}
\put(150, 100){\line(1,1){50}}
\put(260, 50){\line(1,1){50}}
\put(200, 150){\line(1,0){60}}
\put(200,50){\line(-1,1){50}}
\put(310, 100){\line(-1,1){50}}
\end{picture}

\begin{ex} \label{example-hexagon}
The hexagon (H3) lies in the $2$-plane given by $\{x_2=-1, \; x_1+x_3+x_4 =0\}$.
So $c$ has the form $(-c_3-c_4, -1, c_3, c_4)$.
The lines $vw$ and $zy$ are given respectively by $x_1+x_3=1$
and $x_1 +x_3 =-1$. Similarly, the lines $uv$ and $yx$ are given by
$x_3 =1$ and $x_3 =-1$ respectively. The lines $uz$ and $wx$ are given
by $x_1 =-1$ and $x_1 =1$ respectively.

Interchanging $x_1$ and $x_3$ induces the reflection about
the perpendicular bisector of $vw$, while $(x_1, x_2, x_3, x_4) \mapsto
(-x_3, x_2, -x_1, -x_4)$ induces the reflection about $ux$.
These symmetries reduce our consideration to those $c$
lying in regions I-VI. Moreover, (H3) is actually a regular hexagon.
The symmetry $(x_1, x_2, x_3, x_4) \mapsto (-x_4, x_2, -x_3, -x_1)$
induces the reflection about $zw$, which swaps region II with region
IV and region I with region VI. Finally, the symmetry
$(x_1, x_2, x_3, x_4) \mapsto (-x_3, x_2, -x_4, -x_1)$
induces the rotation in $E$ about $t$ taking $x$ to $w$, and
maps region  V to region III. Therefore, we need only consider
$c$ lying in regions I, II, and V.

In the discussion below we again adopt the convention that
$\bar{a}, \bar{a}^{\prime}$ always denote null vectors in $\mathcal C$.

If $c$ lies in region I, then $\bar{u}=\frac{1}{2}(\bar{c}+\bar{a}),$
$\bar{x}=\frac{1}{2}(\bar{c}+\bar{a}^{\prime})$ for some $\bar{a}, \bar{a}^{\prime}$,
and we immediately see that $\bar{a}\bar{a}^{\prime}$ cannot meet
${\rm conv}(\frac{1}{2}(d+{\mathcal W}))$, a contradiction to Cor \ref{meet}.

\smallskip

{\noindent $c$ } {\em lying in region II}:

We have $c_1, c_3, <1$ and $c_1 + c_3 > 1$. The assumption of adjacent
(1B) vertices means that $\bar{v}=\frac{1}{2}(\bar{c}+\bar{a})$
and $\bar{w}=\frac{1}{2}(\bar{c}+\bar{a}^{\prime})$ for some
$\bar{a}, \bar{a}^{\prime} \in E \cap{\mathcal C}$. Hence
$a=(c_3 + c_4, -1, 2-c_3, -2-c_4)$ and $a^{\prime}=(2+c_3+c_4, -1, -c_3, -2-c_4)$.
One checks easily that $\bar{a}^{\prime}$ lies in region IV and $\bar{a}$
lies in region ${\rm IV}^{\prime}$. Moreover, the null conditions for
these vectors yield
$$c_1=\frac{d_3+d_4}{d_1+d_3+d_4}, \ \ c_3=\frac{d_1+d_4}{d_1+d_3+d_4}, \ \
   c_4=-\frac{d_1+d_3+2d_4}{d_1+d_3+d_4}. $$

Let $e:=2u-a$ and $e^{\prime}:=2x-a^{\prime}.$ These lie respectively
in regions ${\rm VII}^{\prime}$ and VII. We can now apply Theorem \ref{Vertex}
to $\bar{a}$ and $\bar{a}^{\prime}$ to obtain the following possibilities:

(i) $\bar{u}, \bar{x} \in {\mathcal C}$ and
     $J(\bar{a}, \bar{u})=0=J(\bar{a}^{\prime}, \bar{x})$;

(ii) $\bar{u} \in {\mathcal C}, \; J(\bar{a}, \bar{u})=0$,
      $\bar{x} \notin {\mathcal C}$, \; $\bar{e}^{\prime} \in {\mathcal C}$
      is null;

(iii) $\bar{x} \in {\mathcal C}, \; J(\bar{a}^{\prime}, \bar{x})=0$,
      $\bar{u} \notin {\mathcal C}, \; \bar{e} \in {\mathcal C}$ is null;

 (iv) $\bar{u}, \bar{x} \notin {\mathcal C}$, \; $\bar{e}, \bar{e}^{\prime}$
      are both null.

We can eliminate (i)-(iii) by noting that the two equations in each case
together with the values of $c_3, c_4$ above imply that $1=\frac{1}{d_1}+
\frac{1}{d_2}+\frac{1}{d_3}$. Using this relation (and the values of
$c_3, c_4$) in the null condition for $\bar{c}$ then leads to a contradiction.

For case (iv) we can again apply Theorem \ref{Vertex} to the null
vertices $\bar{e}$ and $\bar{e}^{\prime}$. The conditions
$J(\bar{e}, \bar{z})=0$ and $J(\bar{e}^{\prime}, \bar{y})=0$
lead, as above, to $1=\frac{1}{d_1}+\frac{1}{d_2}+\frac{1}{d_4}$
and $1 = \frac{1}{d_2} + \frac{1}{d_3} + \frac{1}{d_4}$ respectively.
Using this in the null condition for $\bar{c}$ again leads to a
contradiction. Hence $\bar{z}, \bar{y} \notin {\mathcal C}$
and $\bar{q}:=2\bar{z}-\bar{e}$ and $\bar{q}^{\prime}:=2\bar{y}-\bar{e}^{\prime}$
are null vectors in $E \cap {\mathcal C}$. In fact we now find that
$q=q^{\prime}$, so $caeqe^{\prime}a^{\prime}$ is a hexagon circumscribing
(H3).

Let us consider the pair of null vertices $\bar{c}, \bar{q}$.
We apply the argument in (A) of the proof of Theorem \ref{parallel}
to the wedge with vertex $\bar{c}$ bounded by the rays $\bar{c}\bar{a}$
and $\bar{c}\bar{a}^{\prime}$. All elements of $({\mathcal C} \cap E)
\setminus \{\bar{a}, \bar{a}^{\prime}, \bar{c}\}$ lie below
the line $\bar{a}\bar{a}^{\prime}$. Let $\bar{b}$ be a highest
(with respect to $x_1+x_3$) element among these. Since $\bar{e} \in
{\mathcal C}$, $b_1+b_3 > -1$ and so $\bar{c}+\bar{b}$ cannot
equal $2\bar{u}, 2\bar{t}, 2\bar{x}$. Hence $\bar{c}+\bar{b}=
\bar{a}+\bar{a}^{\prime},$ and we compute that
$b_1+ b_3=\frac{d_1+d_3-2d_4}{d_1+d_3+d_4}.$ The analogous
argument applied to the wedge bounded by the rays $\bar{q}\bar{e}$
and $\bar{q}\bar{e}^{\prime}$ gives a lowest element $\bar{b}^{\prime}$
of $({\mathcal C} \cap E) \setminus \{\bar{q}, \bar{e}, \bar{e}^{\prime} \}$
satisfying $\bar{b}^{\prime} + \bar{q}= \bar{e}+\bar{e}^{\prime}$ and
$b_1^{\prime}+b_3^{\prime} =\frac{2d_4-d_1-d_3}{d_1+d_3+d_4}.$
To avoid a contradiction, we must have $d_1+d_3 \geq 2d_4$.

We can repeat the above argument with the null vertex pairs
$\{\bar{e}, \bar{a}^{\prime} \}$ and $\{\bar{e}^{\prime}, \bar{a} \}$,
obtaining the inequalities $d_3+d_4 \geq 2d_1$ and $d_1+d_4 \geq 2d_3$
respectively. The three inequalities then imply that in fact
$d_1=d_3=d_4$ and $c=(\frac{2}{3}, -1, \frac{2}{3}, -\frac{4}{3})$.
Furthermore, ${\mathcal C}\cap E=\{\bar{a}, \bar{a}^{\prime}, \bar{c},
\bar{e}, \bar{e}^{\prime}, \bar{t}, \bar{q} \}$ and the null condition
for $\bar{c}$ gives $(d_1, d_2)=(3, 9)$ or $(4, 3)$.

By looking at the terms in the superpotential equation corresponding
to the vertices (all of type II), we find that the coefficients
$F_{\bar{c}}, F_{\bar{e}}, F_{\bar{e}^{\prime}}$ have the same sign,
which is opposite to that of $F_{\bar{a}}, F_{\bar{a}^{\prime}}, F_{\bar{q}}$.
Next we note that the only ways to write $d+(\frac{1}{3}, -1, \frac{1}{3},
-\frac{2}{3})$ (resp. $d+(-\frac{1}{3}, -1, \frac{2}{3}, -\frac{1}{3})$)
as a sum of element of $\mathcal C$ are $\bar{t}+\bar{c}=\bar{a}+\bar{a}^{\prime}$
(resp. $\bar{t}+\bar{a}=\bar{c}+\bar{e}$). The superpotential equation then
gives
$ F_{\bar{a}} F_{\bar{a}^{\prime}} J(\bar{a}, \bar{a}^{\prime}) + F_{\bar{t}} F_{\bar{c}} J(\bar{t}, \bar{c}) = 0$
and $ F_{\bar{c}} F_{\bar{e}} J(\bar{c}, \bar{e}) + F_{\bar{t}} F_{\bar{a}} J(\bar{t}, \bar{a}) = 0.$
Since $J(\bar{a}, \bar{a}^{\prime}), J(\bar{c}, \bar{e}), J(\bar{t}, \bar{c})$ and
$J(\bar{t}, \bar{a})$ are all positive, the above equations and facts
imply that $F_{\bar{c}}$ and $F_{\bar{a}}$ have the same sign, a contradiction.

So $c$ cannot lie in region  II.

\smallskip

{\noindent $c$} {\em lying in region V}:

We have $c_3<-1<-c_4 <1< c_1$. The adjacent (1B) vertices assumption
implies that $\bar{w}=\frac{1}{2}(\bar{c}+\bar{a})$ and
$\bar{y}=\frac{1}{2}(\bar{c}+\bar{a}^{\prime})$ for some $\bar{a}, \bar{a}^{\prime}
\in {\mathcal C} \cap E$. It follows that $a=(2+c_3 +c_4, -1, -c_3, -2-c_4)$
and $a^{\prime}=(c_3 +c_4, -1, -2-c_3, 2-c_4)$. The null conditions on these
vectors give
$$ c_1=\frac{(2d_1+d_4)(d_3+d_4)}{d_4(d_1 +d_3+ d_4)}, \ \
    c_3=-\left(1+\frac{d_1}{d_4}\right)\left(\frac{2d_3+d_4}{d_1+d_3+ d_4}\right), \ \
    c_4=\frac{d_3-d_1}{d_1+d_3+d_4}. $$

Since $a_3 = -c_3 > 1$, $a$ lies above the line $uv$. Also, $a_1^{\prime} =
c_3+c_4 =-c_1 < -1$, so $a^{\prime}$ lies below the line $uz$. We can
therefore apply Theorem \ref{Vertex} to $\bar{a}$ and $\bar{a}^{\prime}$
to get the following possibilities:

(i) $\bar{u} \in {\mathcal C}, \; J(\bar{a}, \bar{u})=0=J(\bar{a}^{\prime}, \bar{u})$,

(ii) $\bar{u} \notin {\mathcal C}$, \;  $e:=2u-a, \; e^{\prime}:=2u-a^{\prime}$
   lie in ${\mathcal C } \cap E$ and are null.

If (i) occurs, then the two orthogonality conditions imply that
$d_1=d_3$, so $c_4 =0, \; c_1 =-c_3 = 1+ \frac{d_1}{d_4}$. Substituting these
values of $c_i$ into $J(\bar{a}^{\prime}, \bar{u})=0$ gives
$1 = \frac{2}{d_4} + \frac{1}{d_2}$. But the null condition for
$\bar{c}$ is $1 = \frac{1}{d_2} + \frac{2}{d_1}(1+\frac{d_1}{d_4})^2 >
\frac{1}{d_2} + \frac{2}{d_4} = 1$, which is a contradiction.

Hence (ii) must occur. Note that if the above diagram is rotated so that
the lines $x_1+x_3 = \kappa$ (for arbitrary constants $\kappa$) are horizontal,
then the lines $x_1 -x_3 = \kappa$ would be vertical. $u$ is the only
point in the hexagon lying on $x_1 -x_3 =-2$. Observe that
$a_1 -a_3 =a_1^{\prime}-a_3^{\prime} \geq -2$, otherwise
$\bar{a}\bar{a}^{\prime}$ would not intersect ${\rm conv}(\frac{1}{2}(d+{\mathcal W}))$,
which contradicts Cor \ref{meet}. If, however, $a_1 -a_3 > -2$,
then $\bar{e}\bar{e}^{\prime}$ would not intersect
${\rm conv}(\frac{1}{2}(d+{\mathcal W}))$. So in fact
$a=e^{\prime}, \; e=a^{\prime}$ and $u$ all lie on $x_1 -x_3 = -2$.
In other words, the hexagon is circumscribed by the triangle $caa^{\prime}$
with intersections at $w, u$ and $y$.

It follows easily from the above that $c=(2, -1, -2, 0)$, $d_1=d_3=d_4$,
and the null condition for $\bar{c}$ is $1=\frac{8}{d_1} + \frac{1}{d_2}$.
Also, we have ${\mathcal C} \cap E = \{\bar{c}, \bar{a}, \bar{a}^{\prime}, \bar{t} \}$.
Since $w, u, y$ are type II, by Lemma \ref{Jpos}, we see that
the signs of $F_{\bar{a}}$, $F_{\bar{c}}$, and $F_{\bar{a}^{\prime}}$ in the
superpotential equation cannot be chosen compatibly. We have thus shown that
the hexagon (H3) cannot occur.
\end{ex}

The hexagon (H2) is not regular, but has reflection symmetry
about $uv$ and the perpendicular bisector of $yy^{\prime}$.
It can be eliminated by similar arguments, but we now have
to consider $c$ lying in regions  III and IV as well.
The hexagon (H1) can also be eliminated by the above methods.
Here the hexagon is invariant under the symmetric group
permuting the coordinates $x_1, x_2, x_3$. Together with
Cor \ref{meet}, this fact reduces our consideration to those
$c$ lying in three of the regions formed by extending the
sides of the hexagon.

As mentioned in Remark \ref{subshape}, we also need to rule out subshapes
of the hexagons. For (H2) and (H3) the methods used above can also be applied
to rule out all the sub-parallelograms and trapezia except the rectangle
$y y^{\prime} z^{\prime} z$ of (H2) (see Lemma \ref{subH2} and the
discussion immediately before Ex \ref{example-parallel}).
All sub-triangles will be dealt with at the end of this section.
(There is a triangle with midpoint in (H2) but that can be dealt with
by similar methods.)  For (H2) this leaves the pentagon
$yy^{\prime}vz^{\prime}z$ and the kite $y^{\prime}uz^{\prime}v$,
both of which can still be eliminated using the above methods.

The possible subshapes of (H1) are rather numerous. However, if
$r\geq 4$ we will be able to eliminate all of them in Lemma \ref{hexH1}.
Without this assumption, the above methods can be used to eliminate
those subshapes which do not contain all three type I vectors. Of course
the following discussion will handle the sub-triangles.

\smallskip

Lastly, we consider triangular faces.

\begin{theorem} \label{triangle}
Suppose we have adjacent $($1B$)$ vertices in $\Delta^{\bar{c}}$ corresponding
to a triangular face $\bar{x} \bar{x}^{\prime} \bar{x}^{\prime \prime}$ of
${\rm conv}(\frac{1}{2}(d+ {\mathcal W}))$. Let $E$ be the affine $2$-plane
determined by the triangular face. So there are null vectors
$\bar{a}, \bar{a}^{\prime}$ in $\mathcal C \cap E$ such that
$x = \frac{1}{2}(a+c), \; x^{\prime}= \frac{1}{2}(a^{\prime} + c).$

Suppose the vertices of the triangle are the only elements of $\mathcal W$
in the face. Then we are in one of the following two situations:

$($i$)$ ${\mathcal C} \cap E = \{ \bar{c}, \bar{a}, \bar{a}^{\prime}, \bar{x}^{\prime \prime} \}$,
with $c + x^{\prime \prime} = a + a^{\prime}$
and $J(\bar{x}^{\prime \prime}, \bar{a}) = J(\bar{x}^{\prime \prime}, \bar{a}^{\prime}) =0;$

$($ii$)$ ${\mathcal C} \cap E = \{ \bar{c}, \bar{a}, \bar{a}^{\prime} \}$
where $\frac{1}{2}(a + a^{\prime}) = x^{\prime \prime}$, one  of $x, x^{\prime}, x^{\prime \prime}$
is type I, and the others are either both type I or both type II/III.
\end{theorem}

\begin{proof}

\begin{picture}(400, 130)(-40, 10)
\put(10, 75){\line(2,1){130}}
\put(10, 75){\line(3,-1){130}}
\put(10, 75){\circle{4}}
\put(70, 105){\circle*{4}}
\put(70, 55){\circle*{4}}
\put(130, 135){\circle{4}}
\put(130, 35){\circle{4}}
\put(141.43, 90.72){\circle*{4}}
\put(3, 80){\it c}
\put(70, 110){\it x}
\put(130, 140){\it a}
\put(65, 45){$x^{\prime}$}
\put(126, 22){$a^{\prime}$}
\put(143, 95){$x^{\prime \prime}$}
\put(-30, 75){({\it i})}
\put(230, 75){\circle{4}}
\put(230, 75){\line(2, 1){100}}
\put(230, 75){\line(2, -1){100}}
\put(330, 25){\line(0, 1){100}}
\put(280, 100){\circle*{4}}
\put(280, 50){\circle*{4}}
\put(330, 125){\circle{4}}
\put(330, 25){\circle{4}}
\put(330, 75){\circle*{4}}
\put(222, 75){\it c}
\put(280, 105){\it x}
\put(330, 130){\it a}
\put(275, 40){$x^{\prime}$}
\put(340, 75){$x^{\prime \prime}$}
\put(320, 15){$a^{\prime}$}
\put(190, 75){({\it ii})}
\thicklines
\put(70, 55){\line(0,1){50}}
\put(70, 105){\line(5, -1){71.43}}
\put(280, 50){\line(2, 1){50}}
\put(280, 50){\line(0,1){50}}
\put(70, 55){\line(2,1){71.43}}
\put(280, 100){\line(2, -1){50}}
\end{picture}

\noindent{(A)} We may introduce coordinates in $E$ so that $\bar{x} \bar{x}^{\prime}$
is vertical and to the right of $\bar{c}$. As $\bar{a} \bar{a}^{\prime}$
must meet ${\rm conv}(\frac{1}{2}(d+{\mathcal W}))$, we see
$\bar{x}^{\prime \prime}$ is on or to the right of $\bar{a} \bar{a}^{\prime}$.
Let $\bar{b}$ be any leftmost point of $(\mathcal C \cap E) \setminus \{ \bar{c} \}$.
As in Theorem \ref{parallel}, we see that $\bar{b} + \bar{c} \in d+ \mathcal W$,
 so all elements of $\mathcal C \cap E$ except $\bar{c},\bar{a},\bar{a}^{\prime}$
 are to the right of $\bar{a} \bar{a}^{\prime}$.

\noindent{(B)} Considering $\bar{a} \bar{x}^{\prime \prime}$ and
 $\bar{a}^{\prime} \bar{x}^{\prime \prime}$
we see (using Theorem \ref{Vertex} and Cor \ref{meet}) that either
\begin{enumerate}
\item[(1)] $\bar{x}^{\prime \prime} \in \mathcal C$ and
$J(\bar{x}^{\prime \prime}, \bar{a})=0=J(\bar{x}^{\prime \prime}, \bar{a}^{\prime})$,  or
\item[(2)] $\bar{x}^{\prime \prime} \notin \mathcal C$ and  $x^{\prime \prime} =
 \frac{1}{2}(a + a^{\prime})$.
\end{enumerate}
In case (1), $(\bar{x}^{\prime \prime})^\perp \cap E$ is the line
through $\bar{a} \bar{a}^{\prime}$.
By Prop. \ref{sep-hypl} and Cor. \ref{meet},
observe that all elements of $({\mathcal C}\cap E) \setminus
\{\bar{x}^{\prime \prime} \}$ are left of $\bar{x}^{\prime \prime}$.
 Let $\bar{b}$ be a rightmost
element of $({\mathcal C}\cap E) \setminus \{\bar{x}^{\prime \prime} \}$.
So either $J(\bar{b}, \bar{x}^{\prime \prime})=0$ or
$\bar{b} + \bar{x}^{\prime \prime} \in d + \mathcal W$. Since
$\bar{b}$ is not to the left of $\bar{a} \bar{a}^{\prime}$,
 the second alternative cannot hold and so $\bar{b}$ must lie on
  $\bar{a} \bar{a}^{\prime}$.
Combining this with our results in (A), we see $\mathcal C \cap E$
is as in (i). Also, as $J(\bar{a}, \bar{a}^{\prime})>0$
 and $\bar{a} + \bar{a}^{\prime} \notin  d + \mathcal W$, we see
 $a + a^{\prime}$ must equal $c + x^{\prime \prime}$.

In case (2), by Cor \ref{meet} there are no elements of $\mathcal C \cap E$
right of $\bar{a} \bar{a}^{\prime}$. Hence $\mathcal C \cap E$ is as in (ii).
Now $J(\bar{b}, \bar{e}) >0$ for all $\bar{b} \neq \bar{e}$ in
${\mathcal C} \cap E$, so the last statement of (ii) follows.
\end{proof}

\begin{rmk} \label{tri-mdpt}
We must also consider the case when some midpoints of the sides of
our triangular face lie in $\mathcal W$. (This could happen if two
vertices were  $(1,-1,-1,\cdots), (-1,1,-1, \cdots)$ or
$(1,-2,\cdots), (1,0,-2,\cdots)$ or $(1,-2,\cdots), (-1,0,\cdots)$.) Let us
denote the midpoints of $x x^{\prime}, x x^{\prime \prime}$ and $x^{\prime}
x^{\prime \prime}$ respectively by $z,y,t$.

If $z$ is absent, the arguments of (A) in the proof of Theorem \ref{triangle}
still hold, so we have the alternatives (1),(2) in (B).
If (1) holds then, choosing $\bar{b}$ as above, if $b$ is right of
$aa^{\prime}$, we have $\frac{1}{2}(b + x^{\prime \prime}) \in \mathcal W$.
This gives a contradiction since $\frac{1}{2}(b + x^{\prime \prime})$
cannot be  $y$ or $t$ as $b \neq x, x^{\prime}$. Now $\mathcal C \cap E=
\{\bar{c},\bar{a},\bar{a}^{\prime}, \bar{x}^{\prime \prime} \}$,
and as $c + x^{\prime \prime} \notin 2 \mathcal W$ it must equal
$a + a^{\prime}$. It follows that the midpoints $y,t$ cannot arise.
If instead (2) holds, then $\mathcal C \cap E = \{\bar{c},
\bar{a}, \bar{a}^{\prime} \}$ and again no midpoints can be present.

Suppose now the midpoint $z$ of $x x^{\prime}$ is present. The
argument of (A) shows that to account for $z$, $\frac{1}{2}(\bar{a} + \bar{a}^{\prime})
\in {\mathcal C}$, and all elements of $(\mathcal C \cap E) \setminus \{ \bar{c},
\bar{a}, \bar{a}^{\prime}, (\bar{a} + \bar{a}^{\prime})/2 \}$ are right
of $\bar{a} \bar{a}^{\prime}$.  We still have the alternatives (1) and (2),
but (2) immediately gives a contradiction.

In (1) we see as before there are no elements of $\mathcal C \cap E$
lying to the right of $\bar{a} \bar{a}^{\prime}$, so
$\mathcal C \cap E= \{\bar{c},\bar{x}^{\prime \prime}, \bar{a}, \bar{a}^{\prime},
(\bar{a} + \bar{a}^{\prime})/2 \}$.  Note
that $J(\bar{a},(\bar{a} + \bar{a}^{\prime})/2)$ and
$J(\bar{a}^{\prime},(\bar{a} + \bar{a}^{\prime})/2) >0$.

If $c + x^{\prime \prime}= a+a^{\prime}$, we find
after some algebra that $\bar{a} +
(\bar{a} + \bar{a}^{\prime})/2 \neq 2 \bar{y}$ and also cannot be
written as a different sum of elements of ${\mathcal C}$, giving a
contradiction.

If $\bar{c} + \bar{x}^{\prime \prime} \neq \bar{a} + \bar{a}^{\prime}$
then one sees that $\bar{c} + \bar{x}^{\prime \prime} \notin d+{\mathcal W}$,
and by relabelling $x$ and $x^{\prime}$, $a$ and $a^{\prime}$
we may assume that  $\bar{c} + \bar{x}^{\prime \prime}= \bar{a} +
\frac{1}{2}(\bar{a} + \bar{a}^{\prime})$ and also
$\bar{a} + \bar{a}^{\prime} = 2 \bar{y}$ and
$\bar{a}^{\prime} + \frac{1}{2}(\bar{a} +
\bar{a}^{\prime}) = 2 \bar{t} = \bar{x}^ \prime + \bar{x}^{\prime
  \prime}$. These relations imply $a = x^{\prime}$, a contradiction.
So no triangle with any midpoints present can arise.
\end{rmk}

\begin{rmk} \label{tri-2mdpt}
There are also triangular faces with two points of $\mathcal W$ in
the interior of an edge.  This can only happen if two vertices are
$(-2,1, 0, \cdots)$ and $(1,-2, 0, \cdots)$ (up to permutation). The
other sides of the triangle now have no interior points in
$\mathcal W$ unless the triangle is contained in the hexagon (H1).
We can again modify the proof of Theorem \ref{triangle} to treat
this situation.

If the interior points $z,w$ lie on $x x^{\prime}$, then
$(2 \bar{a} + \bar{a}^{\prime})/3, (\bar{a} + 2 \bar{a}^{\prime})/3$
must be in $\mathcal C$, and all points of $\mathcal C \cap E$ except
for these two and $\bar{c},\bar{a}, \bar{a}^{\prime}$ lie to the right of
$\bar{a} \bar{a}^{\prime}$. By Prop \ref{sep-hypl}, alternative (1) must
now hold. The usual argument shows $\bar{x}^{\prime \prime}$ is the
only element of $\mathcal C \cap E$ on the right of $\bar{a} \bar{a}^{\prime}$.
Now again $J(\bar{a}, \frac{1}{3}(2\bar{a}^{\prime} + \bar{a}))> 0,$
$J(\bar{a}^{\prime}, \frac{1}{3}(\bar{a}+ 2\bar{a}^{\prime})) > 0$,
and the sums $a + (2a^{\prime} +a)/3$ and $a^{\prime} + (a + 2 a^{\prime})/3$
cannot give points in $2 \mathcal W$. Since they also cannot both be
cancelled by $c + x^{\prime \prime}$ in the superpotential equation,
we have a contradiction.

The other possibility for two interior points is, after relabelling
the vertices if necessary, when $z = (2x + x^{\prime \prime})/3$ and
$w = (2x^{\prime \prime} + x)/3$. As usual all elements of
$\mathcal C \cap E$ except for $\bar{c}, \bar{a},\bar{a}^{\prime}$ are
on the right of $\bar{a} \bar{a}^{\prime}$. Alternative (1) must hold,
or else we cannot account for $z,w$. The usual argument shows either
$\bar{x}^{\prime \prime}$ is the only element of $\mathcal C \cap E$ right of
$\bar{a} \bar{a}^{\prime}$, {\em or} $z \in \mathcal C$ is the rightmost
element of $({\mathcal C} \cap E) \setminus \{ \bar{x}^{\prime \prime} \}$
(so $(z + x^{\prime \prime})/2 =w)$. In the former case we cannot get
both $z$ and $w$, as $(c + x^{\prime \prime})/2$ can't equal both $z$ and $w$.
In the latter, considering $\bar{a} \bar{z}$ shows $J(\bar{a}, \bar{z})=0$.
But as $J(\bar{a}, \bar{x}^{\prime \prime})=0$, this means $\bar{a}$
is orthogonal to $\bar{x}$ and hence to $\bar{c}$, a contradiction.

So no triangle with points of $\mathcal W$ in the interior of an edge
can arise (except possibly for a subtriangle of (H1)).
\end{rmk}

Nullity of $\bar{c}, \bar{a}, \bar{a}^{\prime}$ and the conditions in
Theorem \ref{triangle}(i),(ii) again put severe constraints on
$x, x^{\prime}, x^{\prime}$ and the dimensions. The possible triangles for
case (i) are as follows, where (Tr11)-(Tr22) occur only if $K$ is not connected,
and we have also listed the vectors $c, a, a^{\prime}$ for future reference.
Further details of how the following listing is arrived at can be found
in \cite{DW5}.


\[
\begin{array}{|c|c|c|c|} \hline
       & x^{\prime \prime} &        x          &    x^{\prime}       \\ \hline
(Tr1)  & (-2,1,0,0,0)      &  (0,0,-2,1,0)     &  (0,0,-2,0,1)  \\
(Tr2)  & (-2,1,0,0)        &  (0,1,-2,0)       &  (0,1,-1,-1)   \\
(Tr3)  & (0,0,0,-2,1)      &  (-2,1,0,0,0)     & (0,1,-2,0,0)   \\
(Tr4)  &(-2,1,0,0,0,0)     &  (0,0,-2,1,0,0)   & (0,0,0,1,-1,-1)\\
(Tr5)  &(-2,1,0,0,0)       & (0,1,-1,0,-1)     & (0,1,-1,-1,0)  \\
(Tr6)  &(-2,1,0,0,0,0)     &(0,0,1,-1,-1,0)    & (0,0,1,-1,0,-1)\\
(Tr7)  &(-2,1,0,0,0,0)     &(0,0,1,-1,-1,0)    & (0,0,-1,-1,0,1)\\
(Tr8)  &(-2,1,0,0,0,0)     &(0,0,-1,-1,1,0)    &(0,0,-1,-1,0,1) \\
(Tr9)  &(-2,1,0,0,0,0,0)   &(0,0,1,-1,-1,0,0)  &(0,0,1,0,0,-1,-1)\\
(Tr10) &(-2,1,0,0,0,0,0)   &(0,0,-1,1,-1,0,0)  &(0,0,-1,0,0,1,-1)\\ \hline
(Tr11) &(0,0,0,1,-1,-1)    &(-2,1,0,0,0,0)     &(-2,0,1,0,0,0)  \\
(Tr12) &(0,1,0,-1,-1)      &(-2,1,0,0,0)       &(-1,1,-1,0,0)   \\
(Tr13) &(0,0,0,-1,-1,1)    &(-2,1,0,0,0,0)     &(0,1,-2,0,0,0)  \\
(Tr14) &(0,0,0,-1,-1,1)    &(0,1,-1,-1,0,0,0)  &(-2,1,0,0,0,0,0) \\
(Tr15) &(0,0,0,1,-1,-1,0)  &(1,-1,-1,0,0,0,0)  &(1,-1,0,0,0,0,-1)  \\
(Tr16) &(0,0,0,1,-1,-1,0)  &(1,-1,-1,0,0,0,0)  &(0,-1,1,0,0,0,1)  \\
(Tr17) &(0,0,0,1,-1,-1,0)  &(1,-1,-1,0,0,0,0)  &(0,-1,-1,0,0,0,1)  \\
(Tr18) &(0,0,0,1,-1,-1,0)  &(1,-1,-1,0,0,0,0,0)&(1,0,0,0,0,0,-1,-1)  \\
(Tr19) &(0,0,0,1,-1,-1,0)  &(1,-1,-1,0,0,0,0,0)&(0,-1,0,0,0,0,-1,1)  \\
(Tr20) &(-1,-1,1,0,0,0)    &(0,0,1,-1,-1,0)    &(0,0,1,-1,0,-1)    \\
(Tr21) &(-1,1,-1,0,0,0)    &(0,0,-1,1,-1,0)    &(0,0,-1,1,0,-1)    \\
(Tr22) &(-1,1,-1,0,0,0)    &(0,0,-1,-1,1,0)    &(0,0,-1,-1,0,1)  \\ \hline
\end{array}
\]

\medskip

\[
\begin{array}{|c|c|c|c|} \hline
      &  3c                &  3a & 3a^{\prime} \\  \hline
(Tr1) & (2,-1,-8,2,2)      & (-2,1,-4,4,-2)        &  (-2,1,-4,-2,4) \\
(Tr2) & (2,3,-6,-2)        & (-2,3,-6,2)           &  (-2,3,0,-4) \\
(Tr3) &(-4,4,-4,2,-1)      &(-8,2,4,-2,1)          &  (4,2,-8,-2,1) \\
(Tr4) &(2,-1,-4,4,-2,-2)   &(-2,1,-8,2,2,2)        &  (-2,1,4,2,-4,-4)\\
(Tr5) &(2,3,-4,-2,-2)      &(-2,3,-2,2,-4)         &(-2,3,-2,-4,2)\\
(Tr6) &(2,-1,4,-4,-2,-2)   &(-2,1,2,-2,-4,2)       &(-2,1,2,-2,2,-4)\\
(Tr7) &(2,-1,0,-4,-2,2)    &(-2,1,6,-2,-4,-2)      &(-2,1,-6,-2,2,4)\\
(Tr8) &(2,-1,-4,-4,2,2)    &(-2,1,-2,-2,4,-2)      &(-2,1,-2,-2,-2, 4)\\
(Tr9) &(2,-1,4,-2,-2,-2,-2)&(-2,1,2,-4,-4,2,2)     &(-2,1,2,2,2,-4,-4)\\
(Tr10)&(2,-1,-4,2,-2,2,-2) &(-2,1,-2,4,-4,-2,2)    &(-2,1,-2,-2,2,4,-4) \\ \hline
(Tr11)&(-8,2,2,-1,1,1)     &(-4,4,-2,1,-1,-1)      &(-4,-2,4,1,-1,-1)   \\
(Tr12)&(-6,3,-2,1,1)       &(-6,3,2,-1,-1)         &(0,3,-4,-1,-1)     \\
(Tr13)&(-4,4,-4,1,1,-1)    &(-8,2,4,-1,-1,1)       &(4,2,-8,-1,-1,1)   \\
(Tr14)&(-4,4,-2,-2,1,1,-1) &(4,2,-4,-4,-1,-1,1)    &(-8,2,2,2,-1,-1,1) \\
(Tr15)&(4,-4,-2,-1,1,1,-2) &(2,-2,-4,1,-1,-1,2)    &(2,-2,2,1,-1,-1,-4)  \\
(Tr16)&(2,-4,0,-1,1,1,-2)  &(4,-2,-6,1,-1,-1,2)    &(-2,-2,6,1,-1,-1,-4)  \\
(Tr17)&(2,-4,-4,-1,1,1,2)  &(4,-2,-2,1,-1,-1,-2)   &(-2,-2,-2,1,-1,-1,4)  \\
(Tr18)&(4,-2,-2,-1,1,1,-2,-2)&(2,-4,-4,1,-1,-1,2,2)&(2,2,2,1,-1,-1,-4,-4)  \\
(Tr19)&(2,-4,-2,-1,1,1,-2,2)&(4,-2,-4,1,-1,-1,2,-2)&(-2,-2,2,1,-1,-1,-4,4)  \\
(Tr20)&(1,1,3,-4,-2,-2)    &(-1,-1,3,-2,-4,2)      &(-1,-1,3,-2,2,-4)   \\
(Tr21)&(1,-1,-3,4,-2,-2)   &(-1,1,-3,2,-4,2)       &(-1,1,-3,2,2,-4)   \\
(Tr22)&(1,-1,-3,-4,2,2)    &(-1,1,-3,-2,4,-2)      &(-1,1,-3,-2,-2,4)  \\ \hline
\end{array}
\]


\begin{rmk}
In making the above table, it is useful to observe from the nullity
and orthogonality conditions that $x^{\prime \prime}$ cannot be type
I, and that if $x^{\prime \prime}$ is type III, say, $(-2^i, 1^j),$
then $x_i = x^{\prime}_i$ iff $x_j = x^{\prime}_j$.
\end{rmk}

The possibilities for Theorem \ref{triangle}(ii) are as follows
(up to permutation of $x, x^{\prime}, x^{\prime \prime}$ and the
corresponding permutation of $c, a, a^{\prime}$):
\[
\begin{array}{|c|c|c|c|} \hline
&    x^{\prime \prime}       &        x           &    x^{\prime}  \\ \hline
(Tr23)  &   (-1,0,0,0,0)     &  (0,-2,1,0,0)      & (0,0,0,-2,1)  \\
(Tr24)  &   (-1,0,0,0)       &  (0,1,-2,0)        & (0,-1,-1,1)   \\
(Tr25)  &   (-1,0,0,0,0,0)   &  (0,1,-2,0,0,0)    & (0,0,0,1,-1,-1)\\
(Tr26)  &   (-1,0,0,0,0)     &  (0,1,-1,-1,0)     & (0,-1,-1,0,1)  \\
(Tr27)  &   (-1,0,0,0,0,0,0) &  (0,1,-1,-1,0,0,0) & (0,0,0,0,1,-1,-1)\\
(Tr28)  &   (-1,0,0)         &  (0,-1,0)          & (0,0,-1) \\ \hline
\end{array}
\]

\smallskip

\[
\begin{array}{|c|c|c|c|} \hline
       &     c               &      a             &  a^{\prime} \\ \hline
(Tr23) &  (1,-2,1,-2,1)      & (-1,-2,1,2,-1)     & (-1,2,-1,-2,1)\\
(Tr24) &  (1,0,-3,1)         &(-1,2,-1,-1)        &(-1,-2,1,1)\\
(Tr25) &  (1,1,-2,1,-1,-1)   & (-1,1,-2,-1,1,1)   &(-1,-1,2,1,-1,-1)\\
(Tr26) &  (1,0,-2,-1,1)      &(-1,2,0,-1,-1)      &(-1,-2,0,1,1) \\
(Tr27) &  (1,1,-1,-1,1,-1,-1)&(-1,1,-1,-1,-1,1,1) & (-1,-1,1,1,1,-1,-1) \\
(Tr28) &  (1,-1,-1)          & (-1,-1,1)          & (-1,1,-1) \\  \hline
\end{array}
\]

\begin{rmk}
In drawing up the above listing, recall from Theorem \ref{triangle}
that one of the vectors, without loss of generality $x^{\prime \prime}$,
is of type I. We write $x^{\prime \prime} = (-1,0,0,\cdots)$. It now
easily follows from nullity and
the relations between $x, x^{\prime}, x^{\prime \prime}$
and $c, a, a^{\prime}$ that $x_1 = x^{\prime}_1$.

Also, observe that as $x^{\prime \prime}$ is a \ vertex of $\mathcal W$,
no type II vector may have a nonzero entry in the first position.

In contrast to the earlier listing of non-triangular faces, the above
lists result from examining {\em all} triangular faces, including ones
which arise from other faces because certain vertices are absent from $\mathcal W$.
\end{rmk}

\smallskip

The restrictions on the dimensions of the corresponding summands are as follows:

\begin{enumerate}
\item[(Tr1)]  \hspace{1cm} $  (2,1,16,4,4,\cdots) $
\item[(Tr2)]  \hspace{1cm} $  (2,3,12,4,\cdots)  $
\item[(Tr3)]  \hspace{1cm} $  (16,4,16,2,1,\cdots) $
\item[(Tr4)]  \hspace{1cm} $  (2,1,16,4,d_5,d_6,\cdots), \; \  \frac{1}{d_4} + \frac{1}{d_5} = \frac{1}{4}$
\item[(Tr5)]  \hspace{1cm} $  (2,3,6,6,6,\cdots) $
\item[(Tr6)]  \hspace{1cm} $  (2,1,d_3,d_4,4,4,\cdots), \; \  \frac{1}{d_3} + \frac{1}{d_4} = \frac{1}{4} $
\item[(Tr7)]  \hspace{1cm} $  (2,1,12,3,12,12,\cdots) $
\item[(Tr8)]  \hspace{1cm} $  (2,1,d_3,d_4,4,4,\cdots), \; \  \frac{1}{d_3} + \frac{1}{d_4} = \frac{1}{4} $
\item[(Tr9)]  \hspace{1cm} $  (2,1,4,d_4,d_5,d_6,d_7,\cdots), \  \ \frac{1}{d_4} + \frac{1}{d_5} =
\frac{1}{d_6} + \frac{1}{d_7} = \frac{1}{4} $
\item[(Tr10)] \hspace{1cm} $   (2,1,4,d_4,d_5,d_6,d_7,\cdots), \ \ \frac{1}{d_4} + \frac{1}{d_5} =
\frac{1}{d_6} + \frac{1}{d_7} = \frac{1}{4} $
\item[(Tr11)] \hspace{1cm} $  (16, 4, 4, 1, 1, 1, \cdots)  $
\item[(Tr12)] \hspace{1cm} $  (12, 3, 4, 1, 1, \cdots)  $
\item[(Tr13)] \hspace{1cm} $  (16, 4, 16, 1, 1, 1, \cdots) $
\item[(Tr14)] \hspace{1cm} $  (16,4,d_3,d_4,1,1,1,\cdots) \ \ \frac{1}{d_3} + \frac{1}{d_4} = \frac{1}{4} $
\item[(Tr15)] \hspace{1cm} $  (d_1, d_2, 4, 1, 1,1,4, \cdots), \ \ \frac{1}{d_1} + \frac{1}{d_2} = \frac{1}{4}  $
\item[(Tr16)] \hspace{1cm} $ (12, 3,12,1,1,1,12, \cdots)  $
\item[(Tr17)] \hspace{1cm} $  (4,d_2,d_3,1,1,1,4, \cdots), \ \ \frac{1}{d_2} + \frac{1}{d_3} = \frac{1}{4}  $
\item[(Tr18)] \hspace{1cm} $  (4, d_2, d_3, 1, 1, 1, d_7, d_8, \cdots), \ \ \frac{1}{d_2} + \frac{1}{d_3} = \frac{1}{4}
              = \frac{1}{d_7} + \frac{1}{d_8}  $
\item[(Tr19)] \hspace{1cm} $  (d_1, 4, d_3, 1, 1, 1, d_7, d_8, \cdots), \ \ \frac{1}{d_1} + \frac{1}{d_3} = \frac{1}{4}
              = \frac{1}{d_7} + \frac{1}{d_8}  $
\item[(Tr20-22)]   \hspace{1cm}$  (1,1,3,6,6,6, \cdots), (1,2,2,8,8,8, \cdots), \ {\rm or} \  (2,1,2,8,8,8, \cdots)  $
\item[(Tr23)] \hspace{1cm} $   \frac{1}{d_1} + \frac{4}{d_2} + \frac{1}{d_3} + \frac{4}{d_4} + \frac{1}{d_5}=1 $
\item[(Tr24)] \hspace{1cm} $  d_3 = 2d_2 \; : \;\; \frac{1}{d_1} + \frac{9}{d_3} + \frac{1}{d_4} =1 $
\item[(Tr25)] \hspace{1cm} $    \frac{1}{d_1} + \frac{1}{d_2} + \frac{4}{d_3} + \frac{1}{d_4} +
\frac{1}{d_5} + \frac{1}{d_6} =1 $
\item[(Tr26)] \hspace{1cm}$  d_2 = d_3 \; : \;\; \frac{1}{d_1}+ \frac{4}{d_3}  +
 \frac{1}{d_4} + \frac{1}{d_5} = 1. $
\item[(Tr27)] \hspace{1cm}$  \sum_{i=1}^{7} \frac{1}{d_i} = 1 $
\item[(Tr28)] \hspace{1cm} $  \frac{1}{d_1} + \frac{1}{d_2} + \frac{1}{d_3} =1. $
\end{enumerate}

Note that (Tr28) is a subtriangle of (H1), (Tr2) is a subtriangle of
a triangle with midpoints of all sides in $\mathcal W$, and (Tr12)
is a subtriangle of a triangle with the midpoint of one side.

Let us now illustrate by an example how one arrives at the above tables.

\begin{ex} \label{ex-tri}
One possible triangle has vertices $V_1=(0,0,0,-2,1)$, $V_2=(-2,1,0,0,0)$,
$V_3=(0, 1, -2,0,0)$ with the midpoint $V_4=(-1, 1, -1, 0, 0)$ of $V_2V_3$
in $\mathcal W$. The triangle has a symmetry given by interchanging the first
and third entries. It therefore suffices to consider $V_1, V_2, V_4$ as possibilities
for $x^{\prime \prime}$. Of course, by Remark \ref{tri-mdpt} the full triangle
cannot occur. The possible subtriangles $x x^{\prime} x^{\prime \prime}$
are $V_2 V_3 V_1$,  $V_2 V_4 V_1, V_4 V_1 V_2, V_3 V_1 V_2,$ and $V_3 V_1 V_4$.
Now $3c= 2x + 2x^{\prime} - x^{\prime \prime}, 3a=4x-2x^{\prime} + x^{\prime \prime}$,
and $3a^{\prime} = -2x + 4x^{\prime} + x^{\prime \prime}$ can be used to compute these
vectors in each case. For $V_2 V_4 V_1$ one gets $3c=(-6, 4, -2, 2, -1), 3a=(-6, 2, 2, -2, 1)$
and so $\bar{c}$ and $\bar{a}$ cannot be both null. Similarly, for the last three possibilities,
$\bar{a}^{\prime}$ and $\bar{c}$ cannot be both null. That leaves the first case,
 which gives (Tr3). The condition $J(\bar{a}, \bar{x}^{\prime \prime}) =0$
 is $3 =\frac{4}{d_4} + \frac{1}{d_5}$, which implies $(d_4, d_5)=(2, 1)$.
Putting this into the null conditions for $\bar{c}, \bar{a}, \bar{a}^{\prime}$
gives the equations
$$ 6=\frac{16}{d_1} + \frac{16}{d_2} + \frac{16}{d_3}
  =\frac{64}{d_1} + \frac{4}{d_2} + \frac{16}{d_3}
  =\frac{16}{d_1} + \frac{4}{d_2} + \frac{64}{d_3}.$$
The last two equations imply that $d_1=d_3$ and the first two equations give
$d_1 =4d_2$. These in turn give $(d_1, d_2, d_3)=(16, 4, 16)$, as in the tables above.
\end{ex}

Putting all the results in this section together we obtain

\begin{theorem} \label{1B-initial}
If we have two adjacent $($1B$)$ vertices, then the associated $2$-face
of $\mathcal W$ is given by a triangle in the list $(Tr1)-(Tr27)$,
the square with midpoint $($S$)$,  a proper subshape of the hexagonal
face $($H1$)$ containing all three type I vectors,
or the sub-rectangle $y y^{\prime} z z^{\prime}$ of $($H2$)$. \; $\qed$
\end{theorem}

We note for future reference the following properties of the $c$ vector
of the non-triangular faces appearing in the above theorem: for (S), all
nonzero entries have the same absolute value, and there are only $3$
(resp. $2$) nonzero entries for the subfaces of (H1) (resp. (H2)).

\section{\bf More than one type (2) vertex}

In this section we shall now show there is at most one type (2)
vertex in $\Delta^{\bar{c}}$, except in the situation of Theorem \ref{Fano}
and one other possible case.

Suppose we have two type (2) vertices of $V$. Then we have elements
$v,w, v^{\prime}, w^{\prime}$ of $\mathcal W$ with $c,v,w$ collinear
and $c,v^{\prime}, w^{\prime}$ collinear. So we have four coplanar
elements $v,w, v^{\prime}, w^{\prime}$ of $\mathcal W$ where $vw$ and
$v^{\prime} w^{\prime}$ are edges.  Moreover, the edges $vw$ and
$v^{\prime} w^{\prime}$ meet at $c$ outside conv$(\mathcal W)$. Hence $v
w v^{\prime} w^{\prime}$ do not form a parallelogram or a triangle.

\medskip
{}From our listing of polygons in \S 6 and considering their sub-polygons
we see that the possibilities for further analysis are the following:
\begin{itemize}
\item Trapezia (T1)-(T6): We must have $c=2v-s = 2u-w$. Also, we note
for future reference that $sw$ is always an edge of ${\rm conv}(\mathcal W)$
in (T3) and (T5), regardless of whether or not the whole trapezium is a face,
since $sw$ can be cut out by $\{x_2 =1, \ x_1 + x_3 = -2\}$ (cf \ref{facts}(e)).

\item Hexagons (H1)-(H3)

\item Rectangle with midpoints (P17): While the rectangle itself cannot
occur, we need to consider the trapezia obtained by omitting one vertex,
so that the edges are a side of the rectangle and the segment joining
the remaining vertex to the opposite midpoint. As above, note that
the longer of the two parallel sides of the trapezium is always
an edge of ${\rm conv}(\mathcal W)$.

\item Parallelogram with midpoints (P16): This case is similar to (P17).
The sub-polygons to consider are the trapezia obtained by omitting one
vertex of the parallelogram. By symmetry, we are reduced to omitting either
$u$ or $w$. But since $s$ occurs, $w$ cannot be omitted.

\item Triangle with midpoints of all sides: We need to consider the
trapezia obtained by omitting a vertex. By symmetry all three trapezia
are equivalent. This triangle is always a face of ${\rm conv}(\mathcal W)$
as it is cut out by $\{x_2 =1, \ x_1 + x_3 + x_4 = -2\}$ (cf \ref{facts}(e)).

\item Trapezia (T*1),(T*2), (T*3): By Rmk \ref{nonface} these cannot
be faces of ${\rm conv}(\mathcal W)$, so cannot come from {\em adjacent}
type (2) vertices of $V$. For (T*1), besides the full trapezium, we need to
consider the two trapezia obtained by omitting either $s$ or $w$. By symmetry
these are equivalent.
\end{itemize}

For (T1),(T2),(T4),(T6),(T*2),(T*3) we must have $c=2v-s=2u-w$, and so
Lemma \ref{w} applied to $vs$ gives a contradiction.
The same argument works for (P16), as up to permutations,  $c=2v-s=2y-w$.
For (T1*), since Theorem \ref{cedge} rules out  $c=(3v-s)/2 = (3u-w)/2$
(corresponding to the full trapezium), the only other possible $c$ is
$2v-s=2u-r$, and again Lemma \ref{w} rules this out.

For (T3),(T5) and the trapezium coming from the triangle with midpoints,
we need more information from the superpotential equation.
Since $J(\bar{c},\bar{s}), J(\bar{c},\bar{w})>0$, while $A_v, A_u < 0,$
$F_{\bar{s}}, F_{\bar{w}}$ must have the same sign, which must be opposite
to that of $F_{\bar{c}}$. Since $stw$ is always an edge by earlier remarks,
the nullity of $\bar{c}$ implies $J(\bar{s},\bar{w}) > 0$,  contradicting
Prop \ref{edgeCC}(ii).

Essentially the same argument works for (P17), as up to permutations
$c=2s-v=2z-u=(-2,-1,2,0)$.

For (H3) most quadruples cannot give pairs of edges. For we observe
that $u$ (resp. $v, w$) is present iff $x$ (resp. $y, z$) is. Thus, if
$u$ is missing, so is $x$, and $v, w, y, z$ must all be present
(otherwise we do not have a $2$-dimensional polygon). But we now get a
rectangle, which is not admissible.  Hence all vertices are present
and by symmetry we may assume that one of our edges is $uz$ or $uv$.
{}From this we quickly find that the two possible $c$ (up to
permutations) are $(-1,-1,-1, 2) = 2y-x =2z-u$ and
$(1,-1,1,-2)=2v-u=2w-x$. Both cases are ruled out by Lemma \ref{w}.

For (H2), observe that $y$ (resp. $y^{\prime}$) is present iff $z$
(resp. $z^{\prime}$) is. As these four vectors cannot all be absent
(otherwise we do not have a $2$-dim polygon), by the symmetries of
(H2), we can assume $y^{\prime}$ is present. If $v$ is present,
 then all possibilities are eliminated by Theorem \ref{cedge}.
(Note that although $2 \alpha - y^{\prime}= 2z - z^{\prime}$ it is
impossible for $\alpha y^{\prime}$ and $z z^{\prime}$ to both be edges.)
On the other hand, if $v$ is absent, then $y^{\prime} z^{\prime}$ is an
edge. Since the polygon cannot be a parallelogram or a triangle, it follows
that $u$ is present and the polygon is a pentagon. In this case, the only
possibility compatible with Theorem \ref{cedge} is
$c=(0, -1, 2, -2)=2y-u=\frac{3}{2} y^{\prime} -\frac{1}{2}z^{\prime}.$
(This is not a priori ruled out by Theorem \ref{cedge} as
$y^{\prime} z^{\prime}$ has an interior point $\beta$).

To discuss (H1), we write
\[
u=(-2,1,0), \ p=(0,1,-2), \ v=(1,0,-2), \ w=(0,-2,1), \ s=(1,-2,0), \ q=(-2,0,1)
\]
for the vertices,
\[
x=(-1,1,-1),  \  y=(-1,-1,1), \   z=(1,-1,-1),
\]
for midpoints of the longer sides, and
\[
    \alpha=(-1,0,0), \ \beta=(0,0,-1), \ \gamma=(0,-1,0)
\]
for the interior points, with the understanding that the rest of
the components of the above vectors are zero.

As before we consider pairs of vectors which can form edges of an
admissible polygon. We then compute the possibilities for $c$ and
apply Theorem \ref{cedge}. This will eliminate most possibilities.
(For many quadruples of points we can see, as in (H3), that they
cannot all be vertices.) So up to permutations, the remaining possibilities
are as follows.

If no type II is present:

(1) $c = 2 \alpha-u = 2 \beta - p = (0,-1,0,\cdots)$

(2) $c=2u-v = 2q-s = (-5,2,2,\cdots)$

(3) $c= 2u -p = 2q - \gamma = (-4,1,2,\cdots)$

(4) $c=2q-u = 2\alpha - p = (-2, -1, 2, \cdots)$

If all type II are present:

(5) $c=2u-x=2q-y = (-3,1,1,\cdots)$

(6) $c=2u-y = 2x-z= (-3,3,-1,\cdots)$

(7) $c = (3y-z)/2 = 2q -u = (-2,-1,2, \cdots)$

(8) $c = (3p-u)/2 = (3v-s)/2 = (1,1,-3, \cdots)$

Again, we cannot immediately rule out (7) and (8) using \ref{cedge}
because of the presence of interior points.  However for (8) we easily
see using the arguments of \ref{cedge} that the elements of $\mathcal
C$ on the line through $v,s$ are $\bar{c}, \bar{c}_1 = (\bar{v} + \bar{s})/2 =
\bar{z}$ and $\bar{c}_2 = (3 \bar{s} - \bar{v})/2$. Now as $s, v$ are
type III we need $F_{\bar{c}}$ and $F_{\bar{c}_2}$ to have the same
sign, which is the opposite sign to $F_{\bar{c}_{1}}$. But the
superpotential equation now gives a contradiction to the fact that
$A_z < 0$.

In (2)-(7)  Lemma \ref{w} applied respectively to $uv,up,qu, ux,uy,qu$
gives a contradiction. Note that case (4) only occurs when $K$ is
not connected, as the vectors $\beta, \gamma$ are absent (cf \ref{facts}(b)).
We are left with (1), which is precisely the situation of
Theorem \ref{Fano}.  We have therefore proved

\begin{theorem} \label{two2}
Apart from the situation of Theorem \ref{Fano}, the only other possible
case where we can have more than one type $($2$)$ vertex is, up to permutation
of summands, when two type $($2$)$ vertices are adjacent and the $2$-plane
determined by them and $\bar{c}$ intersects ${\rm conv}(\frac{1}{2}(d+ {\mathcal W}))$
in the pentagon with vertices $uyy^{\prime}zz^{\prime}$ contained in the
hexagon $($H2$)$.     \;  $\qed$
\end{theorem}

We will be able to rule this case out in \S 8.

\section{\bf Adjacent (1B) vertices revisited}

We now return to our classification of when adjacent (1B) vertices
can occur.

The idea is as follows: each of the configurations of \S 6
involves, as well as the null vector $\bar{c}$, two new null vectors
$\bar{a}, \ \bar{a}^{\prime}$. Hence the arguments of earlier sections
also apply to $\bar{a}, \bar{a}^{\prime}$. That is, we may consider
the associated polytopes $\Delta^{\bar{a}}$ and $\Delta^{\bar{a}^{\prime}}$.

The following lemma is useful when applied to $\Delta^{\bar{c}},
\Delta^{\bar{a}}$ and $\Delta^{\bar{a}^{\prime}}$.

\begin{lemma} \label{estra}
Suppose we have a $($1B$)$ vertex with exactly $k$ adjacent $($1B$)$ vertices.
Then
\[
 r \leq \# ( (1A) \ {\rm vertices} ) + \#( (2) \ {\rm vertices}) + k + 2.
\]
Suppose we have a $($1B$)$ vertex with no adjacent $($1B$)$ vertices.
Then
\[
r \leq \# ( (1A) \ {\rm vertices} ) + \#( (2) \ {\rm vertices}) + 2.
\]
If there are no $($1B$)$ vertices then
\[
r \leq \# ( (1A) \ {\rm vertices} ) + \#( (2) \ {\rm vertices}) +1.
\]
\end{lemma}

\begin{proof}
By our assumption that $\dim {\rm conv}(\mathcal W)= r-1$, it follows
that  $\Delta^{\bar{c}}$ is a polytope of dimension $r-2$. Any vertex in it
has at least $r-2$ adjacent vertices. So for a (1B) vertex, the
first two statements follow immediately. If there are no (1B) vertices
the third inequality follows because $\Delta^{\bar{c}}$ has at least $r-1$ vertices.
\end{proof}

\begin{lemma} \label{Tr110}
Configurations $($Tr1$)$ - $($Tr22$)$ cannot arise from adjacent
$($1B$)$ vertices.
\end{lemma}

\begin{proof}
The strategy is to count the number of type (1A), (2) and
adjacent (1B) vertices in $\Delta^{\bar{a}}$ or $\Delta^{\bar{c}}$
and apply Lemma \ref{estra} to get a contradiction.

(i) We first observe that for these configurations $c$ and $a$ have at
least four nonzero entries (at least five except for (Tr2)),
so they cannot be collinear with an edge $vw$ with points
of $\mathcal W$ in the interior of $vw$ (see Table 3 in \S 3).
So if $\Delta^{\bar{c}}$ or $\Delta^{\bar{a}}$ has a type (2) vertex,
by Theorem \ref{cedge}, $c$ or $a$ must equal $2v-w$ or $(4v-w)/3$.
It is easy to check that this is impossible except for $c$ in (Tr3),
using the forms of $c$ in Tables 1, 2 in \S 3.

(ii) Next we consider (1A) vertices. For (Tr1)-(Tr10) we have
$| a_i/d_i | \leq \frac{1}{3}$ for all $i$. For Tr(1), Tr(6),Tr(8)
there are three $i$ where equality holds. In these cases one of the
associated $d_i$ equals $1$. Moreover for (Tr1) and (Tr8) two
of these $a_i/d_i$ equal $1/3$ and the third is $-1/3$, wheras for (Tr6)
it is the other way round. For (Tr2)-(Tr5), (Tr7) and (Tr9)-(Tr10)
there are only two $i$ where equality holds. Further, for (Tr2) and
(Tr5) $|c_i / d_i | \leq \frac{1}{3}$ for all $i$, with equality for
just two $i$, and here $c_i/d_i = \frac{1}{3}$. It follows that for
(Tr1), (Tr8) there are at most two (1A) vertices in $\Delta^{\bar{a}}$,
while for (Tr2)-(Tr5), (Tr7) and (Tr9)-(Tr10) there is at most one.
In the case of (Tr6) there are at most three (1A) vertices in general
but at most two if $K$ is connected. For (Tr2) and (Tr5) there are
no (1A) vertices in $\Delta^{\bar{c}}$.

By means of similar considerations, we find that there is one (1A)
vertex (corresponding to $\bar{x}^{\prime \prime}$) in
$\Delta^{\bar{a}}$ for (Tr12), (Tr13),(Tr14) (Tr16), (Tr18), (Tr19)
and at most two (1A) vertices for (Tr11), (Tr17), (Tr20) and the
$d=(1,1,3,6,6,6,\cdots)$ case of (Tr21), (Tr22).
For (Tr15) there are at most four (1A) vertices in $\Delta^{\bar{a}}$, and
for the remaining cases of (Tr21) and (Tr22) there are at most
$r-4$ (1A) vertices ($r-6$ of those correspond to $(-2^3, 1^j)$ where
$j>6$).

(iii) Finally, consider the (1B) vertex $\bar{\xi}$ in $\Delta^{\bar{a}}$
corresponding to $\bar{x}$ in each of the triangles. In order for there
to be an adjacent (1B) vertex, $\bar{a}$ must be (up to permutation)
of the form of the null vector $\bar{c}$ in the $2$-faces in
Theorem \ref{1B-initial}. Now observe that for examples (Tr1), (Tr3),
(Tr4), (Tr6)- (Tr20) the null vector $a$ does not appear in the list
of possible $c$. Hence $\bar{\xi}$ has no adjacent (1B) vertices.
{}From above, type (2) vertices cannot occur, so combining the bounds
for (1A) vertices in $\Delta^{\bar{a}}$ with Lemma \ref{estra} gives
an upper bound for $r$ less than the minimum required by each
configuration, a contradiction.

(iv) Let us now consider (Tr21) and (Tr22). The vector $a$ of (Tr21) has the
same form as $c$ in (Tr22) and vice versa. An adjacent $2$-face containing
$a\xi c$ of (Tr21) can only be a triangle of type (Tr22) containing
$\frac{1}{3}(1, -1, -3, -2, -2, 4, 0, \cdots)$. Thus $\xi$ has at most one
adjacent (1B) vertex, and we get the bound $r-2 \leq 2 + 1 + 0$ in the
$d=(1, 1, 3, 6, 6, 6, \cdots)$ subcase and $r-2 \leq (r-4)+1+0$ in the
other two subcases, a contradiction.

(v) For the remaining two triangles (Tr5) and (Tr2) we consider $\Delta^{\bar{c}}$
instead.

For (Tr5), observe that $c$ determines the plane $x x^{\prime}
x^{\prime \prime}$ and does not occur as a possible null vector for
any other configurations. So we have at most one adjacent pair of
(1B) vertices in $\Delta^{\bar{c}}$. From above, there are no type
(1A) or (2) vertices. But $r \geq 5$, giving a contradiction.

For (Tr2), consider the vertex $\bar{\xi}^{\prime}$ of $\Delta^{\bar{c}}$
corresponding to $x^{\prime}$. If there is a (1B) vertex adjacent to
it, we have a 2-dim face including $c, x^{\prime}$. By Theorem \ref{1B-initial}
the only one is the face including $x$, so there is at most one (1B)
vertex of $\Delta^{\bar{c}}$ adjacent to $\bar{\xi}^{\prime}$. Also,
type (1A) and (2) vertices cannot occur, so $r \leq 3$, a contradiction.
\end{proof}

\begin{lemma} \label{Tr1115}
Configurations $($Tr23$)$-$($Tr27$)$ cannot arise from adjacent
$($1B$)$ vertices.
\end{lemma}

\begin{proof}
Note first that all entries of $c, a, a^{\prime}$ are integers, so
Lemma \ref{Zperp} shows in each case there is at most one (1A) vertex, and
for (Tr23),(Tr24) one checks that there are no (1A) vertices in $\Delta^{\bar{c}}$.
Note also that for all these configurations, as $x^{\prime \prime}$
is a vertex, there are no type II vectors with nonzero entry in place 1.

Observe as in Lemma \ref{Tr110} that there are no type (2)
vertices in $\Delta^{\bar{c}}, \Delta^{\bar{a}}$ or
$\Delta^{\bar{a}^{\prime}}$. (For (Tr24), we need to rule out the
possibility that $c$ has the form (4) in Table 3 in \S 3
with $\lambda=3/2$. This follows since the interior point in that
case would be a type II vector with nonzero entry in place 1.)

So in all cases if we have a (1B) vertex with exactly $k$
(1B) vertices adjacent to it, then by Lemma \ref{estra} we have $r \leq k+3$.
For (Tr23), (Tr24) (using $\Delta^{\bar{c}}$), we have $r \leq k+2$.
We will work with $\Delta^{\bar{c}}$ below.

First consider (Tr23) and look for (1B) vertices adjacent to
$\bar{\xi}$ where $\bar{\xi}$ corresponds to the vertex $x$.
We need a $2$-face including $c, x$. By Theorem \ref{1B-initial} such
a face must be of type (Tr23), and having fixed $c$ and $a$, the only
freedom lies in assigning $1$ in the third null vector to the first
or fifth place. So $k \leq 2$, which gives $r \leq 4$, a contradiction.

Similarly, for (Tr24), since a type (H1) face cannot contain $c$ and $x$,
we need only consider faces of type (Tr24), for which there are again
two possibilities. However, as mentioned above, in one of these possibilities the
vector ``$x^{\prime}$'' has a 1 in place 1 and hence cannot occur.  So there is
at most one (1B) vertex adjacent to $\xi$, and  we deduce
$r \leq 3$, a contradiction.

For (Tr25),(Tr27) we similarly deduce that the only  $2$-face
containing $x,c$ is itself because as above we cannot have any type
II vectors with nonzero entry in place 1.  So $r \leq 4$, a
contradiction.

Finally, for (Tr26) the above argument still works since (Tr24) has
been ruled out (the vectors $a, a^{\prime}$ of (Tr24) are of the same
form as $c, a$ of (Tr26), so a priori (Tr24) could be an adjacent $2$-face).
\end{proof}

\begin{lemma} \label{S}
Configuration $($S$)$ $($square with midpoint$)$ cannot arise from
adjacent $($1B$)$ vertices.
\end{lemma}

\begin{proof}
We refer to \S 6 for the expressions for the vertices $vusw$ of the square.
The null vertex $\bar{c}$ corresponds to $(1, -1, -1, 1, -1, 0,\cdots)$
and the $2$-dimensional face is cut out by $x_2 =-1, x_1 + x_3=0=x_4 + x_5,$
and $x_k=0,$ for $ k > 5.$  Lemma \ref{Zperp}(b) shows that there
is at most one (1A) vertex in $\Delta^{\bar{c}}$. As $r \geq 5$ and
all the nonzero entries of $c$ have the same absolute value, it follows
that there are no type (2) vertices.  Let $\bar{\xi}$ denote the vertex of
$\Delta^{\bar{c}}$ so that $\xi$ is collinear with $u$ and
$a=2u-c=(-1, -1, 1, 1, -1, 0, \cdots)$.
A (1B) vertex adjacent to $\bar{\xi}$ gives a $2$-dimensional face including
$\bar{c}, \bar{u}$. By what we have analysed so far about $2$-faces given by
adjacent (1B) vertices, this face must again be a face of type (S),
and the only possibilities are itself or the face obtained from this
by swapping indices $2$ and $5$. Hence there are at most two (1B)
vertices adjacent to $\bar{\xi}$, and at vertex $\bar{\xi}$, we have
$3 \leq r-2 \leq 1+2$. Thus $r=5$ is the remaining possibility, in which case
$\bar{\xi}$ has exactly two adjacent (1B) vertices and one adjacent
(1A) vertex.

Let us denote by $\bar{\xi}^{\prime}$ the (1B) vertex such that $\xi^{\prime}$ is
collinear with $w$ and $a^{\prime}:=2w-c=(1, -1, -1, -1, 1)$. Let $\bar{\eta}$
denote the other (1B) vertex adjacent to $\bar{\xi}$.
 Then the $2$-face determined by $c,
\xi, \eta$ is cut out by $x_5=-1, x_1+x_3 =0=x_2+x_4.$ The ray
$c\eta$ intersects ${\rm conv}({\mathcal W})$ at $z=(1, 0, -1, 0, -1)$
and $b:=2z-c=(1, 1, -1, -1, -1)$ corresponds to a null vertex. Similarly,
there is a (1B) vertex $\bar{\eta}^{\prime}$ (besides $\bar{\xi}$) adjacent to
$\bar{\xi}^{\prime}$, and the corresponding $2$-face (also of type (S)) is
cut out by $x_3 =-1, x_1+x_2=0=x_4+x_5$. The ray $c\eta^{\prime}$ intersects
${\rm conv}(\mathcal W)$ at $z^{\prime}=(0, 0, -1, 1, -1)$. The vector
$b^{\prime} :=2z^{\prime}-c =(-1, 1, -1, 1, -1)$ corresponds to a null vertex.

Let us examine the (1A) vertex in $\Delta^{\bar{c}}$ more closely. Let
$y \in {\mathcal W}$ such that $J(\bar{y}, \bar{c}) = 0$. As $r=5$,
the null condition for $\bar{c}$ implies that $d_i \geq 2$ with at most
one equal to $2$. Also, for some $j \in \{2, 3, 5\}$ (i.e., $j$ is an
index for which the corresponding entry of $c$ is $-1$) we must have
$y_j =-2,$ so $y$ is type III. Let $i$ be the index such that $y_i =1$.
Then $i \in \{1, 4\}$ (i.e., $i$ is an index for which the corresponding
entry of $c$ is $1$), and the orthogonality condition implies
$(d_i, d_j) = (2, 4)$ or $(3, 3)$. There are thus six possibilities for $y$,
but only one can actually occur.

With the possible exception of the existence of the (1A) vertex,
the above arguments apply equally to the projected polytopes
$\Delta^{\bar{a}},$ $\Delta^{\bar{b}},$ $\Delta^{\bar{a}^{\prime}}$,
and $\Delta^{\bar{b}^{\prime}}$ as the entries of $a, b, a^{\prime}$
and $b^{\prime}$ are just permutations of those of $c$. We claim
that whichever possibility for $y$ occurs in $\Delta^{\bar{c}}$,
there is another projected polytope with no (1A) vertex. Applying the
above arguments to this polytope would result in the contradiction
$r-2 \leq 2$ and complete our proof.

We can use $\Delta^{\bar{a}}$ for the contradiction if $d_1 = 2$
or if any of $(d_1, d_2), (d_3, d_4), (d_1, d_5)=(3, 3)$. If $d_4 =2$
or if $(d_2, d_4)=(3, 3)$ we can use $\Delta^{\bar{a}^{\prime}}$ instead.
Finally, if $(d_1, d_3)=(3, 3)$ we can use $\Delta^{\bar{b}^{\prime}}$
and if $(d_4, d_5)=(3, 3)$ we can use $\Delta^{\bar{b}}$.

For example, when $(d_1, d_3)=(3, 3)$ (so $y=(1^1, -2^3)$), the null
condition for $\bar{c}$ implies that $d_2, d_4$ in particular cannot
equal $2$ or $3$. In order to have a (1A) vertex in $\Delta^{\bar{b}^{\prime}}$,
we must have a type III vector $(1^i, -2^j)$ with $i \in \{2, 4\}$.
But this requires one of $d_2, d_4$ to be $2$ or $3$.

When $(d_4, d_5)=(3, 3)$, then $d_1, d_2$ cannot be $2$ or $3$. But in
$\Delta^{\bar{b}}$ a (1A) vertex corresponds to $(1^i, -2^j)$ with
$i \in \{1, 2\}$, which implies that one of $d_1, d_2$ is $2$ or $3$.

The remaining cases are handled similarly.      \end{proof}

\begin{lemma} \label{subH2}
The subrectangle $y y^{\prime} z z^{\prime}$ of $($H2$)$ cannot arise from
adjacent $($1B$)$ vertices.
\end{lemma}

\begin{proof}
Recall $c = (-2,1,0, \cdots)$, so by Lemma \ref{Zperp} there are no
(1A) vertices of $\Delta^{\bar{c}}$. Moreover, using Tables 1-3 in \S 3,
one may check that there are no type (2) vertices either.
(Note that  type II vectors other than $y,z$ with a nonzero entry in place $1$
cannot occur as then the subrectangle cannot be a face. Similarly the line
through $c, \alpha, \beta$ and $(1,-2, 0, \cdots)$ will not give a type (2)
vertex as this line cannot be an edge.)

Let $\bar{\eta}$
denote the vertex of $\Delta^{\bar{c}}$ collinear with $c$ and $y$.
Any (1B) vertex adjacent to $\bar{\eta}$ will give rise to a face
containing $c$ and $y$, which cannot be of type (H1), and must therefore
be of type (H2), since we have eliminated all other possibilities.
In fact, it must be the face we started with.

So there is just one (1B) vertex adjacent to $\bar{\eta}$, and from above
there are no (1A) or (2) vertices. As $r \geq 4$ for (H2), this
contradicts Lemma \ref{estra}.
\end{proof}

\begin{lemma} \label{hexH1}
If $r \geq 4$, configuration $($H1$)$ or subshapes cannot arise from
adjacent $($1B$)$ vertices.
\end{lemma}

\begin{proof}
We first note some special properties of $\mathcal W$.
Since (H1) is a face, there can be no type II vectors in $\mathcal W$
with nonzero entry in a place $\in \{1,2,3 \}$ and in a place
$\notin \{1,2,3 \}$. Also, if $(-2^i, 1^k)$ with
$i \in \{1,2,3 \}, \; k \notin \{1,2,3 \}$, then $(-1^k)$ must be absent,
which has strong implications, as noted in Remark \ref{facts}(b).

Let $\bar{\xi}$ be a (1B) vertex in the plane. We have at most one (1B)
vertex adjacent to $\bar{\xi}$, as the associated face must again be of
type (H1) and is now determined by $\bar{c}, \bar{\xi}$.
It also readily follows that $c$ cannot be collinear with an edge of
$\mathcal W$ not in the face (assuming as usual we are not in
the situation of Theorem \ref{Fano}).

Now the special properties of $\mathcal W$
in the first paragraph imply that (1A) vertices in
$\Delta^{\bar{c}}$ can correspond only to type III  vectors in $\mathcal W$
which overlap with $c$. A straight-forward check using the null condition
for $\bar{c}$ shows that the possible type IIIs have form $(-2^i, 1^k)$ with
$i \in \{1, 2, 3\}, \; k \notin \{1, 2, 3\}$ and $c_i/d_i = -1/2$.
It follows that $d_i =2$ or $3$ and hence, by nullity, the index $i$ is unique.
So there are at most $r-3$ (1A) vertices. By Lemma \ref{estra},
all $r-3$ (1A) vertices must occur. Applying Cor \ref{orthogonal}
we conclude that $r \leq 4$ (as $d_i \neq 4$ and $r >4$ forces
$i \in {\hat S}_{\geq 2}$).

We will now improve this estimate to $r \leq 3$.
Let the vertices of (H1) be as in \S 7.  If $r=4$ then a (1A) vertex
does exist and we can take it to come from $t =(-2,0,0,1)$ with
$(-1^4)$ absent. It follows that besides $t$ the only other possible
members of $\mathcal W$ lying outside the $2$-plane containing (H1) are
$(0, -2, 0, 1)$ and $(0, 0, -2, 1)$. As noted just before Theorem \ref{triangle}
we may assume the type I vectors $\alpha, \beta, \gamma$ are
all present. (If $K$ is connected, $d_4 =1$ and so this last fact
follows without having first to eliminate those subshapes not
containing one of the type I vectors.)

As noted above, $d_1 = 2$ or $3$, and $c_1 =-1$ or
$-\frac{3}{2}$ respectively.

First consider $c_1=-1$, so $d_1=2$.  Now $c = (-1, c_2, -c_2, 0)$,
and by swapping the $2,3$ coordinates if necessary, we may
take $c_2 >0$. Observe $u, q$ are absent, as if $u$ is present or
if $u$ is absent but $q$ is present, then $u$ (resp. $q$) gives a (1B)
vertex, which contradicts nullity as the associated $a$ would have $a_1 = -3$.
Now the type II vectors $x,y,z$ are absent, as if one is present they
all are, and we have a type (2) vertex. We deduce $\alpha$ gives
a (1B) vertex so $a = (-1,-c_2, c_2, 0)$.  The other (1B) vertex
cannot correspond to $w$ since $\beta, \gamma$
are present. It also cannot be given by $v, s$ as this violates nullity,
so must correspond to $p$ or $\beta$.

If it is $p$, we have  $a^{\prime} = (1, 2-c_2, c_2-4,0)$.
Now Remark \ref{xc} implies $c_2 = (d_3 + 4d_2)/(d_3 + 2d_2)$ so $1 < c_2 < 2$.
But now no entry of $a^{\prime}$ equals $-1$ or $-\frac{3}{2}$. We can
now check that there are no (1A) or (2) vertices with respect to $a^{\prime}$,
so there is at most one vertex of $\Delta^{\bar{a}^{\prime}}$ adjacent to
$p$, a contradiction.

If it is $\beta$ then $p, v$ must be absent. Now
$a^{\prime} = (1, -c_2, c_2-2)$ and Remark \ref{xc} implies $c_2=1$.
Hence  $c =(-1,1,-1), \; a=(-1,-1,1), \; a^{\prime}=(1,-1,-1)$, and
nullity implies $\frac{1}{d_2} + \frac{1}{d_3} = \frac{1}{2}$. It is easy to
check by considering the vertices of $\Delta^{\bar{a}^{\prime}}$ that
 $w,s$ must also be absent, so $\mathcal W$ just contains the three
type I vectors, $t$ and possibly one or both of $(0, -2, 0, 1), (0, 0, -2, 1)$.
But we can check that, if present, these three latter vectors give respectively
vertices with respect to $a, a^{\prime}$  which cannot satisfy any of the
conditions (1A), (1B) or (2). So in fact we have $r=3$.

Similar arguments rule out the case $c_1 = -\frac{3}{2}$.
\end{proof}

\medskip
Lemmas \ref{Tr110}-\ref{hexH1} give the following improvement of Theorem
\ref{1B-initial}.

\begin{theorem} \label{adj1B}
It is impossible to have adjacent $($1B$)$ vertices except possibly
when $r=3$, in which case ${\rm conv}(\frac{1}{2}(d+{\mathcal W}))$
is a proper subface of $($H1$)$ containing all three type I vectors
$($e.g., the tri-warped example $($Tr28$)$$)$.  \; $\qed$
\end{theorem}

We are now in a position to strengthen Theorem \ref{two2} by eliminating
the remaining case of the pentagon.

\begin{theorem} \label{single2}
Let $\bar{c}$ be a null vertex of ${\rm conv}(\mathcal C)$ such that
$\Delta^{\bar{c}}$ contains more than one type $($2$)$ vertex. Then we are in
the situation of Theorem \ref{Fano}.
\end{theorem}
\begin{proof}
We just have to eliminate the case of the pentagon $uy y^{\prime} z z^{\prime}$
in Theorem \ref{two2}. Recall $r \geq 4$ for this configuration,
and $c$ is $(0,-1,2,-2,\cdots)$.

Using the nullity of $\bar{c}$ we check that the only elements of
$\mathcal W$ which can give an element of $\bar{c}^{\perp}$ are
$(-2^2, 1^i)$ where $i > 4$ and we have $d_2 =2$. Note that $(1^1, -2^2)$
cannot be present as then $y^{\prime} z^{\prime}$ is not an edge.
By Cor \ref{orthogonal}, at most one such vector can arise.
So there is at most one (1A) vertex, which occurs only if $r \geq 5$.

If we can show there are no (1B) vertices, then we are done because
if we look at the adjacent vertices of the type (2) vertex associated
to $y$ (in the pentagon), besides one (1A) possibility, the
other possibility is the type (2) vertex associated to $y^{\prime}$
(by Theorem \ref{two2}). As there must be at least $r-2$ adjacent
vertices, we deduce $r-2 \leq 2$, so $r=4$. But now, from above there
is no (1A) vertex so in fact we get $r-2 \leq 1$, a contradiction.

We now use Remark \ref{xc} to make a list of the possible $x \in {\mathcal W}$
associated to (1B) vertices of $\Delta^{\bar{c}}$. These are
$(0, 1, -1, -1, 0, \cdots), (0, -2, 1, 0, \cdots), (1^3, -1^i, -1^j), (-1^4, 1^i, -1^j),
(1^3, -1^4, -1^i),$ $(-1^3, -1^4, 1^i)$ and $(1^3, -2^i)$ where
$i, j \neq 2, 3, 4$. Note that type II vectors with nonzero entries
in places $2,k,m$ cannot occur except for $y^{\prime}, z^{\prime}$
as then $y^{\prime} z^{\prime}$ is not an edge.
For each $x$ in this list, we consider the projected
polytope $\Delta^{\bar{a}},$ where $a=2x-c$. By looking at the form of $a$,
we see from Theorem \ref{two2} that there is at most one type (2) vertex in
$\Delta^{\bar{a}}$. Also, the nonzero components of $a$ are either $\geq 1$ or
$\leq -2$. By Lemma \ref{Zperp}(a), there are no vertices of type (1A) in
$\Delta^{\bar{a}}$. Since $r \geq 4$, by Theorem \ref{adj1B}, the type (1B)
vertex in $\Delta^{\bar{a}}$ corresponding to $\bar{x}$ has no adjacent (1B)
vertices. So we have a contradiction to Theorem \ref{estra}.
\end{proof}

\smallskip
The above result together with Theorem \ref{adj1B} and Lemma \ref{estra}
gives us lower bounds on the number of (1A) vertices.

\begin{theorem} \label{estr}
Let $\bar{c}$ be a null vertex of ${\rm conv}({\mathcal C})$ and
$\Delta^{\bar{c}}$ be the corresponding projected polytope. Suppose
further that $c$ is not type I, i.e., we are not in the case of
Theorem \ref{Fano}.

$($i$)$ If there are no $($1B$)$ vertices in $\Delta^{\bar{c}}$, then
there are at least $r-2$ type $($1A$)$ vertices.

$($ii$)$ If either there is a type $($2$)$ vertex or $r \geq 4$, then
there are at least $r-3$ type $($1A$)$ vertices in $\Delta^{\bar{c}}$.
Hence there are at least $r-3$ elements of $\frac{1}{2}(d + {\mathcal W})$
orthogonal to $\bar{c}$.  \  $\qed$
\end{theorem}

\section{\bf Type (2) vertices}

In this section we consider again type (2) vertices of $\Delta^{\bar{c}}$.
In view of Theorem \ref{single2}, it remains to deal with the case of a unique
type (2) vertex in $\Delta^{\bar{c}}$. By Theorem \ref{adj1B}
there are no adjacent (1B) vertices in this situation.

Let $c$ be collinear with an edge $vw$ of ${\rm conv}(\mathcal W)$.
We first consider the situation where there are no interior points of $vw$
lying in $\mathcal W$. By Theorem \ref{cedge}, we have the two possibilities
$c=2v-w$ and $c=(4v-w)/3$. Moreover, a preliminary listing of the cases appears
in Tables 1 and 2 of \S 3.

\smallskip

{\noindent \bf Case (i):} $c = 2v-w$

We have to analyse cases (1)-(7) in Table 1 of \S 3. The idea is to
determine the number of (1A) and (1B) vertices using respectively
Lemma \ref{Zperp} and Remark \ref{xc}, and then get a contradiction
(sometimes using Theorem \ref{estr}). Note that $J(\bar{w}, \bar{w}) <0$
for (1)-(3).

In (1), (2) and (4)-(7), Lemma \ref{Zperp} shows that there are no
elements of $\frac{1}{2}(d + {\mathcal W})$ orthogonal to $\bar{c}$
(recall $c \notin {\mathcal W}$), so Theorem \ref{estr} shows that $r \leq 3$.
This already gives a contradiction in case (7). (Note that when $r=3$
and $w$ is type I, since $w$ is a vertex there are no type II vectors in
$\mathcal W$.)

In (1) the only $x \in {\mathcal W}$ that could satisfy Eq.(\ref{xceqn})
and give a (1B) vertex with respect to $\bar{c}$ is
$(1,-1,-1)$. But the associated $a=2x-c$ is $(2,-3,0)$ and it easily
follows that $\bar{a},\bar{c}$ cannot both be null.  For (2), the
possible $x \in {\mathcal W}$ which correspond to (1B) vertices are
$(1,-2,0)$ and $(1,-1,-1)$ respectively. In each case we
find the nullity of $\bar{c}$ and Remark \ref{xc}
imply $J(\bar{w}, \bar{w}) >0$, a contradiction to
$F_{\bar{w}}^2 \ J(\bar{w}, \bar{w})=A_w < 0$ (as $w$ is type III).

In (5) with $r=3$, one checks that the only possible $x \in {\mathcal W}$
corresponding to a (1B) vertex is $(0,1,-2)$.
Let us consider the distribution of points of $\mathcal W$
in the plane $x_1 + x_2 + x_3 = -1$. The point $(0, 1, -2)$, if present,
would lie on one side of the line $vw$ while the point $(-1, 0, 0)$
lies on the other side. Now $(-1, 0, 0)$ must lie in $\mathcal W$
as otherwise $v$ cannot be present by Remark \ref{facts}(b). So
since $vw$ is an edge by assumption, $(0, 1, -2)$ cannot lie in $\mathcal W$,
which gives a contradiction to Theorem \ref{estr}(i).


Hence in (1),(2),(5) Theorem \ref{estr} shows $r \leq 2$, which
is a contradiction.

In case (4) the nullity of $\bar{c}$ translates into $1= 9/d_1 + 4/d_2$.
Hence $d_2 \neq 1$, so if $K$ is connected $(0,-1, 0, \cdots)$ is present
and $w$ is not  a vertex, which is a contradiction.
 If $r=3$ and $K$ is not connected, by Remark \ref{facts}(b),
$(1, -2, 0)$ and $(0, -2, 1)$ must be absent, and, from Remark \ref{xc},
the possibilities for $x \in {\mathcal W}$ associated to the (1B) vertex
are $x=(-2, 0, 1)$ and $y=(0, 1, -2)$. In the first case,
${\rm conv}(\mathcal W)$ is the triangle with vertices $v, w, x$
and $a=2x-c=(-1, -2, 2)$.
Now $J(\bar{a}, \bar{w})>0$, contradicting the superpotential equation.
In the second case, $a=(3, 0, -4)$ with $J(\bar{a}, \bar{y}) > 0,$
and $aw$ intersects ${\rm conv}(\mathcal W)$ in an edge.
By Theorem \ref{cedge}, $t=(1, 0, -2) \in {\mathcal W}$ and
${\rm conv}(\mathcal W)$ is a parallelogram with vertices $v,y,w,t$.
Moreover, Remark \ref{2int} implies that $a$ and $w$ are the only elements
of $\mathcal C$ in $aw$. But then the midpoint $(0, 0, -1)$ of $wt$
is unaccounted for in the superpotential equation.

For (6) with $r=3$, there should be at least two vertices in $\Delta^{\bar{c}}$.
But we find there are no (1B) vertices, a contradiction. So $r=2$, and we are in
the situation of the double warped product Example 8.2 of \cite{DW4}.

In case (3), Lemma \ref{Zperp} shows $\frac{1}{2} (d + {\mathcal W}) \cap
\bar{c}^\perp$ has at most one element.  Hence $\Delta^{\bar{c}}$, which has
dimension $\geq 2$ since $r \geq 4$, must contain at least one (1B) vertex.
By Theorem \ref{adj1B}, such a (1B) vertex has at most 2 adjacent vertices.
It follows that $r=4$ and $(1, -2, 0, 0)$ corresponds to the (1A) vertex;
also $d_2=2$. Also, since $(-1, 0, 0, 0) \in {\mathcal W}$,
$(0,-1,0,0)$ cannot be a vertex of ${\rm conv}(\mathcal W)$. But now routine
computations using Eq.(\ref{xceqn}) show there are no (1B) vertices, a
contradiction.

So the only possible case if $K$ is connected is that giving Example 8.2
of \cite{DW4}. If $K$ is disconnected there is the further possibility of (4)
with $r=2$, i.e., $\mathcal W = \{ (-2,1), (-1,0) \}$. This is discussed
in the third paragraph of Example 8.3 of \cite{DW4}. An example in the inhomogeneous
setting is treated there and in \cite{DW2}. An example where the hypersurface
is a homogeneous space $G/K$ is discussed in the concluding remarks
at the end of section \S 10.

\smallskip

{\noindent \bf Case (ii):} $c = (4v -w)/3$

For clarity of exposition let us assume $K$ is connected, using the
assumption as indicated in Remark \ref{Kconn}. We examine the
cases (1)-(11) in Table 2 of \S 3.

Some of these cases can be immediately eliminated. In (3),
Eq.(\ref{4/3cond}) implies $(d_1, d_2) = (3,3)$ or $(4,1)$. In neither
case is $\bar{c}$ null. In (11) Eq.(\ref{4/3cond})
and $J(\bar{v}, \bar{w})>0$ imply $(d_1, d_2, d_3)
= (2,4,4), (2,5,2)$ or $(3,3,3)$, and again $\bar{c}$ is not null.
In (4) and (6) Eq.(\ref{4/3cond}) implies $(d_1, d_2) =(2,1)$ and
$(d_1, d_2) = (3,9)$ or $(4,3)$ respectively. In neither case does the
nullity condition have an integral solution in $d_3$.

Further cases can be eliminated by finding the possible (1A) vertices
(using Lemmas \ref{IX6} and \ref{Zperp}) for the given value of $c$ and
using Theorem \ref{estr}. In particular, we get a contradiction
whenever $r \geq 4$ and there are no (1A) vertices.

In (1), Eq.(\ref{4/3cond}) implies $(d_1, d_2)=(2,1),$ and nullity of $\bar{c}$
implies $\frac{32}{d_3} + \frac{8}{d_4}=3$. But we now find that
$\frac{1}{2}(d + {\mathcal W})\cap \bar{c}^\perp$ is empty, giving a
contradiction as $r \geq 4$.

In (5) Eq.(\ref{4/3cond}) and nullity imply $(d_1,d_2)=(3,3)$ and $\{d_3,
d_4 \}$ = $\{3, 8 \}$.  One can now check that $\frac{1}{2}(d +
{\mathcal W}) \cap \bar{c}^\perp$ is empty, which is a contradiction
as $r \geq 4$.

In (7), Eq.(\ref{4/3cond})  implies $(d_1,d_2)=(2,1)$ and nullity
implies $\frac{1}{d_3} + \frac{1}{d_4} + \frac{1}{d_5} = \frac{3}{8}$.
One can now check that the only possible elements $\bar{u}$ orthogonal
to $\bar{c}$ correspond to $u= (1,0,0,-2,\cdots)$ if $d_4=4$ and
$(1,0,0,0,-2,\cdots)$ if $d_5=4$. The nullity condition means that at
most one of these can occur. This is a contradiction as $r \geq 5$.

In (8) Eq.(\ref{4/3cond}) gives $(d_1,d_2)=(2,3)$ and nullity of
$\bar{c}$ gives $\frac{1}{d_3} + \frac{1}{d_4} = \frac{1}{4}$.
Again one can check that  $\frac{1}{2}(d + {\mathcal W}) \cap \bar{c}^\perp$
is empty, a contradiction as $r \geq 4$.

In (9), Eq.(\ref{4/3cond}) and the nullity of $\bar{c}$ give
$(d_1,d_4)=(2,16)$ and $\{d_2, d_3 \} = \{ 2,3 \}$.
 The only $u$ which can give $\bar{u} \in \bar{c}^\perp$ are
$(0,-2,0,0,1^i,\cdots)$ if $d_2=2$ or $(0,0,-2,0,1^i,\cdots)$ if
$d_3 =2$, where $i \geq 5$. In each case, $i$ is unique since $d_2$
(resp. $d_3$) $\neq 4$. Since $r \geq 4$,  Theorem \ref{estr} now implies
$r = 4$. But now these $u$ are not present (as $i \geq 5$). Therefore
there are actually no (1A) vertices, a contradiction to $r=4$.

In (10) we have $(d_3,d_4)=(2,16)$ and $\{d_1, d_2 \} = \{ 2,3 \}$. The
only $u$ which can give an element of $\bar{c}^\perp$   is
$(0,-2,0,0,1^i,\cdots)$ (for $i$ unique and $\geq 5$) if $d_2=2$.
The final argument in (9) now applies equally here.

Finally, we can eliminate (2) by an analysis of  both the (1A) and (1B)
vertices. First,  Eq.(\ref{4/3cond}) and nullity of $\bar{c}$ force
$(d_1,d_2, d_3) =(6,1,8)$. Next we check that
$\frac{1}{2}(d + {\mathcal W}) \cap \bar{c}^\perp$ is empty,
so $r=3$. Using Remark \ref{xc} we then find there
can be no (1B) vertices, giving a contradiction.

\medskip
So case (ii) cannot occur if $K$ is connected.

\begin{rmk} \label{disconn}
Case (ii) is the only part of this section that relies on the
connectedness of $K$. In fact, the analysis of the cases where
$w$ is type III does not use this assumption. If $K$ is not connected,
using the same methods and with more computation we obtain the
following additional possibilities (all of which are associated to a
$w$ of type II).
\[
\begin{array}{|c|c|c|c|c|c|} \hline
    &         v       &      w     &    c=(4v-w)/3   &    d   &  r   \\ \hline
(9*) & (0, -1, -1, 1) & (1, -1, -1, 0)  &  (-\frac{1}{3}, -1, -1, \frac{4}{3}) &  (1, 2, 6, 8) & 4,5 \\
     &                   &              &            &   (1, 6, 2, 8) & 4, 5\\
(10*) & (1, -1, 0, -1) & (1, -1, -1, 0) & (1, -1, \frac{1}{3}, -\frac{4}{3}) &   (3, 3, 1, 8) & 4 \\
      &              &                 &              &    (6, 2, 1, 8)  & 4,5 \\
(14) & (0, -1, 0, 1, -1) & (1, -1, -1, 0, 0) &   (-\frac{1}{3}, -1, \frac{1}{3}, \frac{4}{3}, -\frac{4}{3}) &
      (1, 3, 1, d_4, d_5)  &  5\\ \hline
\end{array} \]

In (9*) and (10*) there is always a (1B) vertex in $\Delta^{\bar{c}}$,
and $r=4$ or $5$ according to whether the cardinality of
$\bar{c}^{\perp} \cap \frac{1}{2}(d+ \mathcal W)$ is $1$ or $2$.
The dimensions $d_4, d_5$ in (14) must satisfy
$\frac{1}{4}= \frac{1}{d_4} + \frac{1}{d_5}$ (i.e., $\{d_4, d_5\} =
\{5, 20\}, \{6, 12\}, \{8, 8\}$) and again there is always a (1B)
vertex in $\Delta^{\bar{c}}$.
\end{rmk}

\smallskip
{\noindent \bf Interior points}

\smallskip
Finally, we must consider the cases, listed in Table 3, \S 3,
when there may be points of $\mathcal W$ in the interior of $vw$.
As in the earlier cases, we analyse the possible (1A) and (1B)
vertices for these $c$.

For (1) and (2), as $1 < \lambda \leq 2$, the nonzero entries of $c$ are either
$< -2$ or $> 1$. Hence by Lemma \ref{Zperp} there are no (1A) vertices. By
Theorem \ref{estr} we have $r=2$ or $3$.

In case (3) the nullity of $\bar{c}$ implies that a vector
$u \in \mathcal W$ not collinear with $vw$ and with
$\bar{u} \in \bar{c}^\perp$ must be of the form $(-2,0, 1^j)$.
 So $2 \lambda-1 = \frac{d_1}{2}$ and
$c = (-\frac{d_1}{2}, -1 + \frac{d_1}{2},0,\cdots)$. From the
range for $\lambda$ and the nullity of $\bar{c}$, we have $d_1 =3, \ d_2 =1$.
But $d_2 \neq 1$ since $w \in {\mathcal W}$. So again there are no
(1A) vertices and by Theorem \ref{estr} we have $r \leq 3$.

In case (4), a straight-forward preliminary analysis reduces the
possibilities of $u \in {\mathcal W}$ such that $\bar{u} \in
\bar{c}^\perp$ to the choices $u = (-2,0,1,\cdots), (0,-2,1,\cdots)$
or $(-1,-1,1,\cdots)$. Note that $c_1 < -2$ and $c_2=1$, so by Lemma
\ref{Zperp}(a) we see $c_3 <1$, i.e. $\lambda < \frac{3}{2}$.  Now the
second vector cannot occur because the orthogonality equation and
$\lambda \leq 2$ imply that $d_3 =1$, contradicting the presence of
$w$. Since the three vectors are collinear, if two satisfy the
orthogonality equation then all do.  So there is at most one (1A)
vertex and so $r \leq 4$ by Theorem \ref{estr}.  This can be improved
to $r \leq 3$ as follows. If the third vector $(-1, -1, 1,\cdots)$
occurs then the orthogonality relation, the bound on $\lambda$, and
nullity imply that $d_1=5, d_3=2$ and $\lambda = \frac{5}{7} \left(2 +
  \frac{1}{d_2} \right)$.  Now the nullity equation may be written as
a quadratic in $2 + \frac{1}{d_2}$ with no rational root.  If the
first vector $(-2,0,1,\cdots)$ occurs then orthogonality implies
$\lambda = \frac{d_1(d_3+2)}{4d_3 + 2d_1}$, and the bound on $\lambda$
gives $\frac{6}{d_1} + \frac{1}{d_3} >1$. Nullity implies $d_1 \geq 5$
and $d_1 > d_3^2$. We can deduce $d_3 =2$ and $\lambda =
\frac{2d_1}{d_1 + 4}$, and one can check that nullity fails.

For case (5), again a straight-forward analysis of the orthogonality
condition with the help of the nullity of $\bar{c}$ gives the
following $u \in {\mathcal W}$ as possibilities such that $J(\bar{c},
\bar{u})= 0$:
\begin{enumerate}
\item[(a)] $(-2^3, 1^i), \ i \geq 4$ and $d_3 =2$,
\item[(b)] $(1, -2, 0, \cdots)$,
\item[(c)]$(-2^2, 1^i), \  i \geq 4$ and $d_2 = 3$,
\item[(d)] $(0, -2, 1, 0, \cdots)$.
\end{enumerate}
Note $(1,0,-2,\cdots)$ cannot be in $\mathcal W$ as then $vw$ is not an edge.

It follows from Cor \ref{orthogonal} that among (a) only one
vector can occur and among (b), (c), (d) also only one vector can
occur.  (The orthogonality conditions of (b) and (d) are incompatible
with $1 < \lambda \leq 2$.) So $\bar{c}^{\perp} \cap \frac{1}{2}(d+
{\mathcal W})$ contains at most two elements. If it has two elements,
one must then come from (a) and the other from (b)-(d). Together they
give an edge of $\bar{c}^{\perp} \cap \ {\rm conv}(\frac{1}{2}(d+{\mathcal W}))$
with no interior points in $\frac{1}{2}(d+{\mathcal W})$.  Using
Cor \ref{orthogonal} and the null condition, we find that all these
two-element cases cannot occur. Hence $r \leq 4$.

If $r=4$ then the possible $u$ with $J(\bar{c}, \bar{u})=0$ are
given by (a)-(d) with $i=4$. Now, we can show using
techniques similar to those of Theorem \ref{cedge} that $c^* =
(1-2\lambda, \; 2\lambda -1, -1,0)$ also gives a null element of $\mathcal
C$. The possible vectors orthogonal to this element come from
(a) and the vectors (b$^*$), (c$^*$), (d$^*$) obtained from (b), (c), (d) by
swapping places 1 and 2.

If (a) does not give an element in $\bar{c}^{\perp} \cap \frac{1}{2}(d+ {\mathcal W})$,
it is straightforward to show, using the orthogonality and nullity
conditions for $\bar{c}$ and $\overline{c^*}$ together, that the (1A) vertices
for $\bar{c}$ and $\overline{c^*}$ are given by (b) and (b$^*$) respectively.
Also we must have $c = (\frac{4}{3}, -\frac{4}{3}, -1,0)$ and
$d= (4,4,9,d_4)$. We need a (1B) vertex outside $x_4=0$.
{}From Remark \ref{xc}, the only possible (1B) vertices
for $\bar{c}$ correspond to $(1,0,0,-2)$ and $(1,-1,0,-1)$.
In particular there can be no vertices, and hence no elements
of $\mathcal W$,  with $x_4 >0$. Hence $(-1,-1,0,1)$ and therefore
$(1,-1,0,-1)$ are not in $\mathcal W$. So the (1B) vertices for
$\bar{c}$ and $\overline{c^*}$ are given by $(1,0,0,-2)$ and $(0,1,0,-2)$
respectively. Now the line joining the corresponding null vectors
$a,a^*$ misses conv$(\mathcal W)$, a contradiction.

The remaining case is when (a) gives the element in
$\bar{c}^{\perp} \cap \frac{1}{2}(d+{\mathcal W})$
and in ${\overline{c^*}}^{\perp} \cap \frac{1}{2}(d+{\mathcal W})$. Now for
$vw$ to be an edge we need $(-1^4)$ absent, so by Remark \ref{facts}(b)
the only possible members of $\mathcal W$ lying outside $\{x_4 =0\}$ are
the three type IIIs with $x_4 =1$. In particular, all vectors in
${\rm conv}(\mathcal W)$ have $x_4 \geq 0$. As $(-1,-1,1,0) \in \mathcal W$, there must be
(1B) vertices lying in $\{x_4 =0 \}$ for both $\bar{c}$ and $\overline{c^*}$.
We then find that the only possibilities for such a (1B) vertex are
given by (b), (b$^*$) respectively. It follows that $d_1=d_2$,
but now nullity is violated.

We conclude that there are no (1A) vertices, so $r \leq 3$.

\begin{theorem} \label{no2}
Let $\bar{c}$ be a null vertex in $\mathcal C$ such that
$\Delta^{\bar{c}}$ contains a type $($2$)$ vertex  corresponding
to an edge $vw$ of ${\rm conv}({\mathcal W})$. Suppose we are
not in the situation of Theorem \ref{Fano}.

$($i$)$ If there are no points of $\mathcal W$ in the interior of $vw$, then
either we are in the situation of Example 8.2 of \cite{DW4}
or $K$ is not connected and we are in one of cases in the table of
Remark \ref{disconn} or in the situation of the third paragraph
of Example 8.3 of \cite{DW4}.

$($ii$)$ If there are interior points of $vw$ in $\mathcal W$
 then $r \leq 3$.   \; $\qed$
\end{theorem}

For further remarks about the $r=2$ case see the concluding remarks
at the end of \S 10.

\section{\bf Completing the classification}

Throughout this section we will assume that $K$ is connected (and we
are not in the situation of Theorem \ref{Fano}).
Theorems \ref{adj1B} and \ref{no2} then tell us that if $r \geq 4$ there
are no type (2) vertices and no adjacent (1B) vertices in
$\Delta^{\bar{c}}$, for {\em any} null vector $\bar{c} \in {\mathcal C}$.
Since all (1B) vertices lie in the half-space $\{J(\bar{c}, \cdot) > 0\}$
bounded by the hyperplane $\bar{c}^\perp$ containing the (1A) vertices,
we must therefore have in {\em each} $\Delta^{\bar{c}}$ exactly one (1B)
vertex, with the remaining vertices all of type (1A). So if $r \geq 4$
the only remaining task is to analyse such a situation.

As $\dim \Delta^{\bar{c}} = r-2$, we see $\dim (\Delta^{\bar{c}} \cap \bar{c}^{\perp}) =
r-3$. In particular, there must be at least $r-2$ elements of $\mathcal W$
giving elements of $\bar{c}^{\perp} \cap \frac{1}{2}(d + {\mathcal W})$.

Theorem \ref{class-1a} lists the possible configurations of such elements
when $r \geq 3$. The above discussion, together with Remark \ref{suppl-1a},
shows that in cases (1), (2), (3) we can take $m=r-1$,  in  cases (5)(ii),
(5)(iii), (6)(i), (6)(ii) we can take $m=r$, and in (6)(iii) we have
$r=5$.

Finally, since the vectors in (4) are collinear, so that
$\dim (\Delta \cap \bar{c}^{\perp}) = 1$, we have $r =3$ or $4$.
If $r=3$, it follows that $\bar{c}$ and the edge in
${\rm conv}(\frac{1}{2}(d + {\mathcal W}))$ determined
by the vectors in (4) are collinear. This contradicts the
orthogonality condition for the configuration
of vectors in (4) and we conclude that $r=4$.

For each configuration we can consider the possible vectors
$u \in {\mathcal W}$ giving the (1B) vertex. Besides the nullity condition
on $\bar{c}$ and the condition that $\bar{c}$ should be orthogonal to the
(1A) vectors, we have a further relation coming from the null condition
in Remark \ref{xc}. In most cases, routine (but occasionally tedious)
computations show that these relations have no solution.

As a result, we obtain the following possibilities for $u$ (up
to obvious permutations):

\[
\begin{array}{|l|c|c|} \hline
 {\rm case}   &   {\rm (1A) \; vectors} &    {\rm possible \; (1B) \; vector} \\ \hline
(1)   & (-2^1, 1^i),  \; 2 \leq i \leq r-1  &   (-1^r), (1^2,-2^r), (-1^1,-1^2,1^r)\\
   &                    &  (-1^1,1^2,-1^3), (-1^1,1^2,-1^r) \\

(2)   & (1^1, -2^i), \; 2 \leq i \leq r-1   &   (-1^2), (-1^2,-1^3, 1^r) \\

(3)   & (1^1, -2^2); (1^1,-2^3); {\rm possibly} (1^1,-1^2,-1^3);  &  \\
   & (-1^2,-1^3,1^i),\; 4 \leq i \leq r-1     &  (1^4, -2^r), (1^1, -2^r), (-1^r) \\

(4) & (1^1,-2^2), (1^1,-2^3), (1^1,-1^2,-1^3) &  (-1^4), (1^1,-2^4),
(1^1,-1^2,-1^4),  \\

(5)(i)   & (-2^1,1^2),  (-1^1,1^2,-1^3) &   (-1^1), \; (-1^4)      \\

(5)(ii)  & (-2^1,1^2); (-1^1,1^3,-1^i), \;  4\leq i \leq r     &  (1^3,-2^4),
(-1^1,1^2,-1^4),(1^3,-1^4,-1^5) \\
(5)(iii) & (-2^1,1^2); (-1^1,-1^3,1^i), \; 4 \leq i \leq r     &
 (-2^2,1^4),(-1^1,-1^2, 1^4), (-1^3,1^4, -1^5)\\
(6)(i)   & (-1^1,-1^2,1^i), \; 3 \leq i \leq r   &  (-1^1, -1^3, 1^4)
\\
(6)(ii)   & (1^1, -1^2, -1^i), \; 3 \leq i \leq r  & (1^1, -1^3, -1^4), (-1^2, 1^3, -1^4), (1^1, -2^3) \\
(6)(iii) & (1^1, -1^2, -1^i), i=3,4; \ (1^1, -1^3, -1^4) &   (-1^5) \\ \hline
\end{array}
\]
\centerline{Table 7: Unique $1B$ cases}

\begin{rmk} The possibilities for the (1B) vertex in cases (5)(ii)
 and (5)(iii) only apply to the $r \geq 5$ situation. When $r=4$, the
 two cases become the same if we switch the third and fourth summands
 and the possibilities are discussed in Lemma \ref{no5dim4} below.
\end{rmk}

\begin{rmk}
Note that $u$ such as $(-1^2), (-1^3)$ in (4), $(-1^4)$ in (5)(ii)
or $(-1^i)$ with $i >2$ in (6)(ii) cannot arise because they will
not be vertices, due to the presence of the type II vectors in $\mathcal W$.
\end{rmk}

\begin{rmk}
The dimensions must satisfy certain constraints in each case. Some such
constraints were stated in Theorem \ref{class-1a} and Remark \ref{suppl-1a}.
We also have constraints coming from the nullity conditions for $\bar{c}$
and $\bar{a}$. These typically involve the requirement that some expression
in the $d_i$ is a perfect square. The following is a summary of general
constraints in each case:

\noindent{case (1):  $d_1=4$;} case (2): $d_1=1$; case (3): $d_2 =d_3 =2$;
case (4): $d_2 + d_3 \leq 4d_1/(d_1 -1)$ and $d_2, d_3, \geq 2$;
case (5)(i): $(d_1, d_2) = (4, 2), (3, 3)$;
case (5)(ii,iii): $d_1=2$ and if $r \geq 5$ also $d_3 =2$; case (6)(i, ii):
$d_1 = d_2=2$; case (6)(iii): $d_1=d_2=d_3=d_4=2$, $d_5 =25$.
\end{rmk}

Our strategy now is reminiscent of that in \S 8. We have a (1B) vertex
corresponding to $\bar{u}$ and a second null vector $\bar{a}$ satisfying
$a = 2u -c$. Now we may apply our arguments to $\bar{a}$, and conclude that
the vectors in $\bar{a}^\perp \cap \frac{1}{2}(d+ {\mathcal W})$ are also
of the form given in the above table, up to permutation.  The resulting
constraints will allow us to finish our classification.

In some cases we can actually show that $\bar{a}^{\perp} \cap
\frac{1}{2}(d + {\mathcal W})$ is empty and we have a contradiction.
A simple example when this happens is case 6(iii), where we now have
$c=(\frac{6}{5}, -\frac{2}{5}, -\frac{2}{5}, -\frac{2}{5}, -1)$ and
$(a_i/d_i)=(-\frac{3}{5}, \frac{1}{5}, \frac{1}{5}, \frac{1}{5}, -\frac{1}{25})$.
Other cases are treated in Lemma \ref{noaperp} below.

Next we shall show that cases (6)(i)(ii) cannot arise (cf Lemma \ref{no6}),
so in all remaining cases there must be at least one type III vector $w$
with $\bar{w}$ in $\bar{a}^{\perp}$. We now use our explicit formulae
for $c$ and $a$ to derive inequalities on the entries of $a$ and find when
there can be such a type III vector orthogonal to $\bar{a}$. For each such
instance we then check whether $\bar{a}^{\perp} \cap \frac{1}{2}(d + {\mathcal W})$
forms a configuration equivalent to one of those in Table 4.
This turns out to be possible in only two situations (cf Lemmas \ref{case4},
\ref{case1} and Lemma \ref{case3}). These have such
distinctive features that $\mathcal W$ can be completely determined
and judicious applications of Prop \ref{edgeCC} lead to contradictions.
This yields our main classification theorem.

\begin{lemma} \label{noaperp}
The following cases cannot arise:

 case $($4$)$ with $u= (1,-1,0,-1)$,

 case $($6$)$$($ii$)$ with $u = (1,0,-2,\cdots)$,

 case $($4$)$  with $u=(0,0,0,-1)$ except for the case
$c = (\frac{4}{3}, -\frac{2}{3}, -\frac{2}{3}, -1)$
with $d=(4,2,2,9)$.
\end{lemma}

\begin{proof}
$(\alpha)$
 For case (4) with $u=(1,-1,0,-1)$ we find that the nullity and orthogonality
conditions and Remark \ref{xc} leave us with the following possibilities:

\[
\begin{array}{|c | c |c |} \hline
d  & c & (a_i/d_i) \\ \hline
(2,5,3,20) & (1, -\frac{5}{4}, -\frac{3}{4}, 0) & (\frac{1}{2},
-\frac{3}{20}, \frac{1}{4}, -\frac{1}{10}) \\
(2,6,2,12) & (1, -\frac{3}{2},-\frac{1}{2},0) &(\frac{1}{2},
-\frac{1}{12}, \frac{1}{4}, -\frac{1}{6})  \\
(3,4,2,12) & (1, -\frac{4}{3}, -\frac{2}{3}, 0) &
(\frac{1}{3},-\frac{1}{6},\frac{1}{3},-\frac{1}{6}) \\
(5,3,2,15) & (1, -\frac{6}{5}, -\frac{4}{5}, 0)&
(\frac{1}{5}, -\frac{4}{15}, \frac{2}{5}, -\frac{2}{15})\\
\hline
\end{array}
\]
It is easy to see that we can never have $\sum_i \frac{w_i a_i}{d_i} =1$
for $w \in \mathcal W$, so $\bar{a}^{\perp} \cap \frac{1}{2}(d + \mathcal W)$
is empty and we have a contradiction.

$(\beta)$ Similarly,  nullity, orthogonality and Remark \ref{xc} give:
\[
\begin{array}{|c| c |c|} \hline
d  & c & (a_i/d_i) \\ \hline
(2,2,225,d_4, \ldots, d_r) & (\frac{232}{225}, -\frac{29}{30}, -\frac{1}{4},
-\frac{d_4}{900}, \ldots, -\frac{d_r}{900}) &
(\frac{109}{225}, \frac{29}{60}, -\frac{1}{60}, \frac{1}{900},
\ldots, \frac{1}{900})   \\
(2,2,98,d_4, \ldots,d_r) & (\frac{52}{49},-\frac{13}{14},-\frac{1}{2},
-\frac{d_4}{196}, \ldots -\frac{d_r}{196})
&(\frac{23}{49},\frac{13}{28}, -\frac{1}{28}, \frac{1}{196},
\ldots, \frac{1}{196})\\
(2,2,36,d_4, \ldots d_r)) & (\frac{10}{9}, -\frac{5}{6},  -1,
-\frac{d_4}{36}, \ldots, -\frac{d_r}{36}) &
(\frac{4}{9},\frac{5}{12},-\frac{1}{12}, \frac{1}{36}, \ldots \frac{1}{36})
\\
\hline
\end{array}
\]
Moreover $n$ equals 962, 226, 50 respectively.  It is easy to see
that we can never have $\sum_i \frac{w_i a_i}{d_i} =1$ for $w \in
\mathcal W$, so we have a contradiction.

$(\gamma)$ For case (4) with $u=(0,0,0,-1)$ we find that the nullity
and orthogonality conditions and Remark \ref{xc} leave us with the
following possibilities, up to swapping places 2 and 3:
\[
\begin{array}{|c | c |c |} \hline
d  & c & (a_i/d_i) \\ \hline
(2,2,4,25) & (\frac{6}{5}, -\frac{2}{5}, -\frac{4}{5}, -1) & (-\frac{3}{5},
\frac{1}{5}, \frac{1}{5}, -\frac{1}{25}) \\
(2,3,3,25) & (\frac{6}{5}, -\frac{3}{5},-\frac{3}{5},-1) &(-\frac{3}{5},
\frac{1}{5}, \frac{1}{5}, -\frac{1}{25})\\
(3,2,3,121) & (\frac{15}{11}, -\frac{6}{11}, -\frac{9}{11}, -1) &
(-\frac{5}{11},\frac{3}{11},\frac{3}{11},-\frac{1}{121}) \\
(2m-2,2,2,m^2) & ( \frac{2(m-1)}{m}, \frac{1-m}{m}, \frac{1-m}{m},-1)&
(-\frac{1}{m}, \frac{m-1}{2m}, \frac{m-1}{2m}, -\frac{1}{m^2})\\
\hline
\end{array}
\]
It is now straightforward to see that we cannot have $w \in \mathcal W$
with $\sum_{i} \frac{w_i a_i}{d_i} =1$, except in two cases (both associated
to the last entry of the table). One is the case stated in the
Lemma. The other occurs if $m=2$, so $a = (-1, \frac{1}{2}, \frac{1}{2},
-1)$ which is orthogonal to $(-1,0,1,-1),(-1,1,0,-1)$. But as $(0,0,0,-1)$
is a vertex, neither of these vectors can be in $\mathcal W$. So
$\bar{a}^{\perp} \cap \frac{1}{2}(d + {\mathcal W})$ is still empty,
giving the desired contradiction.
\end{proof}

\medskip
As discussed above, we now turn to showing that case (6) cannot occur.
The following remark is useful in finding when type III vectors
can give elements of $\bar{a}^{\perp} \cap \frac{1}{2}(d+ {\mathcal W})$.

\begin{lemma} \label{ineg}
If $w=(-2^i, 1^j)$ and $\bar{w} \in \bar{a}^{\perp}$, then $\frac{a_i}{d_i} <0$
$($assuming we are not in the situation of Theorem \ref{Fano}$)$.
\end{lemma}
\begin{proof}
We need
\begin{equation} \label{aj}
\frac{a_j}{d_j}= 1 + \frac{2 a_i}{d_i},
\end{equation}
so if $\frac{a_i}{d_i} \geq 0$ then $\frac{a_j}{d_j} \geq 1$. Hence
$\frac{a_j^2}{d_j} = d_j \left( \frac{a_j}{d_j} \right)^2 \geq d_j \geq 1$.
As $\bar{a}$ is null this means $a = (-1^j)$ and we are in the situation of
 Theorem \ref{Fano}.
\end{proof}

\begin{lemma} \label{no6}
Configurations of type $($6$)$ cannot arise.
\end{lemma}

\begin{proof}
Recall that we have dealt with $($6$)$$($iii$)$ and
we have $d_1=d_2=2$.

\medskip

\noindent{Case (6)(i):} We have $u = (-1,0,-1,1, \cdots)$, and from
the nullity and orthogonality conditions we deduce that
\[
(c_i /d_i) = \left(-\frac{m+2}{2(m+1)}, \; -\frac{m-2}{2(m-1)}, \; \frac{1}{m^2-1},
 \cdots, \; \frac{1}{m^2-1} \right)
\]
where $\frac{1}{d_3} + \frac{1}{d_4} = \frac{1}{2(m+1)}$ and $n-1 =m^2$
for some positive integer $m$.
We have, therefore,
\[
(a_i/d_i) = \left(\frac{-m}{2(m+1)}, \; \frac{m-2}{2(m-1)}, \;
  -\frac{2}{d_3} - \frac{1}{m^2-1}, \; \frac{2}{d_4} - \frac{1}{m^2-1},
  \; \frac{-1}{m^2-1}, \cdots, \frac{-1}{m^2-1} \right).
\]
Let us estimate the size of the entries in $(a_i/d_i)$. First observe that,
as $d_3, d_4 > 2(m+1)$ from above, we have
\[
4(m+1) < d_3 + d_4 \leq n-d_1 -d_2 = n-4 = m^2-3
\]
so we deduce $m \geq 6$.  Hence we have $\frac{3}{7} \leq
|\frac{a_1}{d_1}| < \frac{1}{2}$, $\frac{2}{5} \leq
|\frac{a_2}{d_2}| < \frac{1}{2}, |\frac{a_3}{d_3}| \leq
\frac{6}{35}, |\frac{a_i}{d_i}| \leq \frac{1}{35}$ for $i \geq 5$.
Also note that $\frac{a_4}{d_4} < \frac{2}{d_4} < \frac{1}{m+1} \leq \frac{1}{7}$.
Finally, $\frac{a_4}{d_4} >0$, else we would have $2n-4 = 2m^2-2 \leq d_4$,
which is impossible.

Consider now a type III vector $w = (-2^i, 1^j)$
with $\bar{w} \in \bar{a}^{\perp}$.
By Remark \ref{ineg}, we need $i = 1, 3$ or $\geq 5$.
If $i=1$ then, by Eq.(\ref{aj}),we have
 $\frac{a_j}{d_j} = \frac{1}{m+1} \in (0, \frac{1}{7}]$.
So we must have $j=4$ and $\frac{1}{m+1} < \frac{2}{d_4}$, contradicting
our above remarks. If $i=3$ then $\frac{a_j}{d_j} \geq \frac{23}{35}$,
which is impossible. Similarly, if $i \geq 5$ then $\frac{a_j}{d_j} \geq
\frac{33}{35}$, which is impossible.

Hence there are no such type III vectors, so we are in case (6)
with respect to $\bar{a}$. We cannot be in 6(ii) as then the null
vector has exactly one positive entry (see below), but $a$ has two
positive entries. For 6(i), the null vector has exactly two negative
entries. Now $a$ has $r-2$ negative entries, so $r=4$. But for 6(i)
the negative entries have modulus $<2$, while
$a_3 = -2 - \frac{d_3}{m^2-1}$, a contradiction.

\medskip
\noindent{Case 6(ii):}  Here there are two possibilities.

\smallskip
\noindent{{\it Subcase} ($\alpha$)}: $u = ( 0,-1,1,-1,\cdots)$. Then, as above, the null
condition for $\bar{c}$ gives
\[
(c_i/d_i) = \left(\frac{(m-1)(m+2)}{2(m+1)^2}, \; -\frac{m+2}{2(m+1)},
\; \frac{-1}{(m+1)^2},\cdots, \; \frac{-1}{(m+1)^2} \right)
\]
where
 $\frac{1}{d_3} + \frac{1}{d_4} = \frac{1}{2(m+1)}$ and $n-1 =m^2$.
So the vector $(a_i/d_i)$ is given by
\[
\left( -\frac{(m-1)(m+2)}{2(m+1)^2}, \; \frac{-m}{2(m+1)},
\; \frac{2}{d_3} + \frac{1}{(m+1)^2}, \; -\frac{2}{d_4} + \frac{1}{(m+1)^2},
\; \frac{1}{(m+1)^2}, \cdots, \; \frac{1}{(m+1)^2} \right).
\]
As before, we have $m \geq 6$. We deduce $|\frac{a_1}{d_1}| \leq
\frac{1}{2}, \frac{3}{7} \leq |\frac{a_2}{d_2}| < \frac{1}{2}$,
$|\frac{a_3}{d_3}| \leq \frac{8}{49}, |\frac{a_4}{d_4}| \leq \frac{1}{7}$,
$|\frac{a_i}{d_i}| \leq \frac{1}{49}$ for $i \geq 5$.
Also $\frac{a_4}{d_4} < 0$, else $d_4 \geq 2(m+1)^2 > 2n$, which is impossible.

We look for vectors $w = (-2^i, 1^j)$ with $\bar{w} \in \bar{a}^{\perp}$.
By Lemma \ref{ineg} we have $i=1,2$ or $4$. If $i=1$ then $\frac{a_j}{d_j}=
\frac{m+3}{(m+1)^2}$, so $j=3$, but now Eq.(\ref{aj}) contradicts
$\frac{2}{d_3} < \frac{1}{m+1}$. A similar argument works if
$i=2$, while if $i=4$, Eq.(\ref{aj}) implies $\frac{a_j}{d_j} > \frac{5}{7}$,
a contradiction.

So 6(ii) must hold for $\bar{a}$, as we have already ruled out 6(i).
But now we need $a$ to have exactly one positive entry, which has modulus
$< 2$. So $r=4$ and this positive entry is $a_3$, but we have
$a_3 > 2$, a contradiction.

\smallskip
\noindent{{\it Subcase} ($\beta$)}: $u = (1,0,-1,-1, \cdots)$. We similarly have
\[
(c_i/d_i)= \left(\frac{(m-2)(m+1)}{2(m-1)^2}, \; -\frac{m-2}{2(m-1)},
\; \frac{-1}{(m-1)^2}, \cdots, \; \frac{-1}{(m-1)^2} \right).
\]
\[
(a_i/d_i) = \left(\frac{m^2 -3m +4}{2(m-1)^2}, \; \frac{m-2}{2(m-1)},
\; \frac{1}{(m-1)^2} -\frac{2}{d_3}, \; \frac{1}{(m-1)^2} -\frac{2}{d_4},
\; \frac{1}{(m-1)^2},\cdots, \; \frac{1}{(m-1)^2} \right)
\]
where $n-1 = m^2$ and $\frac{1}{d_3} + \frac{1}{d_4} = \frac{m+1}{2(m-1)^2}$.
The last two equations easily imply that $m \geq 6$, and
so $|\frac{a_i}{d_i}| < \frac{1}{2}$ for all $i$.

A type III $(-2^i, 1^j)$ giving an element of $\bar{a}^{\perp}$
must have $i=3$ or $4$, by Lemma \ref{ineg}. In both cases we find
from Eq.(\ref{aj}) that $\frac{a_j}{d_j} > \frac{1}{2}$, which is impossible.
So 6(ii) holds for $a$, which is impossible as $a$ has at least two positive entries.
\end{proof}


\begin{lemma} \label{case4}
The only possible example in case $($4$)$ is when
$c = \left(\frac{4}{3}, -\frac{2}{3}, -\frac{2}{3}, -1 \right)$,
$u = (0,0,0,-1),$ and $d =(4,2,2,9)$. $\Delta^{\bar{a}}$ is then in
case $($1$)$ with $a=(-\frac{4}{3}, \frac{2}{3}, \frac{2}{3}, -1)$
and $\bar{a}^{\perp} \cap \frac{1}{2}(d+ {\mathcal W})$ consists of
$(-2, 1, 0,0), (-2, 0, 1,0)$.
\end{lemma}

\begin{proof}
By Lemma \ref{noaperp} we just have to
eliminate the possibility  $u = (1,0,0,-2)$. Now
\[
(a_i/d_i) = \left( \frac{8-2d_1 + d_4}{d_1(4 + 2d_1 +d_4)},
\; \frac{(d_1-1)(2d_1 + d_4)}{2d_1(4 + 2d_1  +d_4)},
\; \frac{(d_1-1)(2d_1 + d_4)}{2d_1(4 + 2d_1  +d_4)}, \;
-\frac{2}{d_4} -\frac{2(d_1-1)}{d_1 (2d_1 + d_4 +4)} \right).
\]
The null condition is
\[
-(d_1 -1)\; d_4^3 - 4(d_1^2 -d_1 -1)\; d_4^2 +4d_1 (3d_1^2 + d_1 +8)d_4
+ 16 d_1^2 (d_1 +2)^2 =0.
\]

For $w=(-2^i, 1^j)$ with $ \bar{w} \in \bar{a}^{\perp}$ we need,
by Lemma \ref{ineg}, $i=1$ or $4$. If $i=1$, then for $j=2$ or $3$,
Eq.(\ref{aj}) can be rewritten as $2d_1^2 + d_1 d_4 + 2d_1 = -5d_4 -32$,
which is absurd. For $j=4$ it can be rewritten as $-1 -\frac{2}{d_4} = \frac{14 +
2d_4 -2d_1}{d_1(2d_1 + d_4 +4)}$. So the right hand side is $< -1$, which
on clearing denominators is easily seen to be false.

If $i=4$, then for $j=1$ Eq.(\ref{aj}) becomes  $\frac{4}{d_4} +
\frac{1}{d_1}=1$. So $(d_1, d_4)=(2,8), (3,6)$
or $(5,5)$, all of which violate the null condition. For $j=2,3$
we obtain from Eq.(\ref{aj}) the equation
$4d_4(d_1 -1) = (2d_1+d_4 +4)(d_1d_4 +d_4 -8d_1)$, which can only
have solutions if $d_4 \leq 9$. On the other hand, the null condition
has no integer solutions if $d_4 \leq 9$.

So no such type III exists, contradicting Lemma \ref{no6}.
\end{proof}

\begin{lemma} \label{no5(i)} Configurations of type
$($5$)$$($i$)$ cannot occur.
\end{lemma}

\begin{proof} It is useful to note that the null condition for
$\bar{c}$ implies that $d_3 \leq 4$ when $(d_1, d_2) = (4,2)$
and $d_3 \leq 3$ when $(d_1, d_2)=(3,3)$. One further finds the following
possibilities:
$$ \begin{array}{|c|c|c|c|} \hline
   u &  d  &  c   &  (a_i/d_i)  \\ \hline
   (-1,0,0,0) & (3,3,1,2)  & (-1, 1, -\frac{1}{3}, -\frac{2}{3}) &
           (-\frac{1}{3}, -\frac{1}{3}, \frac{1}{3}, \frac{1}{3}) \\
        & (3, 3, 2, 1) & (-1, 1, -\frac{2}{3}, -\frac{1}{3}) &
           (-\frac{1}{3}, -\frac{1}{3}, \frac{1}{3}, \frac{1}{3}) \\
      & (4, 2, 1, 3) & (-1, 1, -\frac{1}{4}, -\frac{3}{4}) &
           (-\frac{1}{4}, -\frac{1}{2}, \frac{1}{4}, \frac{1}{4}) \\
      &(4, 2, 2, 2)  & (-1, 1, -\frac{1}{2}, -\frac{1}{2})   &
            (-\frac{1}{4}, -\frac{1}{2}, \frac{1}{4}, \frac{1}{4}) \\
      & (4, 2, 3, 1) &  (-1, 1, -\frac{3}{4}, -\frac{1}{4}), &
            (-\frac{1}{4}, -\frac{1}{2}, \frac{1}{4}, \frac{1}{4}) \\ \hline
   (0, 0, 0, -1) & (3, 3, 2, 121) &  (-\frac{9}{11}, \frac{15}{11}, -\frac{6}{11}, -1)
          &  (\frac{3}{11}, -\frac{5}{11}, \frac{3}{11}, -\frac{1}{121})    \\
      &   (4, 2, 2, 25)  &  (-\frac{4}{5}, \frac{6}{5}, -\frac{2}{5}, -1)
          & (\frac{1}{5}, -\frac{3}{5}, \frac{1}{5}, -\frac{1}{25})    \\
   \hline
\end{array} $$

One easily checks that $\bar{a}^{\perp} \cap \frac{1}{2}(d+ {\mathcal W})$ is
empty in the last two cases, and consists only of type II vectors
in the third to fifth cases, giving a contradiction to Lemma \ref{no6}.
For the first two cases, note that
$\bar{a}^{\perp} \cap \frac{1}{2}(d+{\mathcal W})$ contains $\frac{1}{2}(d+(-1, -1, 1, 0))$
since by hypothesis for 5(i) $(-1,1,-1,0)$ is in $\mathcal W$.
Hence (1),(2) cannot hold with respect to $\bar{a}$. Also the vector $d$ of dimensions
rules out (3),(4) and (5), so we have a contradiction.
\end{proof}

\begin{lemma} \label{case3}
The only possible example for case $($3$)$ is when
$ c = \left(\frac{2}{3}, -\frac{2}{3}, -\frac{2}{3},\frac{2}{3},-1 \right),  \;
u = (0,0,0, 0,-1), \; d =(2,2,2,2,9)$. Then $a=(-\frac{2}{3}, \frac{2}{3}, \frac{2}{3},
-\frac{2}{3}, -1)$ and $\Delta^{\bar{a}}$ is again in case $($3$)$.
\end{lemma}

\begin{proof}
(A) Let $u = (-1^r)$. The null condition for $\bar{c}$ gives
$4d_r = (\delta + d_1 +2)^2$, where $\delta= d_4 + \cdots + d_{r-1}$.
In particular, $d_r$ is a square. Also, $(a_i / d_i) =
\left( \frac{-1}{\sqrt{d_r}}, \frac{1}{2} \left(1 - \frac{1}
{\sqrt{d_r}} \right),\frac{1}{2} \left(1 - \frac{1}{\sqrt{d_r}} \right),
\frac{-1}{\sqrt{d_r}}, \cdots, \frac{-1}{\sqrt{d_r}}, \frac{-1}{d_r} \right)$.

If $d_r=4$, then we find there are no type III vectors in
$\bar{a}^{\perp}$, a contradiction. So $d_r \geq 9$ and we have
$|\frac{a_i}{d_i}|  \leq \frac{1}{3}$
for $i=1,4,\cdots, r-1$, $\leq \frac{1}{9}$ for $i=r$ and $< \frac{1}{2}$
for $i=2,3$. Lemma \ref{ineg} shows $i \neq 2,3$.
 From Eq.(\ref{aj}) and the above estimates, we first get $i \neq r$,
and for the remaining values of $i$, we have $\frac{a_j}{d_j} >0$,
so that $j=2$ or $3$. Also, $d_r =9$ (and hence $d_1 + d_4 + \cdots + d_{r-1}=4$), and
$(a_i/d_i) =( -\frac{1}{3}, \frac{1}{3}, \frac{1}{3}, -\frac{1}{3},
\cdots, -\frac{1}{3}, -\frac{1}{9})$. Upon applying Theorem \ref{class-1a}
to $\Delta^{\bar{a}}$ together with Theorem \ref{estr} and the
above Lemmas, we deduce that we are in case (3)(ii) with $r=5$ and
$d_1=d_4 =2$, giving the example in the statement of the Lemma.

(B) Next let $u = (1,0, \cdots,-2)$.
Now $(a_i/d_i) = (\frac{2}{ d_1} - \alpha, \frac{1- \alpha}{2}, \frac{1 -
\alpha}{2}, -\alpha, \cdots, -\alpha, \frac{(n-2-d_r)\alpha -5}{d_r} )$
where, as a consequence of the null condition for $\bar{c}$, we have
\begin{equation} \label{alpham}
 \alpha = \frac{c_1}{d_1}= \frac{1}{n-2-m} \;\;\; : \;\;\;  m^2 = d_r (n-1).
\end{equation}

Next we get  the identity $(n-2)^2 -m^2 = (n-1)(d_1 + 1 + \delta) +1=
\frac{m^2(d_1 + 1 + \delta)}{d_r}  + 1$, where $\delta$ is as given in
(A) above. We deduce $m < \sqrt{\frac{d_r}{n-3}}(n-2)$, and hence
$\alpha < \frac{2}{d_1 + 1 + \delta} =\frac{2}{n-d_r-3}
\leq \min(\frac{1}{2},\frac{2}{d_1})$.
So $\frac{a_i}{d_i}$ is positive for $i \leq 3$ and negative for $3 <
i \leq r-1$.
Note also that $(n-2-d_r) \alpha <  \frac{2(n-2-d_r)}{n-3-d_r}
=2 \left(1 + \frac{1}{n-3-d_r} \right) \leq \frac{5}{2}$.
In particular, $\frac{a_r}{d_r} < -\frac{5}{2d_r} < 0$.

As usual, we look for $(-2^i, 1^j)$ giving an element of
$\bar{a}^{\perp}$.  By Lemma \ref{ineg}, $i\neq 1,2,3$.
If $4 \leq i \leq r-1$ then Eq.(\ref{aj}) says $\frac{a_j}{d_j} = 1 -2
\alpha >0$, so $j=1,2$ or $3$.  If $j=1$ we obtain $\alpha = 1- \frac{2}{d_1}$.
Comparing this with Eq.(\ref{alpham}) shows $d_1=3$ and $\alpha = \frac{1}{3}$,
but now $\bar{a}$ is not null. If $j=2$ or $3$, we obtain $\alpha = \frac{1}{3}$;
we deduce from Eq.(\ref{alpham}) that $d_1 \leq 5$, and again one can check that
all possibilities violate nullity.

So all type III $(-2^i, 1^j)$ have $i=r$. Therefore we must be in case
(1) or (5) with respect to $\bar{a}$, and $d_r =4$ or $2$ respectively.

If $j=1,2,3$ then $\frac{a_j}{d_j} >0$. Now Eq.(\ref{aj}) combined with
the estimate above for $\frac{a_r}{d_r}$ show that
$-\frac{5}{2} > a_r > -\frac{d_r}{2}$, so $d_r > 5$, a contradiction.

If $4 \leq j \leq r-1$, then in the case $d_r=4$, we find Eq.(\ref{aj})
gives $\alpha = \frac{3}{n-4}$. Combining with Eq.(\ref{alpham}) we get
$n=10$, which is incompatible with $d_r =4$ and $r \geq 5$. If $d_r=2$
we find similarly that $\alpha = \frac{4}{n-3}$ and $m$ satisfies
$3m^2 - 8m-4 =0$; but this has no integral roots.

So $u = (1,0, \cdots -2)$ cannot occur.

(C) For $u = (1^4, -2^r)$, we have
$(a_i/d_i) = ( -\alpha, \frac{1- \alpha}{2}, \frac{1 -
\alpha}{2}, \frac{2}{d_4} - \alpha, -\alpha, \cdots, -\alpha, \frac{(n-2-d_r)
\alpha -5)}{d_r} )$
and Eq.(\ref{alpham}) still holds. The arguments of case (B) carry
over to this case, on swapping indices $1,4$.
\end{proof}

\begin{lemma} \label{elim2}
Case $($2$)$  cannot occur.
\end{lemma}

\begin{proof}
(A) Consider $u=(-1^2)$. Now $(a_i/d_i)=( \frac{2}{d_2}-1,
\frac{-1}{d_2}, \frac{1}{d_2}, \cdots, \frac{1}{d_2},
\frac{1}{d_r} ( 2 - \frac{n+1-d_r}{d_2}))$.

Nullity implies $d_2 \geq 3$, so $-1 < \frac{a_1}{d_1} \leq -\frac{1}{3}$ and
$| \frac{a_i}{d_i} | \leq \frac{1}{3}$ for $2 \leq i \leq r-1$.
Also we have $d_r(n-1) =m^2$, and for this choice of $u$ we have $m=n+1-2d_2$,
so $\frac{a_r}{d_r} = \frac{d_r -m}{d_2 d_r}$ which is positive if $m <0$ and
negative if $m>0$. By Lemma \ref{ineg} and the fact that $d_1=1$,
we only have to consider $(-2^i, 1^j)$ with $i=2$ or $r$.

If $i=2$, then Eq.(\ref{aj}) says $\frac{a_j}{d_j} = 1-\frac{2}{d_2} >0$
so $j \geq 3$.  If $3 \leq j \leq r-1$ then Eq.(\ref{aj}) shows $d_2=3$,
so $m=n-5$ and $d_r = \frac{(n-5)^2}{n-1}= n-9 + \frac{16}{n-1}$.
As $n \geq 7$ we must then have $(d_r, n)=(9,17)$ or $(2,9)$. Imposing
the nullity condition on $\bar{a}$ shows there are only three
possibilities, corresponding to $d = (1,3,2,2,9), (1,3,4,9),(1,3,3,2)$.
In the first two there is only one type III in $\bar{a}^{\perp}$,
as $d_1=1$ and $d_2 \neq 4$, so we are in case (5) with respect to
$\bar{a}$, contradicting the fact that $d_2 \neq 2$. In the last case
we must be in case (4) with respect to $\bar{a}$, but now
$(0,-1,1,-1)$ is present, so $u$ is not a vertex.
If $j=r$ then Eq.(\ref{aj}) becomes $d_2 = \frac{3d_r -n-1}{d_r-2}$
which is less than 3, a contradiction. (We cannot have $d_r=2$ and
$3d_r=n+1$ as $n \geq 7$.)

Hence all such $(-2^i, 1^j)$ have $i=r$ and we are in case (1) or (5)
with respect to $\bar{a}$. For case (1) we need $r-2$ of the
$\frac{a_j}{d_j}$ ($j < r$) equal. This can only happen for our $a$
if $d_2=3$, which is ruled out as in the previous paragraph.
For case (5) we have $d_r=2$ and $\frac{a_j}{d_j} = 3 -\frac{n-1}{d_2}$.
The possibilities on the left-hand side are $\frac{2}{d_2}-1, -\frac{1}{d_2} ,\frac{1}{d_2}$
respectively. On using our relations for $m,n,d_r$ we find that
only the third possibility can occur, and $d_2 =3$. The argument in
the previous paragraph again eliminates this case.

(B) Consider $u=(0,-1,-1, \cdots,1)$. Now
$(a_i/d_i) = ( 2 \beta -1, \beta-\frac{2}{d_2}, \beta-\frac{2}{d_3}, \beta, \cdots, \beta,
\frac{4 - (n+1-d_r) \beta}{d_r} )$
where again from the null condition of $\bar{c}$ we have
\begin{equation} \label{betam}
\beta := \frac{1}{2}(1- c_1) = \frac{2}{n+1-m} = \frac{3 + d_r (\frac{1}{d_2} + \frac{1}{d_3})}
{n +d_r +1} \;\; : \;\; d_r(n-1)=m^2 \; (m>0).
\end{equation}
Now $0 < \beta < \frac{1}{2}$ (the case $\beta = \frac{1}{2}$ leads to
$r=4$, $d_2=d_3=2$ and $a=(0,-1,-1,1)$ which violates nullity).
So for $\bar{w} \in \bar{a}^{\perp}$ we just have to consider
$w=(-2^i, 1^j)$ with $i=2,r$ (as $d_1=1$ we can't have $i=1$;
also by symmetry the case $i=3$ is treated the same way as $i=2$).

If $i=2$ then Eq.(\ref{aj}) says $\frac{a_j}{d_j} = 2 \beta + 1 - \frac{4}{d_2}$.
If $4 \leq j \leq r-1$ we get $\beta = \frac{4}{d_2}-1$; the only possibility
consistent with our bounds on $\beta$ is $d_2 =3, \beta =\frac{1}{3}$ and it is
straightforward to check this is incompatible with the null condition for $\bar{a}$.

If $i=2$ and $j=1$ then Eq.(\ref{aj}) implies $d_2=2$ and again one checks
that nullity for $\bar{c}$ fails. If $j=3$ Eq.(\ref{aj}) says
$\beta = \frac{4}{d_2} - \frac{2}{d_3}-1$, so as $\beta >0$ either $d_2=2$
and $\beta = 1 -\frac{2}{d_3}$ or $d_2=3$ and $\beta = \frac{1}{3} - \frac{2}{d_3}$.
In the former case the bound $\beta < \frac{1}{2}$ shows $d_3 =3$ and
$\beta = \frac{1}{3}$, and now nullity for $\bar{a}$ fails. In the latter case
the bound $\beta >0$ shows $d_3 >6$. Substituting this into the quadratic
which must vanish for nullity of $\bar{c}$, we see
$\delta = d_4 + \cdots + d_{r-1}$ is $<4$.
Checking the resulting short list of cases yields no examples where
nullity holds. If $j=r$ then Eqs.(\ref{aj}) and Eq.(\ref{betam}) imply
$\frac{1}{d_r} + \frac{3}{d_2} = 1 + \frac{1}{d_3}$ and one check that
the possible $(d_2, d_r)$ yield no examples where nullity of $\bar{c}$ holds.

So all such type III have $i=r$, and  case (1) or (5) holds for $\bar{a}$.
For case (1), then, as in (A), $r-2$ of the $\frac{a_j}{d_j}$ ($j \leq r-1$)
must be equal. So either $r=5$ and $d_2 = d_3$ with $\beta
= 1 - \frac{2}{d_2}$ or $r=4$ and one of the preceding equalities holds.
If $\beta = 1-\frac{2}{d_2}$ holds, then the bounds on $\beta$ show
$d_2=3, \beta = \frac{1}{3}$ and as usual nullity for $\bar{a}$ fails.
If $d_2 = d_3$ holds, then using our formulae for $\beta$ and substituting
into the null condition for $\bar{a}$ gives a quadratic with no integer roots.

If case (5) holds, then $d_r =2$. Now Eq.(\ref{aj}) gives
$\frac{a_j}{d_j} = 5 -(n-1) \beta$. If $j=1$, $j=2$, or $4 \leq j \leq r-1$
we get $\beta = \frac{6}{n+1}, \frac{5 + (2/d_2)}{n}, \frac{5}{n}$
respectively. (As usual, the case $j=3$ is treated in just the same
way as $j=2$.) Now using the equations in Eq.(\ref{betam}) relating $n,m$ in each
case gives a quadratic with no real roots.
\end{proof}

\begin{lemma} \label{case1}
For case $($1$)$ the only possibility is when $ c =
(-\frac{4}{3}, \frac{2}{3},\frac{2}{3}, -1 )$ with
$u=(0,0,0,-1)$ and $d= (4,2,2,9).$ $\Delta^{\bar{a}}$
is then in case $($4$)$.
\end{lemma}

\begin{proof}
(A) Consider $u=(0, \cdots, 0,-1)$.  From the null condition
for $\bar{c}$ we see that  $d_r = k^2, n-1 = (k+1)^2$ for some positive
integer $k$ and $(a_i/d_i) =
(\frac{1}{2}( 1-\frac{1}{k}), -\frac{1}{k}, \cdots, -\frac{1}{k}, \frac{-1}{k^2})$.
Note that  since $d_1=4$,  $n >5$ and so $k \neq 1$.

We must consider solutions of Eq.(\ref{aj}). By Lemma \ref{ineg}, $i \neq 1$.
If $i=r$ we have $\frac{a_j}{d_j}= 1- \frac{2}{k^2}$. The resulting equation
has no solution in integer $k >1$ for any choice of $j$.
If $2 \leq i \leq r-1$, we need $\frac{a_j}{d_j} = 1 -\frac{2}{k}$.
We only obtain a solution $k > 1$ if $j=1$; in this case $k=3$, so $n=17$,
$d_r =9$ and we see $r=4$ with $\{d_2, d_3\}=\{ 2,2 \}$ or $\{ 3,1 \}$.
The former case is that in the statement of the Lemma. In the latter
case we can have just one type III and one type II in $\bar{a}^{\perp}$
(since $d_2$ or $d_3$ is 1, one potential type III is missing), so we must be in
case (5) with respect to $\bar{a}$; but no $d_i$ is 2, a contradiction.

(B) Consider $u = (0,1,0,\cdots, 0,-2)$. Now
$(a_i/d_i)=(\frac{1}{2}(1-\beta), \frac{2}{d_2}-\beta, -\beta, \cdots, -\beta,
\frac{(n-d_r-2) \beta -5}{d_r})$
where $\beta = \frac{d_r + 6d_2}{d_2(2n-d_r-4)}$. The nullity condition
for $\bar{c}$ implies $d_2 \geq 3$, $d_2 > \delta$ and $d_r > 2d_2 +4$,
where $\delta$ now denotes $d_3 + \cdots + d_{r-1}$.
We can then  deduce that $\beta < 1$, $0 < \frac{a_1}{d_1} < \frac{1}{2}$,
$0 < \frac{a_2}{d_2} < \frac{1}{3}$. In particular $\beta < \frac{2}{d_2}$.
By Lemma \ref{ineg}, we must consider elements of $\bar{a}^{\perp}$
coming from vectors $(-2^i, 1^j)$ with $i \geq 3$.

If $3 \leq i \leq r-1$, Eq.(\ref{aj}) says $\frac{a_j}{d_j}=
1-2 \beta$. As $\beta < 1$, this immediately
rules out $3 \leq j \leq r-1$. If $j=1$ we get $\beta = \frac{1}{3}$.
Combining this with our formula above for $\beta$ we get
$2d_2( d_2 + \delta -7) + (d_2-3)d_r=0$. The only possibilities
are $d_2=3, \delta=4$ which violates the null condition, or $d=(4,4,2,8)$
which violates the condition that $d_r(n-1)$ should be a square.
If $j=2$, we get $\beta = 1 - \frac{2}{d_2}$. Since we saw
above that $\beta < \frac{2}{d_2}$
we get $d_2 =3$ and $\beta = \frac{1}{3}$, which is ruled out as above.
If $j=r$, Eq.(\ref{aj}) implies $\beta = \frac{d_r +5}{2 + d_2 + \delta + 2d_r}$.
Comparing this with the formula for $\beta$ above leads to a contradiction.

The remaining possibility is for $i=r$. So we are in case (1) or (5)
 with respect to $\bar{a}$, and $d_r = 4$ or $2$ respectively.
But $d_r > 2d_2 + 4$, so this is impossible.

(C) Let $u =(-1,-1,0, \cdots, 0, 1)$. Now
$(a_i/d_i) = (-\frac{1}{2} \beta, -\frac{2}{d_2} - \beta, -\beta,\cdots,
-\beta, \frac{1 + (n-d_r-2) \beta}{d_r})$
where $\beta = \frac{d_r(d_2-4)}{2d_2 (2n + d_r -4)}$, so
 $0 < \beta < \frac{1}{6}$ (noting that the nullity condition
for $\bar{c}$ implies $d_2 \geq 5$).

We look for vectors $(-2^i, 1^j)$ giving elements of $\bar{a}^{\perp}$.
Now Lemma \ref{ineg} rules out $i=r$, while if $3 \leq i \leq r-1$ we need
$\frac{a_j}{d_j} = 1-2 \beta > \frac{2}{3}$. So $j=r$, and Eq.(\ref{aj})
yields $\beta = \frac{d_r-1}{n+d_r-2}$. Equating this to the
expression above for $\beta$ gives an equation which may be
 rearranged
so that it says a sum of positive terms  is zero, which is absurd.

If $i=1$ then Eq.(\ref{aj}) says $\frac{a_j}{d_j} =1 - \beta$. Clearly
this can only possibly hold if $j=r$. The equation then gives
$\beta = \frac{d_r-1}{n-2}$, and equating this with the earlier expression
for $\beta$  leads, as in the previous paragraph, to a contradiction.

So the only possibility is $i=2$, and we are therefore in case (1) or (5)
with respect to $\bar{a}$. But $d_2 \geq 5$ so this is impossible.

(D) Let $u =(-1,1,-1,0,\cdots)$.  Now
$(a_i/d_i)= (-\frac{1}{2} \beta, \frac{2}{d_2} - \beta,
-\frac{2}{d_3} - \beta, -\beta, \ldots, -\beta, \frac{(n-d_r-2)\beta -1}{d_r})$.
and $\beta = \frac{1}{2} - 2 (\frac{1}{d_2} + \frac{1}{d_3})$. An analysis
of the nullity condition for $\bar{c}$ shows that it can only be satisfied
if $\frac{1}{5} < \frac{1}{d_2} + \frac{1}{d_3}
< \frac{1}{4}$, so $d_2, d_3 \geq 5$ and $0 < \beta < \frac{1}{10}$.

Let us now consider solutions to Eq.(\ref{aj}).
If $i=r$, we have $\frac{a_j}{d_j} = 1 -\frac{2}{d_r} + \frac{2(n-d_r-2)\beta}{d_r}$.
If $j \neq 2$, this equation implies that the positive quantity
$(1 + \frac{2(n-d_r-2)}{d_r})\beta$
(or $(\frac{1}{2} + \frac{2(n-d_r-2)}{d_r})\beta$ if $j=1$)
equals a nonpositive quantity (recall $d_r > 1$ as $i=r$).
If $j=2$, we get that it equals $\frac{2}{d_2} + \frac{2}{d_r} -1$. But $d_2 \geq 5$
so $d_r =2$ or $3$, and in each case we find the nullity condition
for $\bar{c}$ is violated.

If $i=1$, Eq.(\ref{aj}) says $\frac{a_j}{d_j} = 1- \beta > \frac{9}{10}$,
so $j=2$ or $r$. But for $j=2$ we get $d_2 =2$, which is impossible as we
know $d_2 \geq 5$, so in fact $j=r$.

If $i=2$, Eq.(\ref{aj}) is $\frac{a_j}{d_j} = 1 + \frac{4}{d_2} - 2 \beta$.
We cannot then have $j=1,3$ or $4 \leq j \leq r-1$ as they lead to
$\beta > \frac{2}{3}, >1 , >1$ respectively. So we must have $j=r$.

If $i=3$, we see $\frac{a_j}{d_j}= 1 -\frac{4}{d_3} -2 \beta$. If
$j=1,2$ or $4 \leq j \leq r-1$, we see in all cases
(using our bounds on $d_2, d_3$) that $\beta > \frac{1}{5}$,
a contradiction. Hence again $j=r$.

If $4 \leq i \leq r-1$, then $\frac{a_j}{d_j} = 1-2 \beta > \frac{4}{5}$
so $j=2$ or $r$. If $j=2$ we obtain $\beta = 1-\frac{2}{d_2} \geq \frac{3}{5}$,
contradicting our earlier inequality for $\beta$; so again we have $j=r$.

We have shown that any $(-2^i, 1^j)$ giving an element of $\bar{a}^{\perp}$
has $j=r$, so we are in case (3), (4) or (5) with respect to $\bar{a}$.
It cannot be case (3) as we know from Lemma \ref{case3} that then
each $d_i$ is $2$ or $9$, and we have $d_1=4$. If we are in case (4),
then Lemma \ref{case4} tells us that $d=(4,2,2,9)$. Moreover,
as $(-2,0,0,1), (0,-2,0,1)$ are the elements of
$\bar{a}^{\perp} \cap \frac{1}{2}(d+ {\mathcal W})$,
we must have $\beta = \frac{4}{d_2}$; but now $\beta > \frac{1}{10}$, a
contradiction. If it is case (5), then we have $d_i=2$ for some $i$, which we
can take to be 4. Now $\bar{a}$ must be orthogonal to vectors
associated to $(-1^4, 1^5, -1^k)$ or $(-1^4, -1^5, 1^k)$, and either
case is incompatible with our expressions for $a_i/d_i$.

(E) Consider $u=(-1,1, 0, \cdots,0, -1)$. Now
$(a_i/d_i) = (-\frac{1}{2} \beta, \frac{2}{d_2}-\beta, -\beta, \cdots, -\beta,
\frac{((n-d_r-2)\beta -3)}{d_r} )$
and $\beta = \frac{8d_2 - d_r(d_2-4)}{2d_2(2n-d_r-4)}$. It is easy to check
that $\beta < \frac{2}{d_2}$. Also, the nullity condition for $\bar{c}$
implies $d_2 \geq 3$ and $(\frac{4}{d_2}-1)d_r + 8 >0$; hence $\beta > 0$.

The analysis is similar to that in (D). If $i=r$ then Eq.(\ref{aj})
implies that a positive quantity times $\beta$ equals a positive
linear combination of reciprocals of $d_i$, minus 1. This sum
of reciprocals is therefore $>1$, which gives us upper bounds on $d_r$.
The only case where Eq.(\ref{aj}) and the null condition can hold is
if $j=2$ and $d_2=7, d_r=4, \; d_3 + \cdots + d_{r-1}=11$

If $i=1$ then Eq.(\ref{aj}) says $\frac{a_j}{d_j} = 1 -\beta > 0$, so
$j=2$ or $r$.  But $j=2$ implies $d_2=2$, which from above cannot hold, so $j=r$.

Now Lemma \ref{ineg} rules out $i=2$. If $3 \leq i \leq r-1$, we have
$\frac{a_j}{d_j} = 1 - 2 \beta$. If $j=1$ then we get $\beta =\frac{2}{3}$,
which cannot hold.
If $j=2$ then $\beta = 1 - \frac{2}{d_2}$, and as
$\beta < \frac{2}{d_2}$ we deduce $d_2 =3$ and $\beta = \frac{1}{3}$,
which violates the null condition for $\bar{c}$.
If $3 \leq j \leq r-1$, then $\beta=1$, which is impossible. So we have $j=r$.

So in all cases we have $j=r$, except in the exceptional case discussed
above where we can have $i=r$ and $j=2$. But our list (1)-(6) of possible
configurations in $\bar{a}^{\perp}$ shows that if the $(i,j)=(r,2)$
case occurs then no other type III can be in $\bar{a}^{\perp}$. So
we are in case (5), which is impossible as $d_r =4 \neq 2$ for this example.
Hence the exceptional case cannot arise.

We see therefore that $j=r$ in all cases.  So, as in (D), we must be
in case (3), (4) or (5) with respect to $\bar{a}$. As before, the fact
that $d_1=4$ rules out case (3). For case (4) we need the $d_k$ to be
$4,2,2,9$ and $d_r$ to be $4$ (as the $(-2^i,1^j)$ in $\bar{a}^{\perp}$
have $j=r$), but this contradicts $d_1=4$.

So we are in case (5). Now  the orthogonality condition for the
family of type II vectors leads to $\beta > \frac{2}{3}$, which is
impossible.
\end{proof}

\begin{lemma} \label{no5a}
Case $($5$)$$($ii$)$ cannot occur if $r \geq 5$.
\end{lemma}

\begin{proof}
(A) Consider $u = (0,0,1,-2, 0,\cdots)$. We have

$(a_i/d_i) = (1 - \beta, \; 1-2 \beta,
\; \frac{n-2d_2 -3}{n-d_2 -2} - (\frac{n-6-3d_2}{n-d_2-2}) \beta,
\; \frac{2(d_2 +2)\beta}{n-d_2 -2}  - \frac{d_2+1}{n-d_2-2} -\frac{4}{d_4},
\; \frac{2(d_2 +2)\beta}{n-d_2 -2}  - \frac{d_2+1}{n-d_2-2},\cdots)$

{\noindent where} all terms from the fifth onwards are equal and where

\centerline{$\beta := 1 + \frac{c_1}{2} = \frac{8(n-d_2-2) + d_4 (n+d_2)}{2d_4 (n+d_2 +2)}$.}

The nullity condition for $\bar{c}$ implies $d_4 > 52$, so $n > 56$
and we deduce $0 < \beta < \frac{15}{26} < \frac{3}{5}$. Hence $\frac{a_1}{d_1} >0$.
It is also easy to show that $\frac{a_i}{d_i} >0$ for $i \geq 5$
and $\frac{a_3}{d_3} >0$.

So if $(-2^i, 1^j)$ gives an element of $\bar{a}^{\perp}$ we need $i=2$ or 4.
As $d_4 >52$, Lemmas \ref{no6} - \ref{case1} show that
 case (5) must hold with respect to $\bar{a}$. In particular $d_i =2$,
so we cannot have $i=4$. Hence $i=2$ and $d_2=2$. Now
Eq.(\ref{aj}) implies
 $\beta = \frac{2}{3}, \frac{2n-5}{3n-4}, \frac{3n-9}{4(n-2)} +
\frac{n-4}{d_4(n-2)}$
or $\frac{3n-9}{4(n-2)}$, depending on whether $j=1,3,4$ or $\geq 5$.
In all cases this contradicts the bound $\beta < \frac{3}{5}$ and $n > 56$.

(B) Consider $u = (-1,1,0,-1, 0, \cdots)$. We have $(a_i/d_i) =$

\centerline{$ (- \beta, \; \frac{2}{d_2} + 1-2 \beta,
 \; -\frac{d_2 + 1}{n-d_2 -2} - (\frac{n-6-3d_2}{n-d_2-2})\beta,
\; \frac{2(d_2 +2)\beta}{n-d_2 -2}  - \frac{d_2+1}{n-d_2-2} -\frac{2}{d_4},
\;  \frac{2(d_2 +2)\beta}{n-d_2 -2}  - \frac{d_2+1}{n-d_2-2}, \cdots)$}

{\noindent where all} terms from the fifth onwards are equal and

\centerline{$\beta := 1 + \frac{c_1}{2}= \frac{n+d_2-2}{2(n+d_2 +2)} + \frac{n-2}{d_2(n+d_2 +2)} +
\frac{n-d_2-2}{d_4(n+d_2 +2)}$.}

The nullity condition for $\bar{c}$ implies $d_2 \geq 9$ and $d_4 \geq 4$.
It is now easy to check that $\frac{3}{7} < \beta <  \frac{31}{36}$,
and that $\frac{a_i}{d_i} >0$ for $i \geq 5$.

As in (A), case (5) must hold with respect to $\bar{a}$. Now if
$(-2^i, 1^j)$ gives an element of $\bar{a}^{\perp}$ we need $d_i=2$.
This, combined with Lemma \ref{ineg}, means $i=1$ or $3$.

If $i=1$, Eq.(\ref{aj}) immediately shows $j$ cannot be 2. Moreover,
if $j=3$ or $\geq 5$, Eq.(\ref{aj}) yields a value for $\beta$ that violates
the nullity condition for $\bar{c}$. If $j=4$ we obtain
$\beta = \frac{n-1}{2n} + \frac{n-d_2-2}{d_4 n}$. As we are in case
(5) with respect to $\bar{a}$, we need to consider the elements in
$\bar{a}^{\perp} \cap \frac{1}{2}(d+ {\mathcal W})$ corresponding to
type II vectors. Their number and pattern, as stipulated by
Theorem \ref{class-1a}, together with orthogonality to $\bar{a}$,
imply further linear relations among the components of $(a_i/d_i)$
and small upper bounds for $r$ (usually of the form $r=5, 6$). In all
cases these additional constraints can be shown to be incompatible with
the above values of $\beta$.

As an illustration of the above method, note that our
type III vector is $(-2, 0,0,1,0,\cdots)$. If $\Delta^{\bar{a}}$ is
in case (5)(iii), Theorem \ref{class-1a} says that the possible type
IIs must have a $-1$ in place $1$ and a $0$ in place $4$. Since
$r \geq 5$, the remaining $-1$ must be in a place whose corresponding
dimension is $2$. As $d_3=2$ we can have $(-1, *, -1, 0, *, \cdots)$
where $*$ indicates a possible location of the $1$ in the type II.
The other possibility is for $-1$ to be in place $k$ for
some $k \geq 5$. After a permutation we can assume $k=5$, and $d_5=2$
must hold. The type II is then of the form $(-1, *, *, 0, -1, *, \cdots)$
where $*$ again indicates possible positions for the $1$. In the first
case, the orthogonality conditions imply $\frac{a_2}{d_2}=\frac{a_5}{d_5}$,
which gives $\beta=\frac{n-1}{2n} + \frac{n-d_2-2}{d_2 n}$.
Comparing with the value of $\beta$ from Eq.(\ref{aj}), we get $d_2=d_4$.
Using this in the first value of $\beta$ in (B) gives a contradiction
after some manipulation. In the second case, the argument we just
gave implies that we can only have $r=5$ and the orthogonality condition
implies $\frac{a_2}{d_2}=\frac{a_3}{d_3}$, which gives
$\beta=\frac{n-1}{n+d_2 +2} + \frac{2(n-d_2 -2)}{d_2(n+d_2 +2)}$.
After a short computation, one sees that the two values of $\beta$ are
again incompatible. If $\Delta^{\bar{a}}$ is in case (5)(ii), the argument
is essentially the same, as we only have to switch the places of the second
$-1$ and the $1$ in the type IIs.

Let us now take $i=3$. If $j=1$, Eq.(\ref{aj}) implies
$\beta = \frac{3d_2 +4 -n}{5d_2+10-n}$.
If the denominator is negative, then  $\beta >1$, which is
a contradiction. If it is positive we find
 that this is incompatible with the
inequality $ \beta > \frac{n+d_2-2}{2(n+d_2 +2)}$ which comes from
the displayed expression for $\beta$ above.

If $j=2$, we get $\beta=\frac{d_2 +1}{2(d_2 +2)} + \frac{n-d_2-2}{2d_2(d_2 +2)}$.
As above, we can rule this out by considering the vectors in
$\bar{a}^{\perp} \cap \frac{1}{2}(d+{\mathcal W})$ associated to type II
vectors. A similar argument works for $j \geq 5$, where we find $\beta =
\frac{n-2d_2-3}{2(n-2d_2-4)}$, and for $j=4$, where we have $\beta = \frac{n-d_2-2}
{d_4(n-2d_2-4)} + \frac{n-2d_2-3}{2(n-2d_2 -4)}$.

(C) Next let $u=(0,0,1,-1,-1, 0, \cdots)$. We have $(a_i/d_i) =
(1 - \beta, \; 1-2 \beta, \; \frac{n-2d_2 -3}{n-d_2 -2} -(\frac{n-6-3d_2}{n-d_2-2}) \beta,
\; \frac{2(d_2 +2)}{n-d_2 -2} \beta - \frac{d_2+1}{n-d_2-2}-\frac{2}{d_4},
\; \frac{2(d_2 +2)}{n-d_2 -2} \beta - \frac{d_2+1}{n-d_2-2}-\frac{2}{d_5},
\; \frac{2(d_2 +2)}{n-d_2 -2} \beta - \frac{d_2+1}{n-d_2-2}, \cdots)$
where all terms from the sixth on are equal and

\centerline{$\beta := 1 + \frac{c_1}{2}= \frac{n+d_2}{2(n+d_2 +2)} +
\left(\frac{1}{d_4} + \frac{1}{d_5} \right) \frac{n-d_2 - 2}{(n+d_2 +2)}$.}
The nullity condition for $\bar{c}$ implies $d_4$ and  $d_5 \geq 13$.
It is now easy to check that $0 < \beta <  \frac{2}{3}$, so $\frac{a_1}{d_1} >0$.
As in (A), we find also that  $\frac{a_i}{d_i}>0$ for $i \geq 6$, and that
$\frac{a_3}{d_3} >0$.

As in (A) again, case (5) must hold with respect to $\bar{a}$, so if
$(-2^i, 1^j)$ gives an element of $\bar{a}^{\perp}$ we need $d_i=2$.
This, combined with Lemma \ref{ineg}, means $i=2$. In this situation
in all cases Eq.(\ref{aj}) gives a value of $\beta$ incompatible with
our bounds on $n$ and  $\beta$.
\end{proof}

\begin{lemma} \label{no5b}
Case $($5$)$$($iii$)$ cannot arise if $r \geq 5$.
\end{lemma}

\begin{proof} This is similar to the proof of the previous Lemma so we will
be brief.

(A) Consider $u = (0,-2,0,1,0,\cdots)$. Now $(a_i/d_i)$ is given by

\centerline{$ (1-\beta,  -\frac{4}{d_2} + 1-2 \beta,
\; \frac{(n-3)+(n+d_2-2)(\beta-1)}{n-d_2-2},
\; \frac{2}{d_4} + \frac{2d_2 \beta -(d_2+1)}{n-d_2-2},
\; \frac{2d_2 \beta -(d_2+1)}{n-d_2-2}, \cdots,
 \frac{2d_2 \beta -(d_2+1)}{n-d_2-2})$}

\noindent{where}

\centerline{$\beta := 1+ \frac{c_1}{2} =\frac{1}{2} - \frac{(d_2 +4d_4)(n-2-d_2-d_4) +4d_4^2}
{2d_2 d_4 ( 2n -d_2-4)}$.}

The nullity condition for $\bar{c}$ implies $d_4 \geq 8, d_2 \geq 27$ and
$d_2 > 2d_4$,  and it readily follows that
$\frac{1}{4} < \beta < \frac{1}{2}$. In particular, $\frac{a_1}{d_1}>0$.

As before, we see that case (5) holds with respect to $\bar{a}$, so for
$\bar{a}^{\perp}$ we have to consider type III vectors $(-2^i, 1^j)$
where $d_i = 2$. So we need only consider $i=3$ or $i \geq 5$.

In either situation, we proceed as in part (B) of the proof of
Lemma \ref{no5a}, and obtain inconsistencies in the equations
involving $\beta$ or contradiction to the bounds on $\beta$
or the dimensions.

(B) Let $u=(0,0,-1,1,-1,0,\cdots).$ The nullity condition on
$\bar{c}$, which has a symmetry in $d_4$ and $d_5$, now implies
$d_4, d_5 \geq 46$ and $> 28d_2$. Now $(a_i/d_i)$ is given by

\centerline{$(1-\beta, \; 1-2\beta, \; (\frac{n+d_2 -2}{n-d_2-2})\beta -\frac{n-1}{n-d_2 -2},
\; \frac{2}{d_4} + \frac{2d_2 \beta -(d_2 +1)}{n-d_2-2},
\; \frac{-2}{d_5} + \frac{2d_2 \beta -(d_2 +1)}{n-d_2 -2},
\; \frac{2d_2 \beta -(d_2 +1)}{n-d_2 -2}, \cdots)$}

{\noindent where} all terms from the sixth on are equal and

\centerline{$\beta:=1 + \frac{c_1}{2}=\frac{1}{2}+ \frac{1}{n+d_2 -2} +
     (\frac{1}{d_4}+\frac{1}{d_5})(\frac{n-d_2 -2}{n+d_2 -2}).$}

We deduce that $\frac{1}{2} < \beta < \frac{3}{5}$ and so
$\frac{1}{2} > \frac{a_1}{d_1} > 0$.
$\Delta^{\bar{a}}$ must be in case (5), and for $(-2^i, 1^j)$ associated
to elements of $\bar{a}^{\perp} \cap \frac{1}{2}(d+{\mathcal W})$, we have
$d_i=2$ and we need only consider $i=2, 3$ or $\geq 6$.

If $i=2$, Eq.(\ref{aj}) and the upper bound on $\beta$ imply
$\frac{a_j}{d_j} > \frac{3}{5}$. As $d_2 =2$, one checks that
this never holds.

If $i \geq 6$, Eq.(\ref{aj}) implies $\frac{a_j}{d_j} > \frac{47}{48}$.
The bound $d_4, d_5 > 28d_2$ can be used to show that this never happens.

If $i=3$, Eq.(\ref{aj}) and the expression for $\beta$ above are seen
to be incompatible if we use the bounds on $d_4, d_5$ and $\beta$.

(C) Consider $u=(-1,-1,0,1,0,\cdots).$ The nullity condition on
$\bar{c}$ gives $d_4 \geq 4, d_2 \geq 5$. The vector $(a_i/d_i)$ is

\centerline{$(-\beta, 1-2\beta-\frac{2}{d_2},
\; \frac{(n+d_2 -2)\beta -(d_2 +1)}{n-d_2 -2},
\; \frac{2}{d_4} + \frac{2d_2 \beta -(d_2 +1)}{n-d_2 -2},
\; \frac{2d_2 \beta - (d_2 +1)}{n-d_2 -2}, \cdots)$}

\noindent{where} all terms from the fifth on are equal and

\centerline{$\beta:=1 + \frac{c_1}{2} =
\frac{1}{2}- \frac{d_4^2 +(d_2 +d_4)(n-d_2-d_4-2)}{d_2 d_4 (3n-d_2-6)}$.}
Now as $\frac{1}{d_2} + \frac{1}{d_4} < \frac{1}{2}$, we see that $\frac{1}{3} <
\beta < \frac{1}{2}$.

Again $\Delta^{\bar{a}}$ is in case (5) and we consider vectors
$(-2^i, 1^j)$ associated to elements of $\bar{a}^{\perp} \cap \frac{1}{2}(d+ {\mathcal W})$,
where we must have $d_i=2$, so $i \neq 2, 4$.

If $i=1$, Eq.(\ref{aj}) becomes $\frac{a_j}{d_j} = 1-2 \beta > 0$. This immediately
means $j \neq 2, 5, \cdots, r$. If $j=3$, the value of $\beta$ from Eq.(\ref{aj})
and the above expression for $\beta$ lead to $d_4 \leq 10/3$. For $j=4$, we
obtain a contradiction by the method of part (B) in the proof of Lemma \ref{no5a}.

If $i \geq 5$, Eq.(\ref{aj}) says $\frac{a_j}{d_j} = \frac{4d_2 \beta+(n-3d_2-4)}{n-d_2-2}$.
If $i=3$, Eq.(\ref{aj}) say $\frac{a_j}{d_j} = 1 -\frac{2(d_2+1)}{n-d_2-2} +
  \frac{2(n+d_2-2)\beta}{n-d_2-2}$.
In both situations, we can again apply the method of part (B) in
the proof of Lemma \ref{no5a} to obtain contradictions.
\end{proof}

\smallskip

The last case to consider is case (5)(ii) with $r=4$, which is the same
as case (5)(iii) with $r=4$ if we interchange the third and fourth summands.

\begin{lemma} \label{no5dim4}
No configurations for case $($5$)$$($ii$)$ with $r=4$ can occur.
\end{lemma}

\begin{proof}
When $r=4$ we no longer have $d_3=2$, but the nullity condition for
$\bar{c}$ implies that $\frac{1}{d_3}+\frac{1}{d_4} \geq \frac{1}{2}$.
Hence either $\{d_3, d_4\}$ is one of $\{3, 3\}, \{3, 4\}, \{3, 5\},
\{3, 6\}, \{4, 4\}$ or one of $d_3$ or $d_4$ is $2$. Using this
together with the nullity conditions for $\bar{a}, \bar{c}$ and
the orthogonality conditions, we see that we only need to consider
$u=(0, -2, 1, 0), (0, 0, 1, -2), (-1, 1, 0, -1)$ and $(-1, -1, 1, 0)$.

(A) Let $u=(0, -2, 1, 0)$. From the nullity condition for $\bar{c}$
we deduce that $d_4=2$, $d_2 \geq 13$, and $d_3 \geq 3$. Now

\centerline{$(a_i/d_i) = (1-\beta, \; 2\beta -1 -\frac{4}{d_2},
\; \frac{2}{d_3} + \frac{2d_2 \beta-(d_2 +1)}{d_3 + 2},
\; \frac{(2d_2 +d_3 +2)\beta-(d_2 +1)}{n-d_2 -2})$}

{\noindent where} $\beta:=1+\frac{c_1}{2}=\frac{1}{2}-
\frac{d_2 + 4d_3 +2d_3^2}{d_2 d_3 (d_2 + 2d_3 + 4)}$. We find
that $\frac{11}{26} < \beta < \frac{1}{2}$ and so
$\frac{a_1}{d_1} > 0$.

The above facts imply that $\Delta^{\bar{a}}$ is again in case (5)(ii),
and for $(-2^i, 1^j)$ associated to an element of $\bar{a}^{\perp}$
we must have $i=4$. If $j=1,2$, Eq.(\ref{aj}) leads to a contradiction
to the above dimension restrictions. The case $j=3$ can be eliminated using
the method of part (B) in the proof of Lemma \ref{no5a}.

(B) Next let $u=(0,0,1,-2)$. The nullity condition for $\bar{c}$
implies that $d_3=2$ and $d_4 > 32d_2 + 14 \geq 46$. Now $(a_i/d_i)$
is given by

\centerline{$(1-\beta, \; 1-2\beta,
\; \frac{d_4 -d_2 +1}{d_4 +2} + (\frac{2d_2 + 2-d_4}{d_4 +2})\beta,
\; -\frac{4}{d_4} + \frac{2(d_2 +2)\beta-(d_2 +1)}{d_4 +2})$}

{\noindent where} $\beta:=1+ \frac{c_1}{2} = 1- \frac{d_4^2 + 2d_2 d_4 -16}{2d_4(2d_2+d_4 +6)}$.
One easily sees that $\frac{1}{2} < \beta < 1$, so that
$0< \frac{a_1}{d_1} < \frac{1}{2}$. Since $d_4 > 2d_2 +2$ we obtain
$0< \frac{a_3}{d_3} < \frac{3}{5}.$

Therefore, $\Delta^{\bar{a}}$ is in case (5)(ii) and for $(-2^i, 1^j)$
associated to an element of $\bar{a}^{\perp}$ we must have (by Lemma \ref{ineg})
$i=2$ and so $d_2=2$. Putting this value of $d_2$ into the nullity condition
for $\bar{c}$ gives a cubic equation in $d_4$ with no integral roots,
a contradiction.

(C) Consider now $u=(-1, 1, 0, -1)$. From the nullity condition for $\bar{c}$
we deduce that $d_2 \geq 4, d_3 =2,$  $d_4 \geq 3$ and $4d_2 > 3d_4$. Also,

\centerline{$(a_i/d_i)=(-\beta, \frac{2}{d_2}+1-2\beta,
  \; \frac{(2d_2 +2-d_4)\beta-(d_2 +1)}{d_4 +2},
  \; -\frac{2}{d_4} + \frac{2(d_2 +2)\beta -(d_2 +1)}{d_4 +2})$}

{\noindent where} $\beta:=1+\frac{c_1}{2} =
   \frac{1}{2}+\frac{2d_2 +2d_4 +d_4^2}{d_2 d_4(2d_2 + d_4 +6)}$.
It follows that $\frac{1}{2} < \beta < \frac{3}{4}$ and $\frac{a_2}{d_2} > 0$.

Now we see that $\Delta^{\bar{a}}$ is either in case (1) or (4) or (5)(ii).
In the first two instances, by Lemmas \ref{case1}, \ref{case4} $d$ is a
permutation of $(4, 2, 2, 9)$. Since $3d_4 < 4d_2$ we have $d_2 = 9, d_4 =4$.
But then the null condition for $\bar{c}$ is violated. So we are in case (5)(ii).
For $(-2^i, 1^j)$ associated to an element of $\bar{a}^{\perp}$, as $d_i=2$,
we have $i=1$ or $3$.

If $i=1$, then Eq.(\ref{aj}) becomes $\frac{a_j}{d_j} = 1-2\beta < 0$, so $j=3, 4$.
When $j=3$ the value of $\beta$ given above together with Eq.(\ref{aj}) imply
that $d_4 = 3$ or $4$. But then the null condition for $\bar{c}$ is violated.
For $j=4$ we may use the argument of part (B) of the proof of Lemma \ref{no5a}.

If $i=3$, using $\beta> \frac{1}{2}$ in Eq.(\ref{aj}), we see that $j=2, 4$.
In either case, applying our bounds for the dimensions in Eq.(\ref{aj}) lead
to contradictions.

(D) Let $u=(-1,-1,1,0)$. The null condition for $\bar{c}$ implies that
$d_4 =2, d_3 \geq 3, d_2 \geq 5$. With $\beta:=1+ \frac{c_1}{2}$, we have

\centerline{$(a_i/d_i)= (-\beta, 1-2\beta-\frac{2}{d_2},
\; \frac{2}{d_3} + \frac{2d_2 \beta -(d_2 +1)}{d_3 +2},
\; \frac{(2d_2 + d_3  +2)\beta-(d_2 +1)}{d_3 +2})$.}

One computes that $\beta=\frac{1}{2}-\frac{2d_2 + d_3^2 +2d_3}{d_2(3d_3^2 +2d_2 d_3 + 6d_3)}$,
and from the dimension bounds one gets $\frac{5}{12} \leq \beta < \frac{1}{2}$.
$\Delta^{\bar{a}}$ cannot be in case (1) or (4), otherwise as $d_2 \geq 5$, we
must have $d_2 =9, d_3 =4$, and the null condition for $\bar{c}$ is violated.
So $\Delta^{\bar{a}}$ is in case (5)(ii). For $(-2^i, 1^j)$ associated to an
element of $\bar{a}^{\perp}$, we must then have $i=1, 4$.

If $i=1$, then Eq.(\ref{aj}) is $\frac{a_j}{d_j} = 1-2\beta > 0$, so $j=3$ or
$4$. In either situation, we may apply the argument of part (B) of the
proof of Lemma \ref{no5a} to get a contradiction. If $i=4$, Eq.(\ref{aj})
together with the dimension bounds above show first that we can only
have $j=3$. In that case, a more detailed look at Eq.(\ref{aj}) leads to
a contradiction.
\end{proof}

\medskip
We can summarise our discussions thus far by

\begin{theorem} \label{unique1B} Let $r \geq 4$ and $K$ be connected.
Suppose that we are not in the situation of Theorem \ref{Fano}.
Assume that $\bar{c} \in {\mathcal C}$ is a null vector such that
$\Delta^{\bar{c}}$ has the property that there is a unique vertex of
type $($1B$)$ and all other vertices are of type $($1A$)$. Then the
only possibilities are given by Lemmas \ref{case4}/\ref{case1} and \ref{case3},
up to interchanging $\bar{a}$ and $\bar{c}$ and a permutation of the
irreducible summands.  \;   $\qed$
\end{theorem}

We will now sharpen the above Theorem using Proposition \ref{edgeCC}.

\begin{corollary} \label{finish}
Let $r \geq 4$. Assume that $K$ is connected and we are not in the
situation of Theorem \ref{Fano}. Then the possibilities given
by Lemmas \ref{case4} and \ref{case3} $($and hence Lemma \ref{case1}$)$
cannot occur.
\end{corollary}

\begin{proof}
We will discuss the $r=4$ case (i.e. that in Lemmas \ref{case4},
 \ref{case1}) in detail and leave the details of
the $r=5$ case (from Lemma \ref{case3})
to the reader, as the arguments are very similar.

In the $r=4$ case, first observe that $\mathcal C$
has exactly two null vectors, $\bar{c}$ and $\bar{a}$ in the notation of
 Lemma \ref{case4}, as the entries of $a,c$ are determined by the vector
$d$ of dimensions. Hence, by Prop. \ref{sep-hypl}, these are the
only elements of $\mathcal C$ outside conv$(\frac{1}{2}(d + {\mathcal W}))$.

Since $(-1^4)$ is a vertex of $\mathcal W$, all type II vectors
in $\mathcal W$ must be zero in place $4$. As $(1, -1,-1, 0)$
is associated to an element of $\bar{c}^{\perp}$, it, together
with $(-1, 1, -1, 0), (-1, -1, 1, 0)$ are the only type II
vectors in $\mathcal W$.

Next we analyse vectors in $\mathcal W$ and see if
they are associated to elements of $\mathcal C$, this last property
being important for applying Prop \ref{edgeCC}. Recall that
$(1,-2,0,0),(1,0,-2,0),(-2,1,0,0),(-2,0,1,0)$ must be in $\mathcal W$.
The first two give elements of $\bar{c}^{\perp}$, the last two give
elements of $\bar{a}^{\perp}$.

First consider $v=(1, -2, 0, 0)$. Now
 $\bar{v}$ is  a vertex of conv$(\frac{1}{2}(d+{\mathcal W}))$.
By the superpotential equation, $d+v = 2\bar{v}$ can be written
as $\bar{c}^{(\alpha)} + \bar{c}^{(\beta)}$,
with $\bar{c}^{(\alpha)}, \bar{c}^{(\beta)} \in {\mathcal C}$.
Since $\bar{v}$ is a vertex, every such expression must involve
$\bar{a}$ or $\bar{c}$, unless it is the trivial expression
$\bar{v} + \bar{v}$ and $\bar{v} \in {\mathcal C}.$  By computing
$2v-a, 2v-c$ we find that these cannot lie in conv$(\mathcal W)$
(it is enough to exhibit one component $<-2$ or $>1$).
Thus $\bar{v} \in {\mathcal C}$. Now an analogous argument shows
that if $w=(0, -2, 1, 0)$ is in $\mathcal W$
then $\bar{w}$ also lies in $\mathcal C$. But $vw$ is an edge of ${\rm conv}(\mathcal W)$
with no interior points in $\mathcal W$. So Prop \ref{edgeCC}
gives $1 - \frac{4}{d_2} =4 J(\bar{v}, \bar{w}) = 0$, a contradiction
to $d_2=2$. Hence $w \notin {\mathcal W}$. Similarly we see
$(0,-2,0,1) \notin \mathcal W$.

Next consider $z=(-1, -1, 1, 0)$. By Remark \ref{facts}(e) and the above, $z$
is a vertex of ${\rm conv}(\mathcal W)$. As above we can show that
$\bar{z} \in {\mathcal C}$. Now $v, z$ are the only elements of
the face $\{x_1 + 2x_2 = -3\} \cap {\rm conv}(\mathcal W)$
(cf. proof of Prop. 4.3 in \cite{DW4}). So applying
Prop \ref{edgeCC} to $vz$ we obtain $0 = 4J(\bar{v}, \bar{z}) = \frac{1}{4}$,
 a contradiction.

To handle the $r=5$ case (from Lemma \ref{case3}), first observe that null
elements of $\mathcal C$ must have entries $\frac{2}{3}$ in two
of the places $1,\cdots,4$ and $-\frac{2}{3}$ in the other two
places. For $a,c$ as in Lemma \ref{case3}, we can take
 $(1,-2,0,0,0), (1,0,-2,0,0), (0,1,0,-2,0), (-2,1,0,0,0)$ to lie
in $\mathcal W$.
The argument above to show that $\bar{v}$ is in $\mathcal C$
still works for such type III vectors $v$. As above, we can use
Prop \ref{edgeCC} to show the other type III vectors
$(-2^i, 1^j) \; (i \leq 4)$
do not lie in $\mathcal W$; hence the $\bar{a},\bar{c}$ of Lemma \ref{case3}
are the only null elements of $\mathcal C$.

Let $z=(-1, 1, -1, 0, 0)$
(it lies in $\mathcal W$ since $(1, -1, -1, 0, 0)$ is associated to
an element of $\bar{c}^{\perp}$) and $v=(1, 0, -2, 0, 0)$.
As above we find $\bar{z}$ is in $\mathcal C$, and the arguments
of Prop 4.3 in \cite{DW4} show $vz$ is an edge of conv$(\mathcal W)$.
A contradiction results as above from applying Prop \ref{edgeCC} to
$vz$.
\end{proof}

\medskip

The discussion at the beginning of this section now tells us that
if $K$ is connected and $r \geq 4$ the only case when we have a superpotential
of the kind under discussion is that of Theorem \ref{Fano}.
The proof of Theorem \ref{classthm} is now complete.

\medskip

\noindent{\bf Concluding remarks.} {}

\vspace{0.2cm}

1. When $r=2$, then $c$ is collinear with the elements of $\mathcal W$.
In other words, the projected polytope $\Delta^{\bar{c}}$ reduces
to a single vertex, which must be of type (2). The possible elements of
$\mathcal W$ are $(-2,1),(-1,0),(0,-1),(1,-2)$. If  $\mathcal W$ has just
two elements then Theorem \ref{no2} tells us we are either in the situation
of Theorem \ref{Fano} (the B{\'e}rard Bergery examples), or in Example 8.2 or
the third case of Example 8.3 in \cite{DW4}. In fact one can show that
this last possibility can be realised in the class of homogeneous hypersurfaces
exactly when $(d_1, d_2) = (8, 18)$. An example for these dimensions
is provided by $G = SU(2)^9 \ltimes {\rm Sym}(9)$
(where Sym$(9)$ acts on $SU(2)^9$ by permutation) and $K$ is the
product of the diagonal $U(1)$ in $SU(2)^9$ with Sym$(9)$. The arguments
of \cite{DW2} show that this in fact gives an example where the
cohomogeneity one Ricci-flat equations are fully integrable.

If  $\mathcal W$ has three elements, we may adapt the proof
of Theorem \ref{cedge} to derive a contradiction. Here the essential point
is that whenever we had to check that a sum of two elements of $\mathcal C$
does not lie in $d+{\mathcal W}$, such a fact remains true because the
interior point of $vw$ is the midpoint.

If $\mathcal W$ contains all four possible elements, then
$\kf \subset \g$ is a maximal subalgebra (with respect to inclusion).
We suspect that this case also does not occur. In any event, it is
of less interest because the only way to obtain a complete cohomogeneity
one example is by adding a $\Z/2$-quotient of the principal orbit
as special orbit.

2. The only parts of this paper which depend on $K$ being connected
(or slightly more generally, on the condition in Remark \ref{disconnK})
are parts of \S 5,  Case (ii) of  \S 9, and all of \S 10. To remove this
condition, the main task would be generalizing Theorem \ref{class-1a} by
getting a better handle on the type II vectors associated to (1A)
vertices (cf Lemmas (\ref{II/II})-(\ref{II/III})).

\end{document}